\documentclass[final,3p,times]{article}
\usepackage{fullpage}
\usepackage[latin9]{inputenc}
\usepackage{array}
\usepackage{verbatim}
\usepackage{bm}
\usepackage{multirow}
\usepackage{amsmath}
\usepackage{amssymb}
\usepackage{graphicx}
\usepackage{esint}
\usepackage{amsfonts}
\usepackage{booktabs}
\usepackage{siunitx}
\newtheorem{lemma}{Lemma}[section]

\usepackage{tikz}
\usetikzlibrary{shapes.geometric, arrows, positioning, matrix}

\makeatletter

\usepackage{tabularx}\usepackage{subfigure}\usepackage{subfigure}\usepackage{color}

\usepackage{setspace}

\usepackage{bm}%% The amsthm package provides extended theorem environments
\begin{document}
%\title{A fully three-dimensional multi-moment algorithm for Euler equation on
%hybrid unstructured grids}

\title{THINC scaling method that bridges VOF and level set schemes}

\author{Ronit Kumar$^{(\star)}$, Lidong Cheng$^{(\star,\dagger)}$, Yunong Xiong$^{(\star)}$, Bin Xie $^{(\ddagger)}$, \\
R\'emi Abgrall$^{(\sharp)}$ and Feng Xiao$^{(\ddagger)}$\footnote{Corresponding author}\\
$(\star)$ School of Engineering, Department of Mechanical Engineering, \\Tokyo Institute of Technology,
Tokyo, 152-8550, Japan\\
$(\dagger)$ School of Aerospace Engineering, Beijing Institute of Technology,\\
Beijing, 100081, China\\
{$(\ddagger)$}School of Naval Architecture, Department of Ocean and Civil Engineering, \\Shanghai Jiaotong University, 
Shanghai, 200240, China\\
$(\sharp)$ Institute of Mathematics, University of Z\"urich, Z\"urich, \\CH8057, Switzerland}
\date{}

%\cortext[cor]{Corresponding authors: Dr. F. Xiao (Email: xiao.f.aa@m.titech.ac.jp)}
\maketitle
\begin{abstract}
\indent We present a novel interface-capturing scheme, THINC-scaling, to unify the VOF (volume of fluid) and the level set methods, which have been developed as two different approaches widely used in various applications. The key to success is to maintain a high-quality THINC reconstruction function using the level set field to accurately retrieve geometrical information and the VOF field to fulfill numerical conservativeness. The interface is well defined as a surface in form of a high-order polynomial, so-called the polynomial surface of interface
 (PSI). The THINC reconstruction function is then used to update the VOF field via a finite volume method, and the level set field via a semi-Lagrangian method. Seeing the VOF field and the level set field as two different aspects  of the THINC reconstruction function, the THINC-scaling scheme preserves at the same time the advantages of both VOF and level set methods, i.e. the mass/volume conservation of the VOF method and the geometrical faithfulness of the level set method, through a  straightforward solution procedure. The THINC-scaling scheme allows to represent an interface with high-order polynomials and has algorithmic simplicity which largely eases its implementation in unstructured grids. Two and three dimensional algorithms in both structured and unstructured grids have been developed and verified. The numerical results reveal that the THINC-scaling scheme, as an interface capturing method,  is able to provide high-fidelity solution comparable to other most advanced methods, and more profoundly it can resolve sub-grid filament structures if the interface is represented by a polynomial higher than second order.   
\end{abstract}

\bigskip
%\begin{keyword}
\noindent \textbf{Keywords:} Moving interface, multiphase flow,  VOF,  THINC,  level set,  high-order interface representation, 
 unstructured grids.
%\end{keyword}

\section{Introduction}

Representing and updating interfaces that move in 3D space poses a challenging task to scientific and engineering computing. Being the mainstream approaches, the VOF (volume of fluid) and level set are the two most popular methods used in capturing moving interface and find their applications in diverse fields, such as numerical simulation of multiphase fluid dynamics, graphic processing, topological optimization, and many others. Historically, these two approaches were independently developed based on completely different concepts and solution methodologies, and thus have their own superiority and weakness.  

VOF method uses the volume fraction of one out of multiple fluid species in a control volume (mesh cell) to describe the distribution of the targeted fluid in space. Rigorous numerical conservation can be guaranteed if a finite volume method is used to transport the VOF function, which is found to be crucial in many applications. 

The VOF function by definition has a value between 0 and 1 for each mesh cell, and the interface can be identified as surface (3D) or line (2D) segments cutting through mesh cells with given VOF values. There is a freedom in choosing  geometrical information to represent the interface within a mesh cell. For example, the simplest version of this sort is the SLIC (Simple Line Interface Calculation) method that uses line segments aligned
to grid lines \cite{slic1976,hirt1981volume}. An improved and more accurate way for interface reconstruction is to add the normal direction or the orientation of the interface as another geometrical information, which results in a large family of the VOF method  known as the PLIC (Piecewise Linear Interface Calculation) schemes, the most representative geometrical VOF method. The PLIC VOF methods involve geometrical computations, and use a plane to represent the interface segment embedded in a mesh cell. Efforts have been devoted in the past decades to devise  efficient and accurate numerical formulations for geometrical reconstructions in PLIC schemes \cite{youngs1982time,lafaurie1994modelling,rider1998reconstructing,scardovelli2000analytical,pilliod2004second}, which make the PLIC VOF method mature in structured grids with  adequate numerical accuracy for many applications. The analytical relation developed in \cite{scardovelli2000analytical} provides a very efficient way for PLIC algorithm in Cartesian grid.  However, the geometrical reconstruction becomes quite challenging in the case of high order surface representation, rather than the currently used plane representation which is in the form a first-order polynomial. There are only a few geometrical reconstructions using quadratic interface reconstruction reported in 2D \cite{renardy2002,scardovelli2003interface,lopez2004volume,diwakar2009quadratic}. The algorithmic complexity of geometrical VOF method also increases when applied to unstructured grids. Another class of ease-to-use VOF schemes, so-called algebraic VOF, have been devised and found particular popularity in  unstructured grids, like those in \cite{rudman97,ubbink99,darwish06,heyns13,zhang14}
among others.  In spite of simplicity,  it is usually observed that  the  algebraic VOF schemes are less appealing in numerical accuracy compared to the geometrical VOF schemes using PLIC reconstructions.  
 
Moreover, the VOF function is usually characterized by a large jump or steep gradient across the interface. So, directly using the VOF field to retrieve the geometrical information of the interface, such as normal and curvature, may result in larger errors in compared to the level set method. 

The level set method \cite{osher1988fronts,sethian1999level,osher2003implicit}, on the other hand, defines the field function as a signed distance function (or level set function) to the interface, which possesses a uniform gradient over the whole computational domain and thus provides a perfect field function to retrieve the geometrical properties of an interface. However, the level set function is not conceptually nor algorithmically conservative. The numerical solution procedure, including both transport and reinitialization, does not guarantee the conservativeness of numerical solution. It may become a fatal problem in many applications, like multiphase flows involving bubbles or droplets.  In comparison with the VOF method, level set method is less popular in real-case multiphase flow simulations. Efforts have been made to enhance the numerical conservation of level set method, such as the mass/volume conservation correction with global or local  constraints \cite{sussman1999level,enright2002hybrid,Manuel2019CLS}, or the so-called conservative level set schemes which are in spirit equivalent to the phase field method \cite{olsson2005CLS,sun2007PF,Chiu2011PF}.             

Another natural practice to retain the pros while overcome the cons of the two methods is to combine the VOF method and level set method, which leads to the coupled level set/VOF methods (CLSVOF) \cite{sussman2000coupled,menard2007coupling, yang2006adaptive, sun2010coupled}. The PLIC  type VOF method is blended with the level set method so as to improve both conservativeness and geometrical faithfulness in numerical solution \cite{aniszewski2014VOF}. In a CLSVOF scheme of this type, some extra algorithmic efforts are required for adjustment between the VOF and level set fields, in addition to the standard computing operations of the PLIC and level set schemes. 

In this paper, we propose a new scheme that unifies the VOF and level set methods, based on the observation that the VOF field can be seen as a scaled level set field using the THINC (Tangent of Hyperbola Interface Capturing) function, which has been used in a class of schemes for capturing moving interface \cite{xiao2005simple,xiao2011revisit,ii2012interface,ii2014interface,xie2014efficient,xie2017toward,qian2018}. The resulting scheme, so-called THINC-scaling scheme, converts the field function from level set to VOF by the THINC function, and converts the VOF field back to the corresponding level set field via an inverse THINC function. So, the VOF field and the level set field can be seen as two different aspects of the THINC function. The THINC-scaling method can make use of the advantages of both VOF and level set at different stages of a single solution procedure, which eventually realizes the high-fidelity computation of moving interface regarding both numerical conservativeness and geometrical representation. Without explicit geometrical computation, the THINC-scaling method is algorithmic simple and can be straightforwardly extended to 3D unstructured grids. More profoundly, the THINC-scaling method is able to use high-order polynomials to represent the interface without substantial difficulties. 

This paper is organized as follows. Section 2 describes some basic formulations that connect the VOF field and level set field through the THINC function. The THINC-scaling scheme is detailed in section 3. We present numerical results of  benchmark tests on both structured and unstructured grids in two and three dimensions to verify the THINC-scaling method in section 4, and end the paper with some summary remarks in section 5.   

\section{The connection between level set and VOF functions}

We consider an interface $\partial \Omega$ separating two kinds of fluids, fluid 1 and fluid 2, occupying volumes $ \Omega^1$ and $ \Omega^2$ respectively in space. We introduce the following two indicator functions to identify the different fluids and the interface. 
\begin{figure}[htbp]
	\centering
	\subfigure[] {
		\centering
		\includegraphics[width=0.4\textwidth]{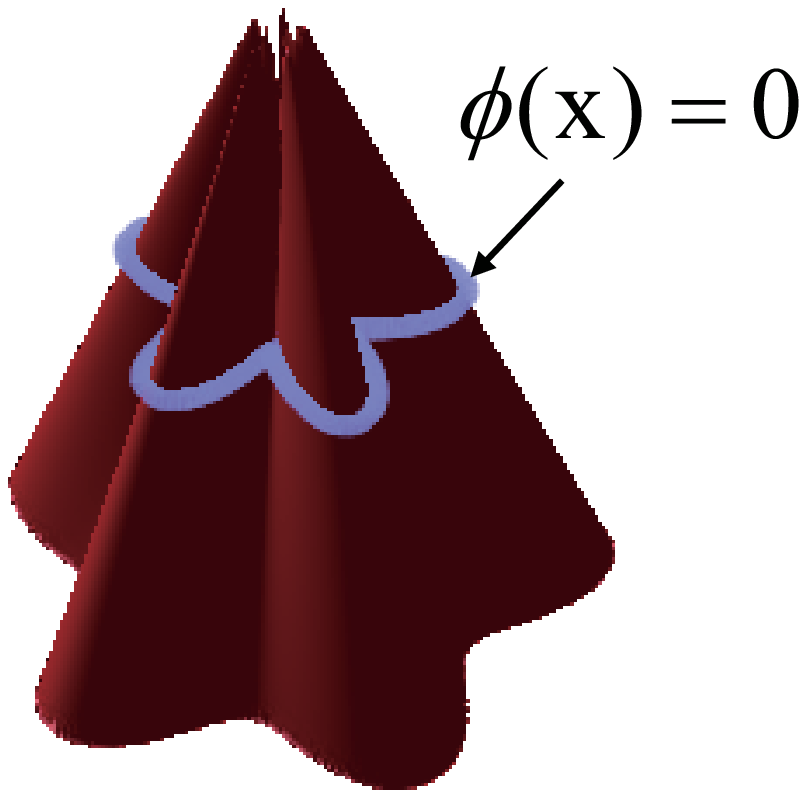}} \hspace{0.5cm}
	\subfigure[] {
		\centering
		\includegraphics[width=0.4\textwidth]{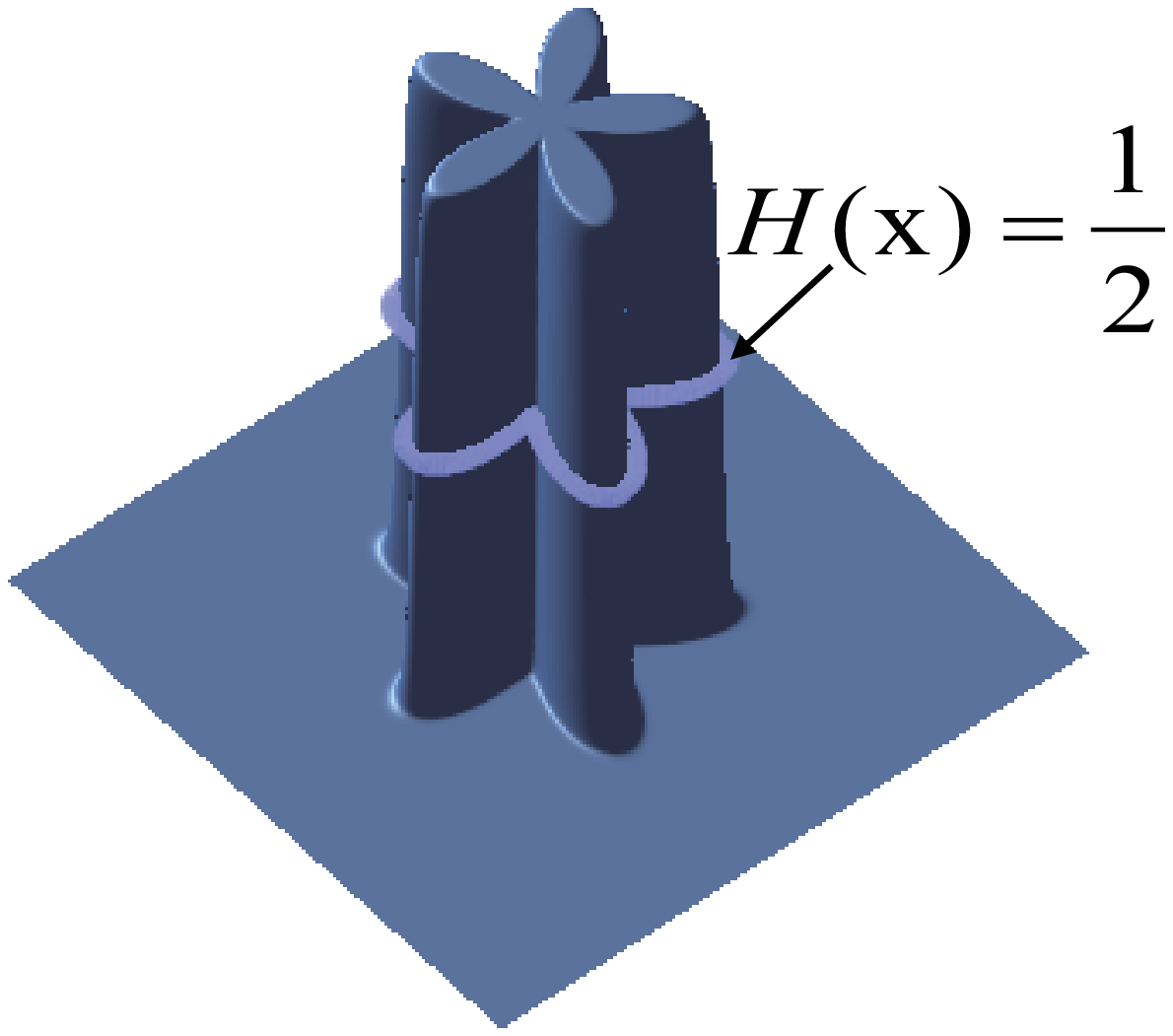}}
	\caption{Two indicator functions to identify multi-materials and interface: (a) The level set function; (b) The THINC function (a VOF function).}
\label{fun}
\end{figure}

\begin{itemize}
\item Level set function (Fig.\ref{fun}(a)): \\
The level set function is defined as a signed  distance from a point $\textbf{x}=(x,y,z)$ in three dimensions to the interface $\partial  \Omega$ by
\begin{equation}
\phi(\textbf{x})=\left\{
\begin{array}{rcc}
\displaystyle{\inf_{\textbf{x}_I\in \partial  \Omega}\|\textbf{x}-\textbf{x}_I\|} & & \text{if } \textbf{x}\in  \Omega^1 \\
0\quad\quad\quad & & \text{if } \textbf{x}\in\partial  \Omega \\
-\displaystyle{\inf_{\textbf{x}_I\in \partial  \Omega}\|\textbf{x}-\textbf{x}_I\|} & & \text{if } \textbf{x}\in  \Omega^2
\end{array}, \right.
\label{ls}
\end{equation}
where $\textbf{x}_I$ represents any point on the interface, and referred to as interface point. 

\item VOF function (Fig.\ref{fun}(b)): \\
The VOF function in the limit of an infinitesimal control volume is the Heaviside function. Assuming the VOF function as the abundance of fluid 1, we have the Heaviside function in its canonical form as
\begin{equation}
H^c(\textbf{x})=\left\{\begin{array}{ll}{1} & \text{if } \textbf{x}\in  \Omega^1 \\  
{\frac{1}{2}} & \text{if } \textbf{x}\in\partial  \Omega \\ 
{0} & \text{if } \textbf{x}\in  \Omega^2 \end{array}. \right.
\label{heaviside}
\end{equation}
\end{itemize}

Recall that 
\begin{equation}
H^c(\textbf{x})=\frac{1}{2} \lim_{\beta \rightarrow \infty}\left(1+\tanh \left(\beta\textbf{x})\right)\right), 
\label{hv-def}
\end{equation}
we get a continuous Heaviside function 
\begin{equation}
H(\textbf{x})=\frac{1}{2} \left(1+\tanh \left(\beta\textbf{x})\right)\right) 
\label{hv-beta}
\end{equation}
with a finite steepness parameter $\beta$. 

Given a computational mesh with cells of finite size, the cell-wise VOF function is defined by 
\begin{equation}
\bar{H}_{\Omega_i}=\frac{1}{|\Omega_i|}\int_{\Omega_i}H(\textbf{x})d\Omega,
\label{vof-def}
\end{equation}
where $\Omega_i$ is the target cell with a volume $|\Omega_i|$. In practice, we recognize the interface cell where an interface cuts through, in terms of the VOF value, by $\varepsilon \le \bar{H}_{\Omega_i} \le 1-\varepsilon$ with $\varepsilon$ being a small positive, e.g. $\varepsilon={10}^{-8}$.

From the definition of level set function \eqref{ls}, we know that \eqref{hv-beta} scales a  level set function to a Heaviside  function.

Now, we establish a connection between the level set function and the Heaviside function. 

\begin{itemize}
\item {THINC scaling ($\phi\mapsto H)$}: 

We convert the level set field by  
\begin{equation}
H(\textbf{x})=\frac{1}{2} \left(1+\tanh \left(\beta\left(\mathcal {P}(\textbf{x}) \right)\right)\right), 
\label{thinc-p}
\end{equation}
where $\mathcal {P}(\textbf{x})$ is a polynomial 
\begin{equation}
\mathcal {P}(\textbf{x})=\sum_{r,s,t=0}^{p}{a}_{rst}{x}^{r}{y}^{s}{z}^{t}
\label{ls-poly}
\end{equation}
whose coefficients are determined from the level set function through the following constraints.
\begin{equation}
    \frac{\partial^D \mathcal {P}(\textbf{x})}{\partial x^{d_x}\partial y^{d_y}\partial z^{d_z}}=    \frac{\partial^D \phi (\textbf{x})}{\partial x^{d_x}\partial y^{d_y}\partial z^{d_z}}, \ \ ({d_x}, {d_y},  {d_z})=0,1,2\cdots \ \text{and } \ {d_x}+{d_y}+{d_z}=D.
    \label{ls-constraint}
    \end{equation}
In practice, we calculate the coefficients of $\mathcal {P}(\textbf{x})$  via numerical approximations using the discrete level set field available in the computational domain.   

We refer to \eqref{thinc-p} as the THINC scaling formula, and \eqref{ls-poly} as the surface polynomial or level set polynomial. Formulae   \eqref{ls-poly}  and \eqref{ls-constraint}  imply that the level set field is approximated by a polynomial function which can be of arbitrary order in principle. For example, 2nd-order (quadratic) polynomial \cite{ii2012interface,xie2017toward}, as well as 4th- and 6th-order polynomials \cite{qian2018}, have been used in previous versions of THINC method. 

\item {Inverse THINC scaling ($H\mapsto \phi)$}: 

Given the THINC function, we can directly compute the corresponding level set function by 
\begin{flalign} 
{\phi}(\textbf{x})=\frac{1}{\beta}{\tanh}^{-1}\left(2{H}(\textbf{x})-1\right) \  \ {\rm or} \ \ {\phi}(\textbf{x})=\frac{1}{\beta}\ln\left(\frac{{H}(\textbf{x})}{(1-{H}(\textbf{x})}\right). 
\label{THINC_inv}
\end{flalign}
Formula \eqref{THINC_inv} is referred to as the inverse THINC scaling that converts the Heaviside function to the level set function, which provides a continuous level set field in the target cell $\Omega_i$, and thus the value to transport the level set field via the semi-Lagrangian step in the THINC-scaling scheme described later. 

\end{itemize}

\begin{description}
\item {Remark 1.} Formulae \eqref{thinc-p} and \eqref{THINC_inv} provide analytical relations to uniquely convert between a level set field and a VOF field, which unifies the two under a single framework handleable with conventional mathematical analysis tool, and more importantly allows us to build interface-capturing schemes that take advantages from both VOF and level set methods. 

\item {Remark 2.} The interface is defined by  
\begin{equation}
\mathcal {P}(\textbf{x}) =0,  
\label{p0}
\end{equation}
where the surface polynomial $\mathcal {P}(\textbf{x})$ defined in \eqref{ls-poly} enables to accurately formulate the geometry of the interface. In principle, we can use arbitrarily high order surface polynomial to represent the interface straightforwardly without substantial difficulty. 

\item {Remark 3.} A finite value of the steepness parameter ${\beta}$ modifies the Heaviside function to a continuous and differentiable function \eqref{thinc-p}, which  serves an adequate approximation to the VOF function with adequately large steepness parameter $\beta$. 

\end{description}

The above observations implies the possibility that both level set and VOF fields can be seen as the two faces of the THINC function. So, we can expect to unify the two methods into a single scheme with superior numerical solution if a high-quality THINC function can be maintained in the numerical procedure, which is described in the next section.

\section{The THINC-scaling scheme for moving interface capturing}

We assume that the moving interface is transported by a velocity field $\textbf{u}$. Thus, the two indicator functions  are advected in a Eulerian form by the following equations, i.e. 
\begin{flalign}
\frac{\partial H}{\partial t}+ \nabla\cdot\left(\textbf{u}H\right)=H\nabla\cdot\textbf{u}
\label{h-adv}
\end{flalign}
for the THINC function ${H\left(\textbf{x},t\right)}$, and 
\begin{flalign}
\frac{\partial \phi}{\partial t}+ \textbf{u}\cdot\nabla\phi=0 
\label{ls-adv}
\end{flalign}
for the level set function ${\phi\left(\textbf{x},t\right)}$. 

It is noted that \eqref{h-adv}, with the left hand side being of the flux form, allows the direct use of the finite volume formulation to ensure the conservativeness of advection transport. Thus, for incompressible flow ($\nabla\cdot\textbf{u}=0$), the VOF field is numerically conserved.   

Next, we present the THINC-scaling scheme to simultaneously solve \eqref{h-adv} and \eqref{ls-adv}. %We focus on two dimensional case which can be extended to three dimensions straightforwardly.   
The computational domain is composed of non-overlapped discrete grid
cells $\Omega_{i}\,(i=1,2,\ldots,N)$ of the volume $\left|\Omega_{i}\right|$, which can be either structured or unstructured grids.  For any target cell element $\Omega_i$, we denote its mass 
center by $\textbf{x}_{ic}=(x_{ic},y_{ic},z_{ic})$, and its $J$ surface segments of areas $\left|\Gamma_{ij}\right|$ by $\Gamma_{ij}$ with $j=1,2,\ldots,J$. The outward unit normal is denoted by  ${\mathbf{n}}_{ij}=(n_{xij},n_{yij},,n_{zij})$. 

Assume that we know at time step $n$ ($t=t^n$) the VOF value   
\begin{equation}
\bar{H}_i^n=\frac{1}{|\Omega_i|}\int_{\Omega_i}H(\textbf{x},t^n)d\Omega,
\label{vofn}
\end{equation}
for each cell, and the level set value 
\begin{equation}
\phi_i^n=\phi(x_{ic},y_{ic},z_{ic},t^n)
\label{lsn}
\end{equation}
at each cell center, we use the third-order TVD Runge-Kutta scheme \cite{shu88} for time
integration to update both VOF and level set values, $\bar{H}_i^{n+1}$ and $\phi_i^{n+1}$, to the next time step $n+1$ ($t=t^{n+1}=t^{n}+\Delta t$).

We hereby summarize the solution procedure of the THINC-scaling scheme for one  Runge-Kutta sub-step that advance $\bar{H}_i^{m}$ and $\phi_i^{m}$ at sub-step $m$ to $\bar{H}_i^{m+1}$ and $\phi_i^{m+1}$ at sub-step $m+1$.   

\begin{description}
    \item {\bf Step 1.} Compute the level set polynomial of $p$th order for cell $\Omega_{i}$ from the level set field using the constraint condition \eqref{ls-constraint},
\newline
\begin{equation}
\mathcal{P}_{i}\left({\textbf x} \right)=\sum_{r,s,t=0}^{p}{a}_{rst}{X}^{r}{Y}^{s}{Z}^{t}
\end{equation}
\par where $\left(X,Y,Z\right)$ is the local coordinates with respect to the center of cell ${\Omega}_{i}$, i.e.  $X={x}-{x}_{ic}$, $Y={y}-{y}_{ic}$, $Z={z}-{z}_{ic}$. The coefficients ${a}_{rst}$ are computed by Lagrange interpolation or least square method using the level set values $\phi_i^{m}$ in the target and nearby cells. 

  \item {\textbf {Step 2}}:
Construct the cell-wise THINC function {\color{black} under the constraint of volume conservation using the VOF value by }
\newline
\begin{flalign} \label{mass_constraint}
\frac{1}{|{\Omega}_{i}|}{\int}_{{\Omega}_{i}}\frac{1}{2}\left(1+\tanh\left(\beta\left(\mathcal{P}_{i}\left({\textbf x}\right)+{\phi}_{i}^{\Delta}\right)\right)\right)d \Omega={\bar H}_{i}^{m}. 
\end{flalign}
With a pre-specified $\beta$ and the surface polynomial obtained at step 1, the only unknown ${\phi}_{i}^{\Delta}$ can be computed from \eqref{mass_constraint}. In practice, we use the numerical quadrature detailed in \cite{xie2017toward}, and the resulting nonlinear algebraic function of ${\phi}_{i}^{\Delta}$ is solved by the Newton iterative method. See appendix A of this paper for details. 

We then get the THINC function 
\begin{flalign} \label{thincf_mass_constraint}
H^m_{i}({\bf x})=\frac{1}{2}\left(1+\tanh\left(\beta\left(\mathcal{P}_{i}\left({\textbf x}\right)+{\phi}_{i}^{\Delta}\right)\right)\right),
\end{flalign}
which satisfies the conservation constraint of the VOF value, and ${\phi}_{i}^{\Delta}$ is a correction to the interface location due to the conservation constraint. 
The piece-wise interface in each interface cell is defined by  
\begin{flalign} \label{cell-p}
\psi_{i}\left({\textbf x}\right)\equiv\mathcal{P}_{i}\left({\textbf x}\right)+{\phi}_{i}^{\Delta}=0. 
\end{flalign}
We refer to \eqref{cell-p} as the Polynomial Surface of the Interface (PSI) equation in cell ${\Omega}_{i}$. 

\item{\textbf {Step 3}}: Update the VOF function by solving \eqref{h-adv} through the following finite volume formulation, 
\begin{equation}
\bar{H}^{m+1}_{i}=\bar{H}^{m}_{i}-\frac{\Delta t}{\left|\Omega_{i}\right|}\sum_{j=1}^{J}\left(\int_{\Gamma_{ij}}\left(({\bf u} \cdot {\bf n}) H^m_{i}({\bf x})_{iup}\right)d\Gamma\right)+\frac{\bar{H}^m_{i}}{\left|\Omega_{i}\right|}\sum_{j=1}^{J}\left({\bf u} \cdot {\bf n})_{ij}\left|\Gamma_{ij}\right|\right)\Delta t,\label{fvm_semid-1}
\end{equation}
where the upwinding index $iup$ is determined by 
\begin{equation}
\begin{split}iup=\begin{cases}
=i,\ {\rm for}\ ({\bf u} \cdot {\bf n})_{ij}>0;\\
=ij,\ {\rm otherwise},
\end{cases}\end{split}
\label{eq:upwind_index}
\end{equation}
and $ij$ denotes the index of the neighboring cell that shares cell boundary 
$\Gamma_{ij}$ with target cell $\Omega_{i}$. The integration on cell surface is computed by Gaussian quadrature formula. 

\item{\textbf {Step 4}}: Update level set value at cell center using a semi-Lagrangian method as follows. 
\begin{description}
    \item {\bf Step 4.1}: 
We first find the departure point ${\textbf x}_{id}$ for each cell center ${\textbf x}_{ic}$ by solving the initial value problem, 
\begin{equation}
\begin{split}\begin{cases}
\displaystyle {\frac{{d}{\textbf x}}{d\tau}=-{\textbf u}\left({\textbf x},t^n+\tau\right)} \\
{\textbf x}(0)={\textbf x}_{ic}
\end{cases}\end{split}
\label{ivp}
\end{equation}
up to $\tau=\Delta t =t^{n+1}-t^n$, which leads to ${\textbf x}_{id}={\textbf x}(\tau)$. In the present work, we use a second-order Runge-Kutta method to solve the ordinary differential equation \eqref{ivp}.
    \item {\bf Step 4.2}: Update the level set value  ${\phi}_{i}$ at cell center ${\textbf x}_{ic}$   using the Lagrangian invariant solution, 
\begin{flalign} \label{ls_sl}
{\phi}_{i}^{m+1}=\tilde{\phi}_{id}^{m}\left({\bf x}_{id}\right),   
\end{flalign}
where $\tilde{\phi}_{id}^{m}$ is the level function on cell $\Omega_{id}$ where the departure point ${\textbf x}_{id}$ falls in. 
Using the inverse THINC-scaling formula \eqref{THINC_inv}, we immediately get the level set function from \eqref{thincf_mass_constraint},
\begin{flalign} \label{ls-THINC}
\tilde{\phi}_{id}^{m}\left({\bf x}\right)=\frac{1}{\beta}{\tanh}^{-1}\left(2H^m_{id}({\bf x})-1\right), \end{flalign}
\end{description}
which gives the level set value everywhere in cell $\Omega_{id}$ that includes the departure point ${\textbf x}_{id}$. 
%%%%%%%%%%%
\item{\textbf {Step 5}}:
Reinitialize the level set field. 
We fix the level set values computed from \eqref{ls-THINC} for the interface cells which are identified by  ${\epsilon}_{1}{\leq}{\bar H}^{m}_{i}{\leq}{1-{\epsilon}_{2}}$ with ${\epsilon}_{1}$ and ${\epsilon}_{2}$ being  small positive numbers.  The level set values at the cell centers away from the interface region are reinitialized to satisfy the Eikonal equation, 
\begin{equation}
|\nabla\phi|=1. 
\end{equation}

We use the fast sweeping method (FSM) \cite{zhao2005fast} on structured grid and the iteration method \cite{sussman1994level,Dianat2017unstructured} on unstructured grid for reinitializing level set values. 

%%%%%%%%%%%
\item{\textbf {Step 6}}:
Go back to Step 1 for next sub-time step calculations.
%%%%%%%%%%%

\end{description}{}

\begin{description}
\item {Remark 4.} The THINC-scaling scheme shown above unifies the VOF method and level set method. Eq.\eqref{thincf_mass_constraint} retrieves the VOF field from the level set field, while \eqref{ls-THINC} retrieves the level set field from the VOF field with numerical conservativeness.

\item {Remark 5.} The inverse THINC-scaling \eqref{ls-THINC} facilitates a semi-Lagrangian solution without any spatial reconstruction or interpolation, such as those used in \cite{strain99}. This step essentially distinguishes the present scheme from the coupled THINC/level set method in \cite{qian2018}, where the level set function is updated by a fifth-order Hamilton-Jacobi WENO scheme with a 3rd-order TVD Runge-Kutta time-integration scheme. {\color{black} See appendix B for a detail comparison between THINC-scaling and  THINC/level set methods. }

\item {Remark 6.} Even without explicit geometrical reconstruction in THINC-scaling method, the interface is clearly defined by the PSI equation \eqref{cell-p}, $\psi_i({\textbf x})=0$, within the interface cells, which provides the sub-cell interface structure with geometrical information, such as position, normal direction and curvature to facilitate the computation of so-called sharp-interface formulation. It distinguishes the present scheme with superiority from any other algebraic interface-capturing methods. 

\item {Remark 7.} As discussed in \cite{xiao2011revisit}, the steepness parameter ${\beta}$ can be estimated by the thickness of the jump transition across the interface using 
\begin{equation}
 \beta=\frac{1}{\eta}{\tanh}^{-1}\left(1-2\varepsilon\right)   
\end{equation}
where $\eta$ denotes the normalized half thickness of the jump with respect to the cell size, and $\varepsilon$ is a small positive number to define the range of interface transition layer in terms of the VOF value, i.e. $\varepsilon \le H(\textbf{x})\le 1-\varepsilon$, we use $\varepsilon={10}^{-8}$ in this work. {\color{black}For example, in order to keep a 3-cell thickness for the interface, $\eta$ can be set as $1.5\Delta$ with $\Delta$ being the cell size, which then approximately results in $\beta\approx6/\Delta$. As $\beta\propto 1/\Delta $, refining grid resolution effectively increases  $\beta$ and makes the THINC function approach to the Heaviside function of VOF.  Our numerical experiments show that using such a $\beta$ the THINC method can resolve interfaces with about two or three mesh cells.  }

\item {Remark 8.} With a reasonably large $\beta$, the THINC function maintains a steepness across the interface transition layer and effectively removes numerical diffusion (smearing). As discussed above, a compact thickness of 2 or 3 cells can be always maintained during the computation by simply updating the VOF field by the finite volume scheme \eqref{fvm_semid-1},  without any extra artificial compression or anti-diffusion manipulation.   

{\color{black}
\item {Remark 9.} The numerical procedure to determine ${\phi}_{i}^{\Delta}$ needs to solve \eqref{mass_constraint}. In current formulation, the integration of multi-dimensional THINC function is approximated by a numerical quadrature for simplicity,  which might cause the degradation of numerical accuracy, and remains an unsolved open problem worth further efforts. 

\item {Remark 10.} In order to preserve the location of interface from being changed during the reinitialization process, we fix the level set values in the interface cells while reinitializing the level set field for the outside region without any special treatment. So, the level set field within the interface cells might not rigorously satisfy $|\nabla \phi=1|$. Our numerical tests show that the  level set field  near the interface remains acceptable quality even without any extra treatment.
}

\end{description}

%%%%%%%%%%%%%%%%%%%%%%%%%%%%%%%%%%%%%%%%%%%%%%%%%%%%%%%%%%%%%%%%%%%%%%%%%%%%%%%%%%%%%%%%%%%%%%%%%%%%
\section{Numerical tests}

We verify the THINC-scaling scheme to capture moving interfaces using some advection benchmark tests. We focus on the numerical tests in two and three dimensions on both structured (Cartesian) grid and unstructured (triangular/tetrahedral) grid. 
The numerical errors in terms of the VOF field are quantified via the $L_1$ error norm 
\begin{align}
E(L_1) =\sum_{i=1}^{N_e}|\bar{H}_{i}-\bar{H}_{i}^{e}||\Omega_{i}| \label{L1}
\end{align}
and relative error $E_r$
\begin{align}
E_r =\frac{\sum_{i=1}^{N_e}\left|\bar{H}_{i}-\bar{H}_{i}^{e}\right|\left|\Omega_{i}\right|}{\sum_{i=1}^{N_e}\left|\bar{H}_{i}^{e}\right||\left|\Omega_{i}\right|},\label{Lr}
\end{align}
where $N_e$ is the total number of elements in the computational domain and, $\bar{H}_{i}$ and $\bar{H}^e_{i}$ are respectively the numerical and exact VOF values.

\begin{figure}[htbp]
	\centering
	\subfigure[$\text{2D Cartesian mesh}$] {
	\includegraphics[width=0.30\textwidth]{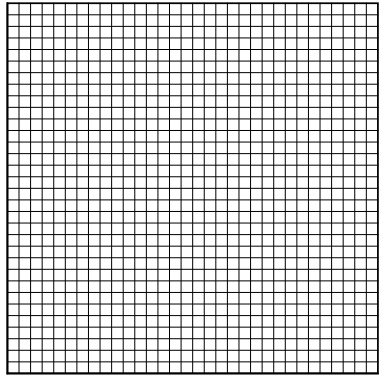}} \hspace{0.4cm}
	\subfigure[$\text{2D triangular mesh}$] {
	\includegraphics[width=0.30\textwidth]{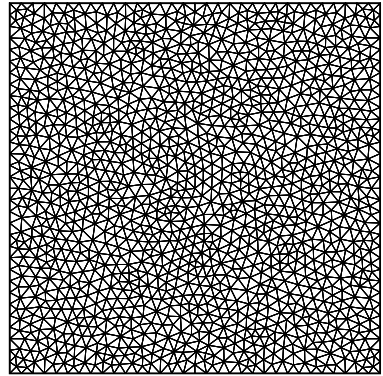}} \\ 
	\subfigure[$\text{3D Cartesian mesh}$] {
		\centering
		\includegraphics[width=0.30\textwidth]{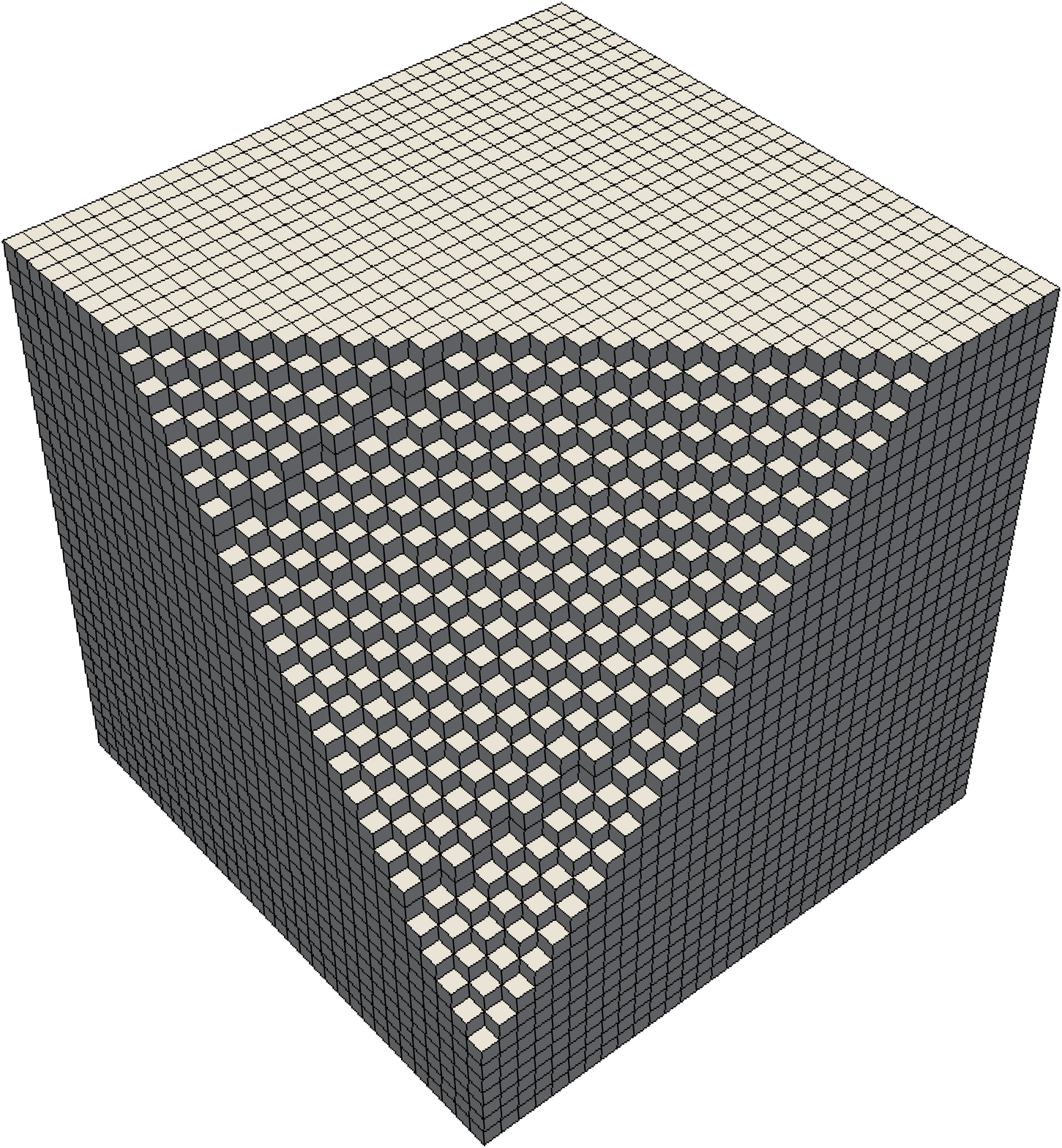}} \hspace{0.4cm}
	\subfigure[$\text{3D tetrahedral mesh}$] {
		\centering
		\includegraphics[width=0.30\textwidth]{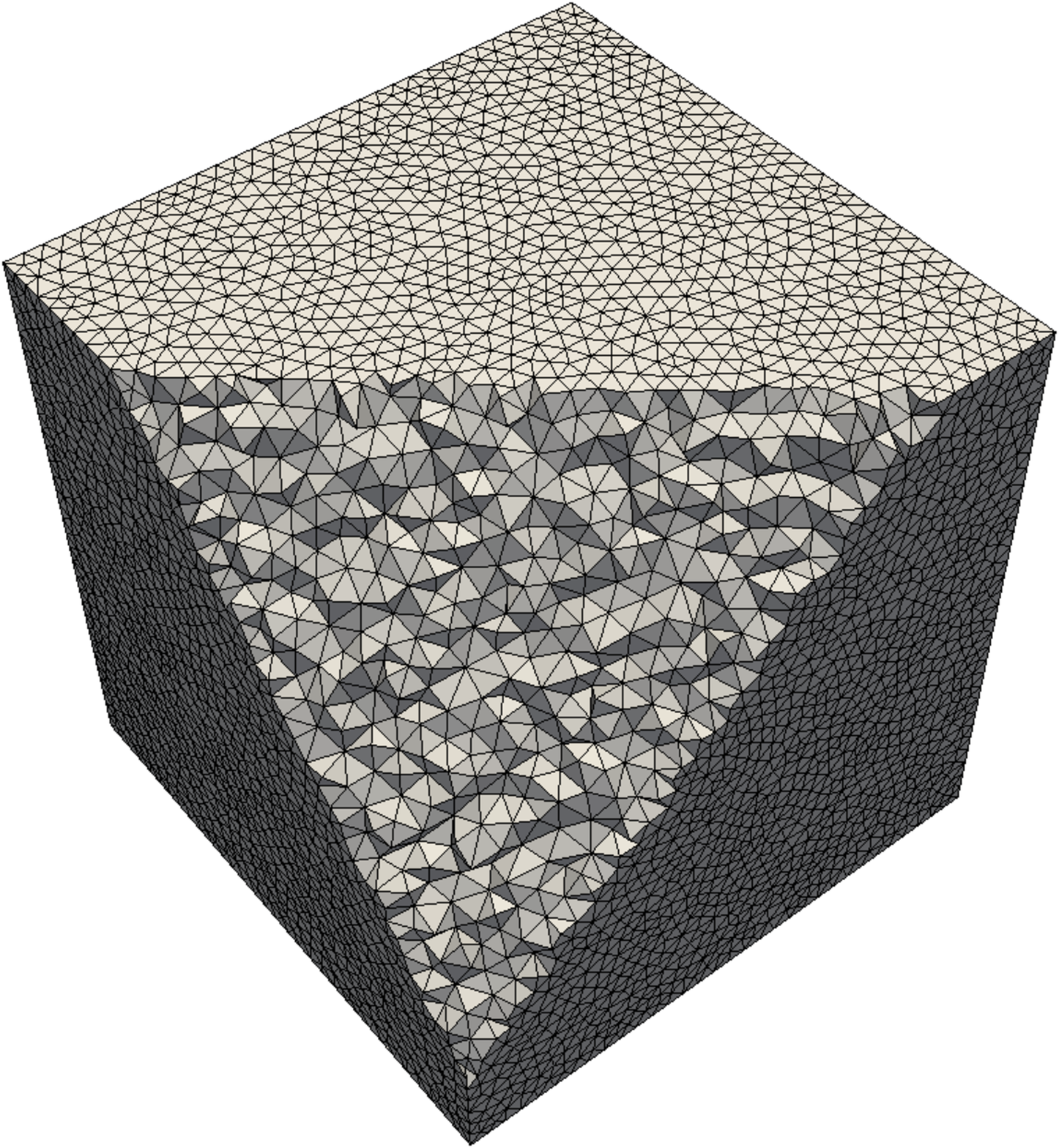}} %\hspace{0.4cm}
	\caption{\small{Mesh configurations used for benchmark tests.}}
\label{meshes}
\end{figure}

Unstructured 2D triangular and 3D tetrahedral meshes are generated using open source mesh generator software Gmsh version 4.3.0\cite{gmsh}. The quality of tetrahedral elements are optimized using Netgen option available in the software. We generate the mesh by taking uniformly distributed nodes on each boundary of computational domain. For convenience, we denote the resolution of both Cartesian and unstructured meshes by N, the number of uniformly distributed nodes on each boundary edge of computational domain. Fig.\ref{meshes} shows different mesh configurations used for benchmark tests. The number of elements with respect to mesh resolutions used in different test cases are shown in Table \ref{tri-ele} (2D) and Table \ref{tet-ele} (3D).

\begin{table}[h]
\centering
\small
\caption{Number of elements for different unstructured triangular grid resolutions used in 2D solid body rotation and vortex deformation transport benchmark tests}\label{tri-ele} 
\begin{tabular}[t]{|c|c|c|c|} 
\hline
\multicolumn{2}{|c|}{2D solid body rotation} & \multicolumn{2}{|c|}{2D vortex deformation}\\ 
\hline
{N} & {Number of elements} & {N} & {Number of elements} \\
\hline
$50$  & $6{,}670$   & $32$  &  $2{,}654$ \\
$100$ & $26{,}436$  & $64$  &  $10{,}776$  \\
$200$ & $105{,}724$ & $128$ &  $43{,}202$\\
\hline
\end{tabular}
\end{table}
\begin{table}[h]
\centering
\small
\caption{Number of elements for different unstructured tetrahedral grid resolutions used in 3D deformation flow and 3D shear flow benchmark tests}\label{tet-ele} 
\begin{tabular}[t]{|c|c|c|c|} 
\hline
\multicolumn{2}{|c|}{3D deformation flow} & \multicolumn{2}{|c|}{3D shear flow}\\ 
\hline
{N} & {Number of elements} & {N} & {Number of elements} \\
\hline
$32$  & $131{,}383$   & $32$  & $292{,}482$\\
$64$  & $1{,}163{,}273$ & $64$  & $2{,}314{,}437$\\
\hline
\end{tabular}
\end{table}

{\color{black}
%As suggested by the reviewer, we explicitly give the values of $\beta$ for all test cases.  

As discussed before, in order to make the jump thickness of reconstructed THINC function to be within 3 mesh cells for a given grid resolution, $\beta$ for each test case is determined by $\beta=6 / \Delta$ with $\Delta$ being the average size of mesh cells. Consequently, a finer mesh results in  a thinner transition layer for the interface. It guarantees that the THINC function converges to the exact Heaviside function as $\Delta\rightarrow 0$. 

We show below in Tables \ref{beta-table2d} and \ref{beta-table3d} the $\beta$ values for the structured and unstructured grids used in the numerical tests presented in this paper. All $\beta$ values are large enough to allow the THINC function to adequately mimic a jump-like profile. Moreover, $\beta$ becomes larger and generates steeper jump on finer grids.   

\begin{table}[h]
\centering
\caption{Values of $\beta$ for 2D grids. } 
\begin{tabular}{|c|c|c|c|c|}
\hline Mesh type& \multicolumn{2}{c|}  { Cartesian mesh } & \multicolumn{2}{c|} { Triangular unstructured mesh } \\
\hline $\mathrm{N}$ & $\Delta$ & $\beta$ & $\Delta$ & $\beta$ \\
\hline 32 & $3.12 \times 10^{-2}$ & 192 & $1.67 \times 10^{-2}$ & 361.31 \\
\hline 64 & $1.56 \times 10^{-2}$ & 384 & $8.24 \times 10^{-3}$ & 727.88 \\
\hline 128 & $7.81 \times 10^{-3}$ & 768 & $4.11 \times 10^{-3}$ & 1456.57 \\
\hline
\end{tabular}
\label{beta-table2d}
\end{table}
\begin{table}[h]
\centering
\caption{Values of $\beta$ for 3D grids. } 
\begin{tabular}{|c|c|c|c|c|}
\hline Mesh type& \multicolumn{2}{c|} { Cartesian mesh } & \multicolumn{2}{c|} { Tetrahedral unstructured mesh } \\
\hline $\mathrm{N}$ & $\Delta$ & $\beta$ & $\Delta$ & $\beta$ \\
\hline 32 & $3.12 \times 10^{-2}$ & 192 & $3.35 \times 10^{-2}$ & 179.18 \\
\hline 64 & $1.56 \times 10^{-2}$ & 384 & $1.69 \times 10^{-2}$ & 355.16 \\
\hline
\end{tabular}
\label{beta-table3d}
\end{table}

}

%~~~~~~~~~~~~~~~~~~~~~~~~~~~~~~~~~~~~~~~~~~~~~~~~~~~~~~~~~~~~~~~~~~~~~~~~~~~~~~~~~~~~~~~~~~~~~~~~~~~~~~~~~~~~~
\subsection{Two dimensional benchmark tests}

\subsubsection{Solid body rotation test}
In this test, so-called Zalesak's slotted disk test \cite{zalesak1979fully}, initially a circle with a radius of 0.5 centered at $\left(0.5,0.75\right)$ in a unit square computational domain is notched with a slot defined by $\left(|x-0.5|\leq0.025 \: {\text {and}} \: y\leq0.85\right)$. The slotted circle is rotated with the velocity field given by $\left(y-0.5,0.5-x\right)$. 

In order to maintain sharp interface jump, we use $\beta=6.0/\Delta$, where $\Delta$ represents the cell size and can be defined straightforwardly as $\Delta=\min\left(\Delta x,\Delta y\right)$ in case of 2D Cartesian grid. In case of unstructured triangular grid, we define the cell size by hydraulic diameter $\Delta=4A/P$, where A and P are the area and perimeter of the triangular cell element respectively. For interface reconstruction, a quadratic polynomial 
\begin{equation}
\mathcal{P}_{i}\left({\textbf x} \right)=\sum_{r,s=0}^{2}{a}_{rs}{X}^{r}{Y}^{s}
\end{equation}
is used in this test. 

\begin{figure}[htbp]
	\centering
	\subfigure[$\text{N}=50$] {
		\centering
		\includegraphics[width=0.3\textwidth]{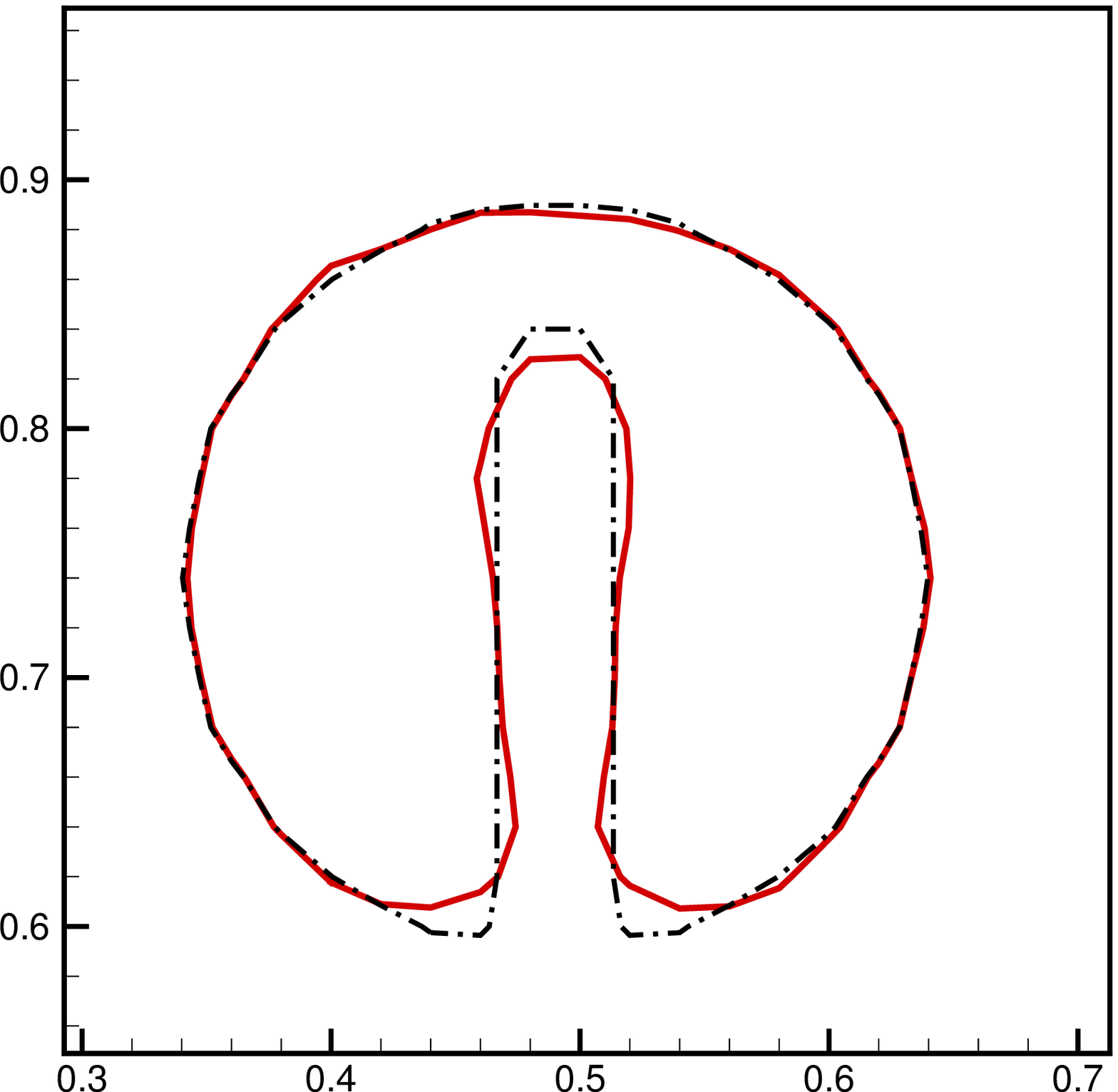}} \hspace{0.4cm}
	\subfigure[$\text{N}=100$] {
		\centering
		\includegraphics[width=0.3\textwidth]{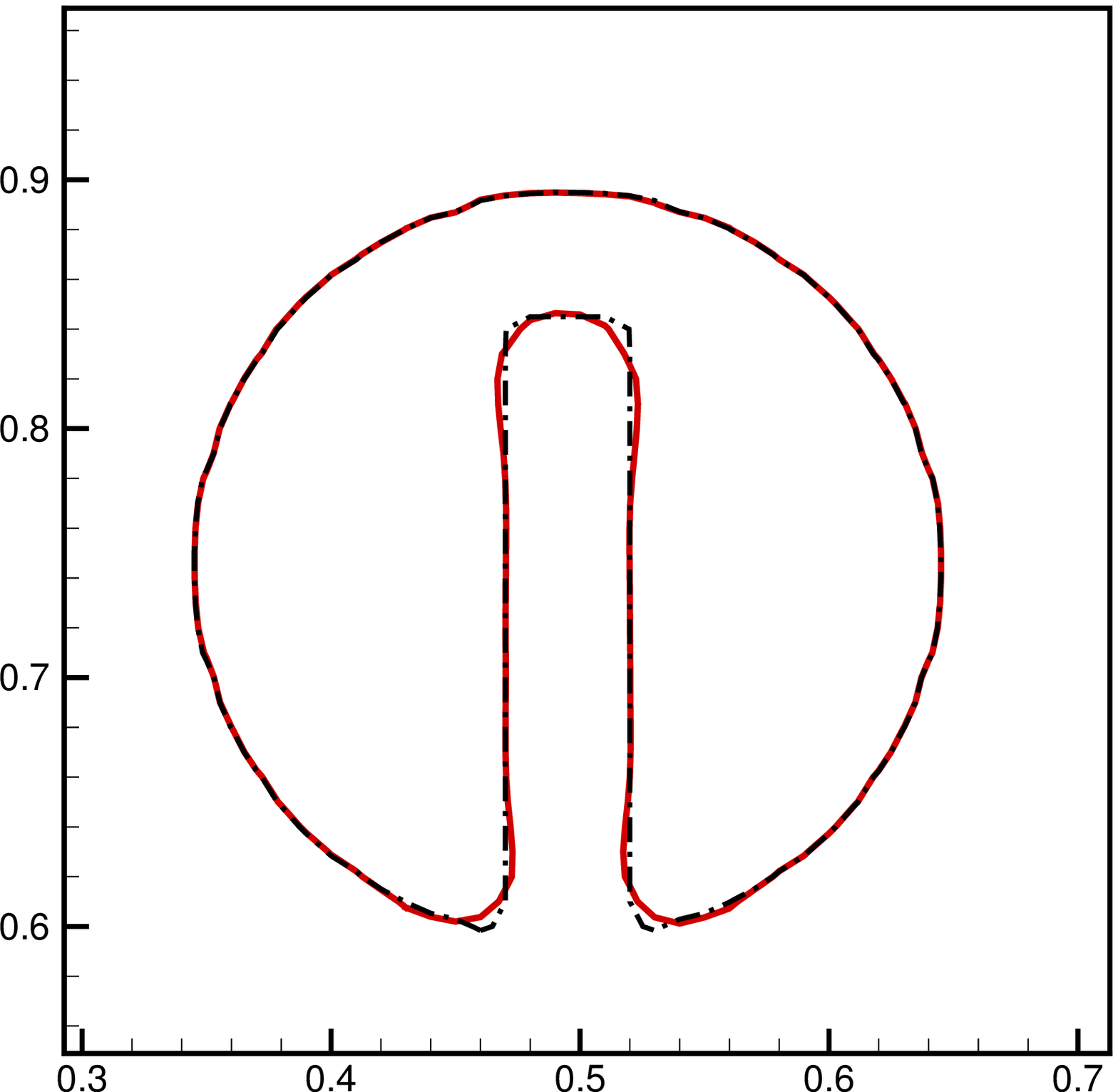}} \hspace{0.4cm}
	\subfigure[$\text{N}=200$] {
		\centering
		\includegraphics[width=0.3\textwidth]{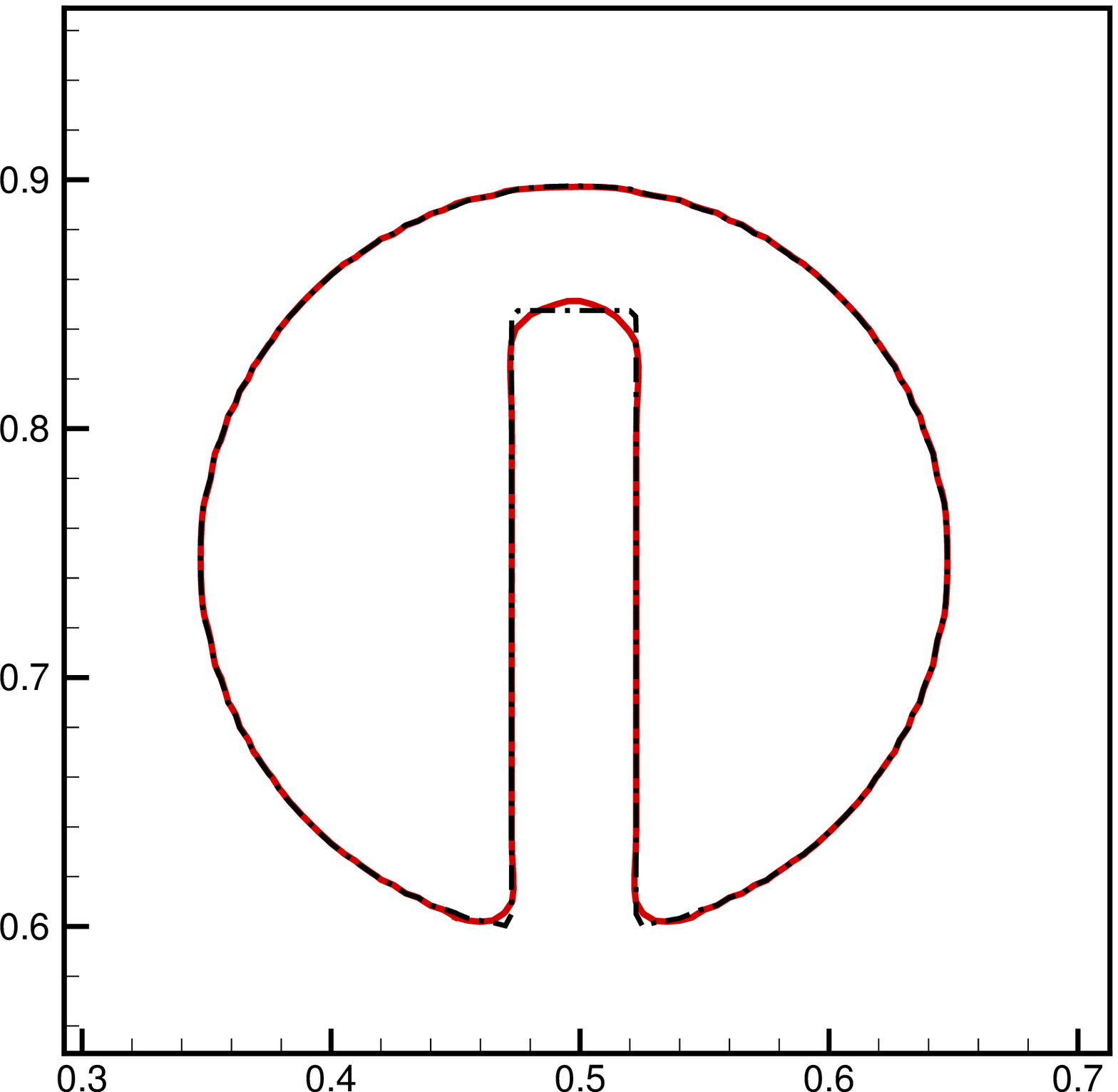}}
	\caption{VOF 0.5-contour line in Zalesak solid body rotation test after one revolution on meshes of (a) $50\times50$, (b) $100\times100$ and (c) $200\times200$ cells. The black dashed line stands for the exact solution, and the red solid line for the numerical solution. }
\label{zalesak_results}
\end{figure}

\begin{figure}[htbp]
	\centering
	\subfigure[$\text{N}=50$] {
		\centering
		\includegraphics[width=0.3\textwidth]{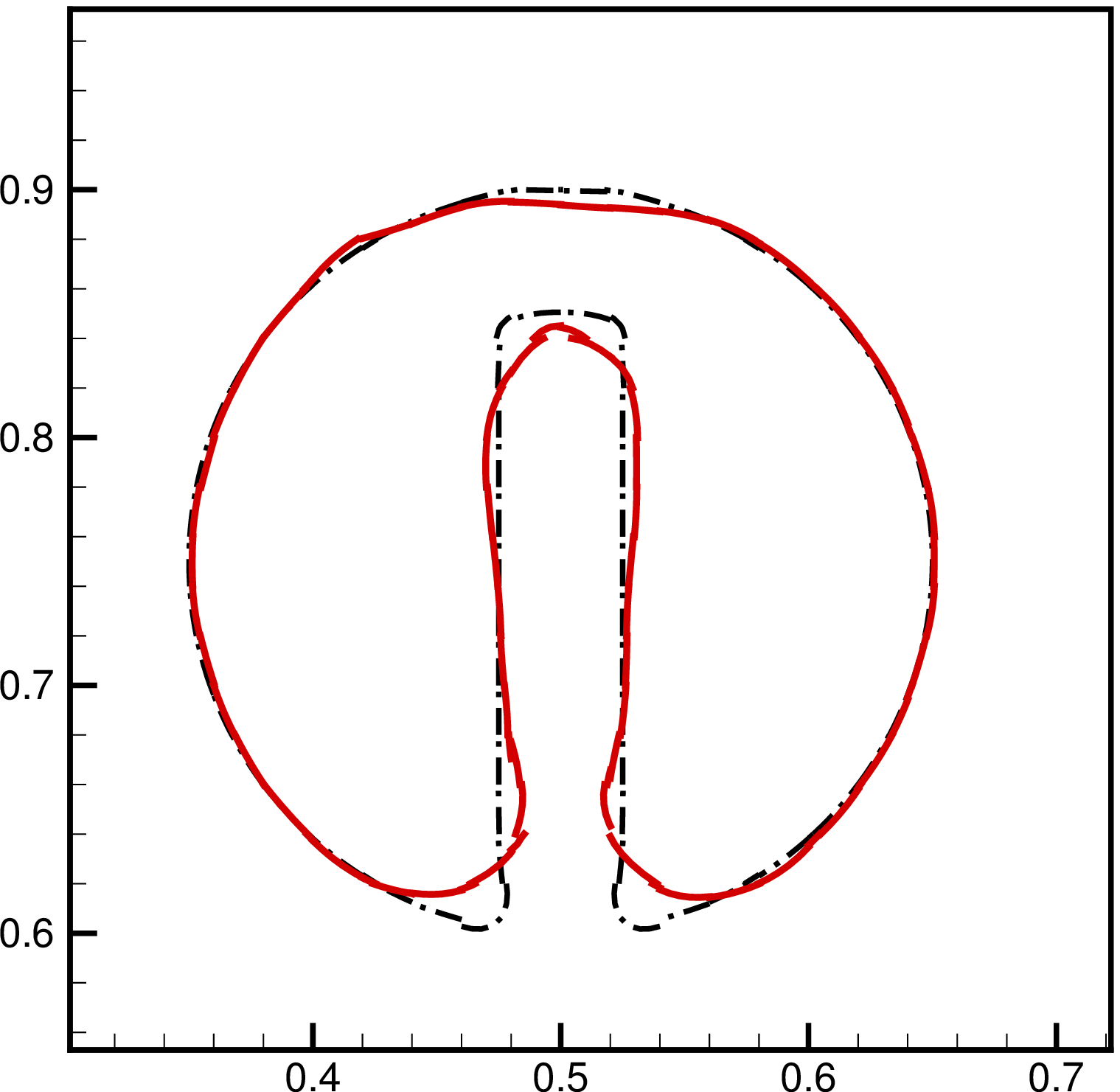}} \hspace{0.4cm}
	\subfigure[$\text{N}=100$] {
		\centering
		\includegraphics[width=0.3\textwidth]{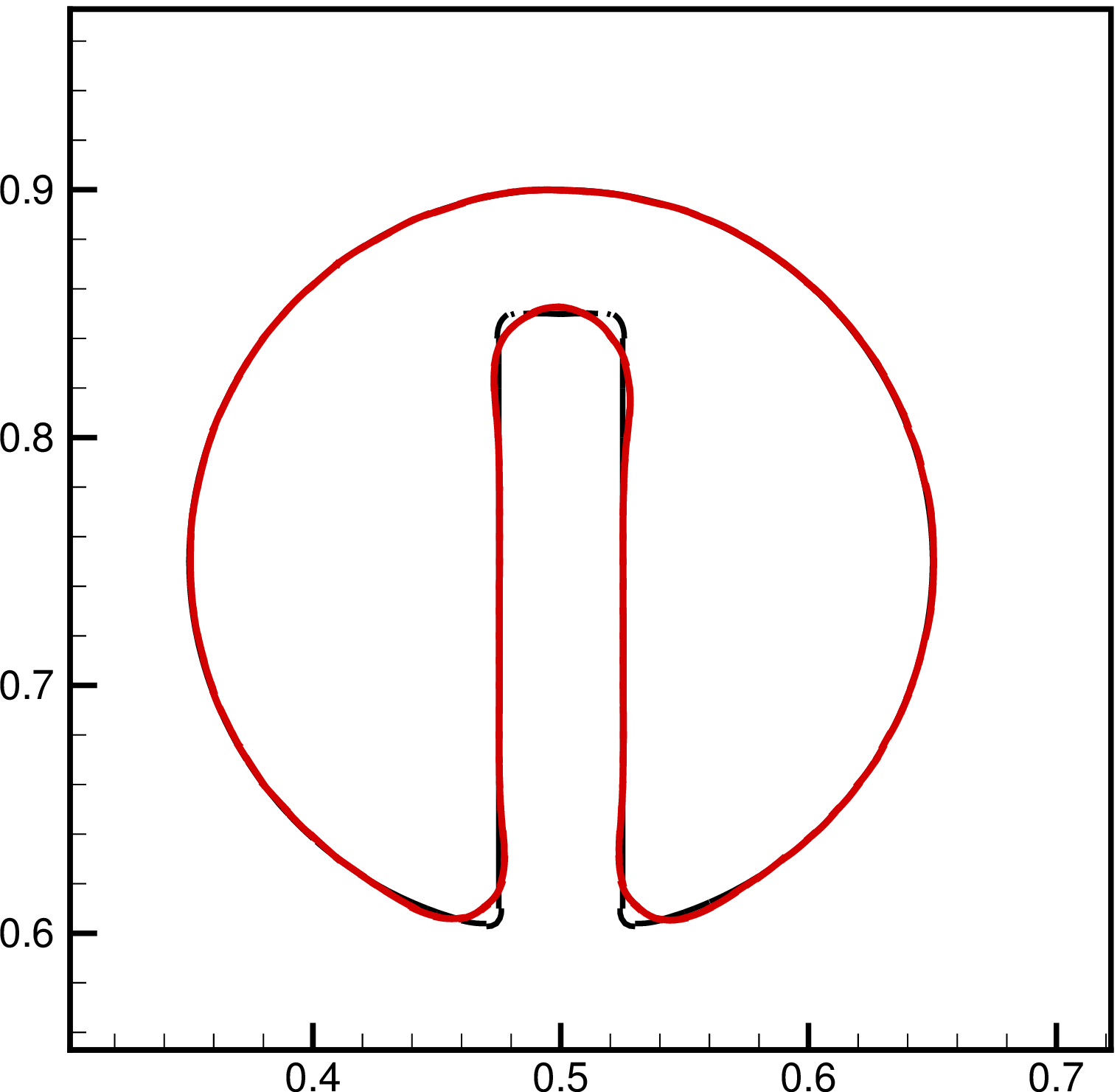}} \hspace{0.4cm}
	\subfigure[$\text{N}=200$] {
		\centering
		\includegraphics[width=0.3\textwidth]{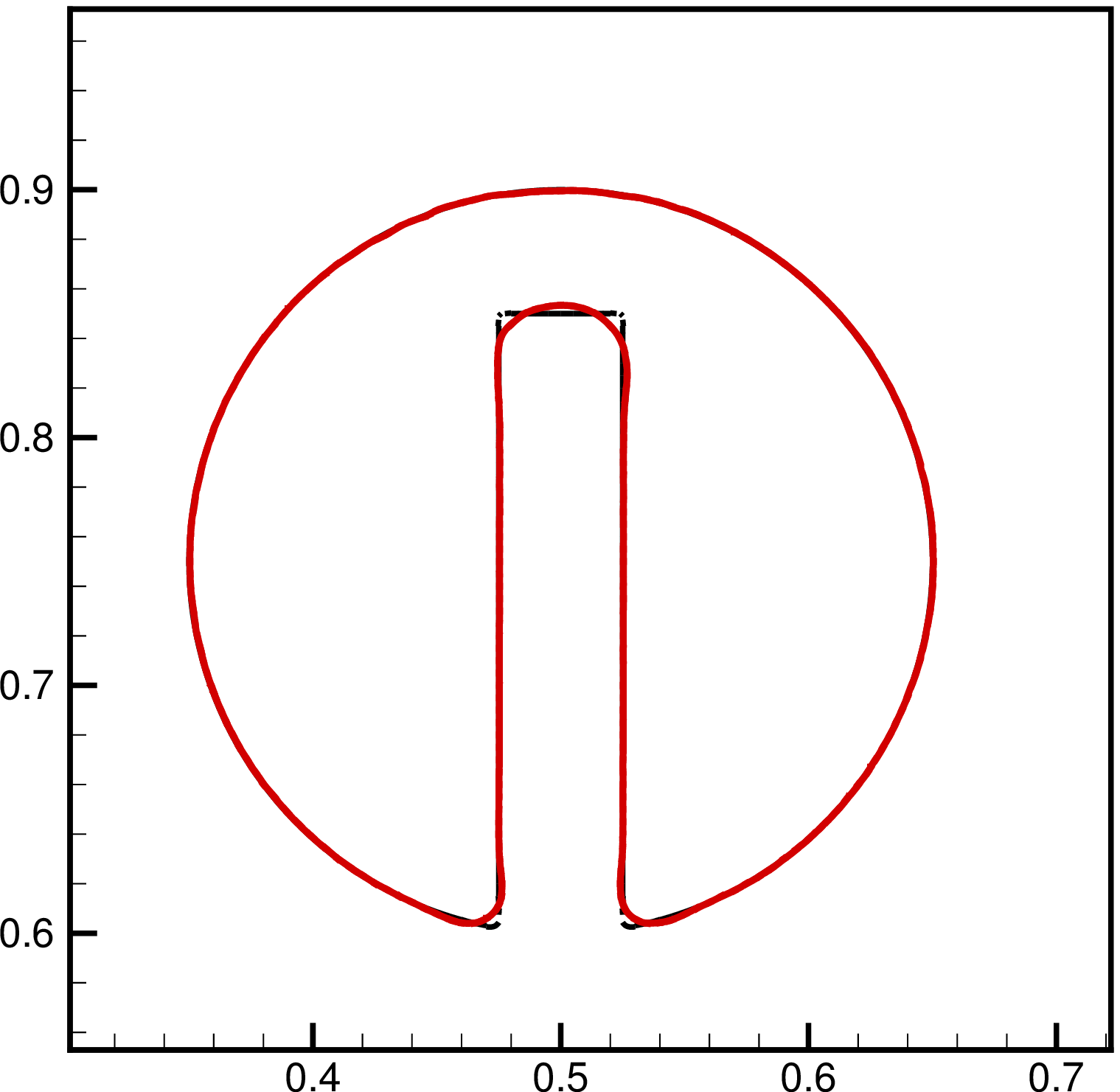}}
	\caption{The numerical results of Zalesak solid body rotation test after one revolution showing the reconstructed interfaces on meshes of (a) $50\times50$, (b) $100\times100$ and (c) $200\times200$ cells. The PSI defined by \eqref{cell-p} of the interface cells are plotted by red solid line against the exact solution (black dashed line). }
\label{zalesak_reconstructed_results}
\end{figure}

We compute this test on structured grid for different grid sizes with 50, 100 and 200 vertices evenly distributed on each edge of computational domain. Fig. \ref{zalesak_results} shows the numerical results of interface identified by VOF 0.5 contour line on Cartesian grid.  As observed from these results, the interface in the slot region is well captured. As demonstrated in  \cite{xie2017toward} and \cite{qian2018}, the quadratic polynomial representation of the interface preserves the geometrical symmetry of the solution, which deteriorates significantly if a linear function (straight line) is used as commonly observed in the results of VOF method using PLIC reconstructions.  We also show the PSI of the interface cells in Fig. \ref{zalesak_reconstructed_results}. The interface is retrieved and represented by the cell-wise quadratic curves in the interface cells, which provides a clearly defined surface in form of high-order polynomial  within each interface cell. 

\begin{figure}[htbp]
	\centering
	\subfigure[$\text{N}=50$]{
		\centering
		\includegraphics[width=0.3\textwidth]{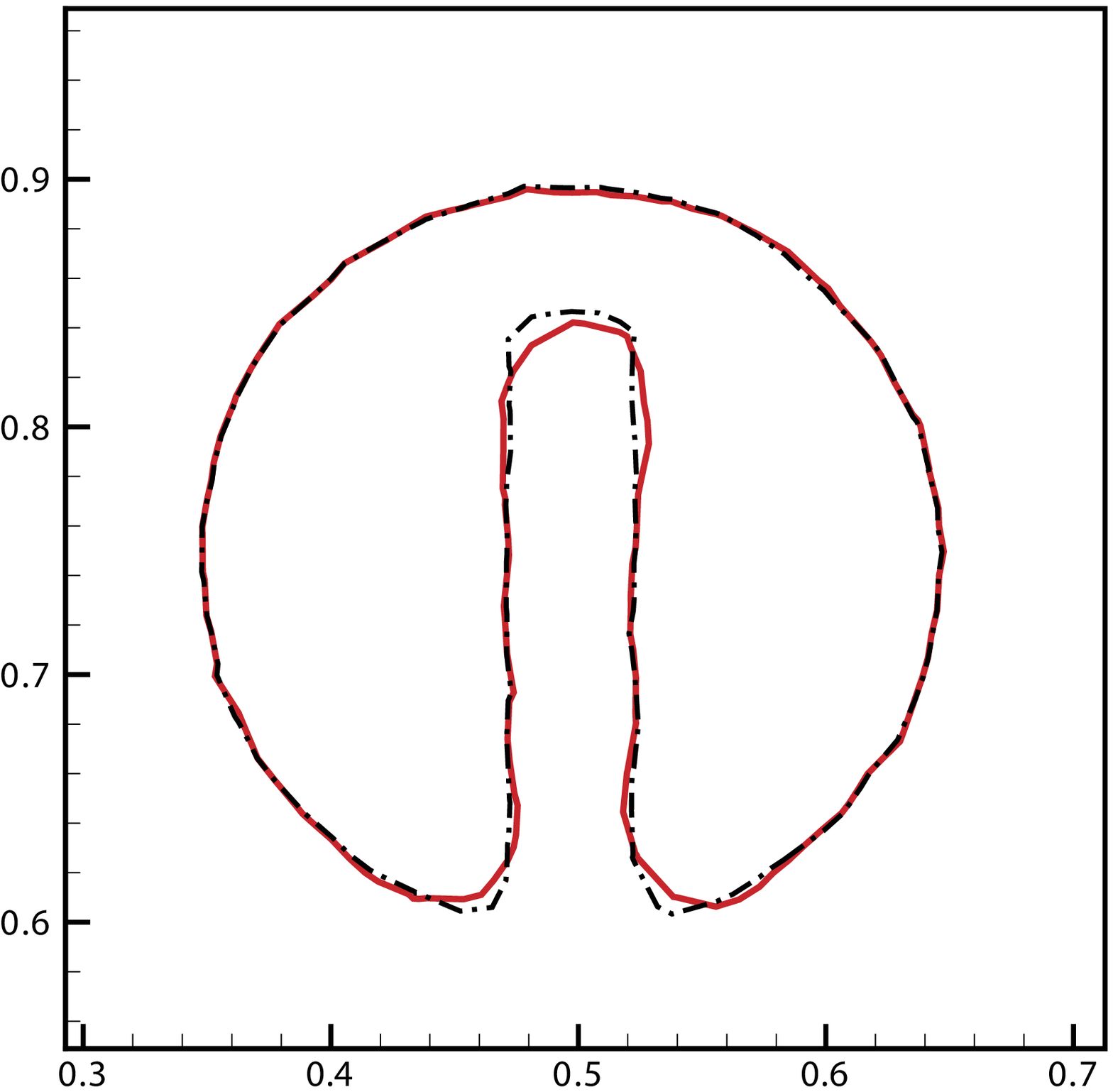}} \hspace{0.4cm}
	\subfigure[$\text{N}=100$]{
		\centering
		\includegraphics[width=0.3\textwidth]{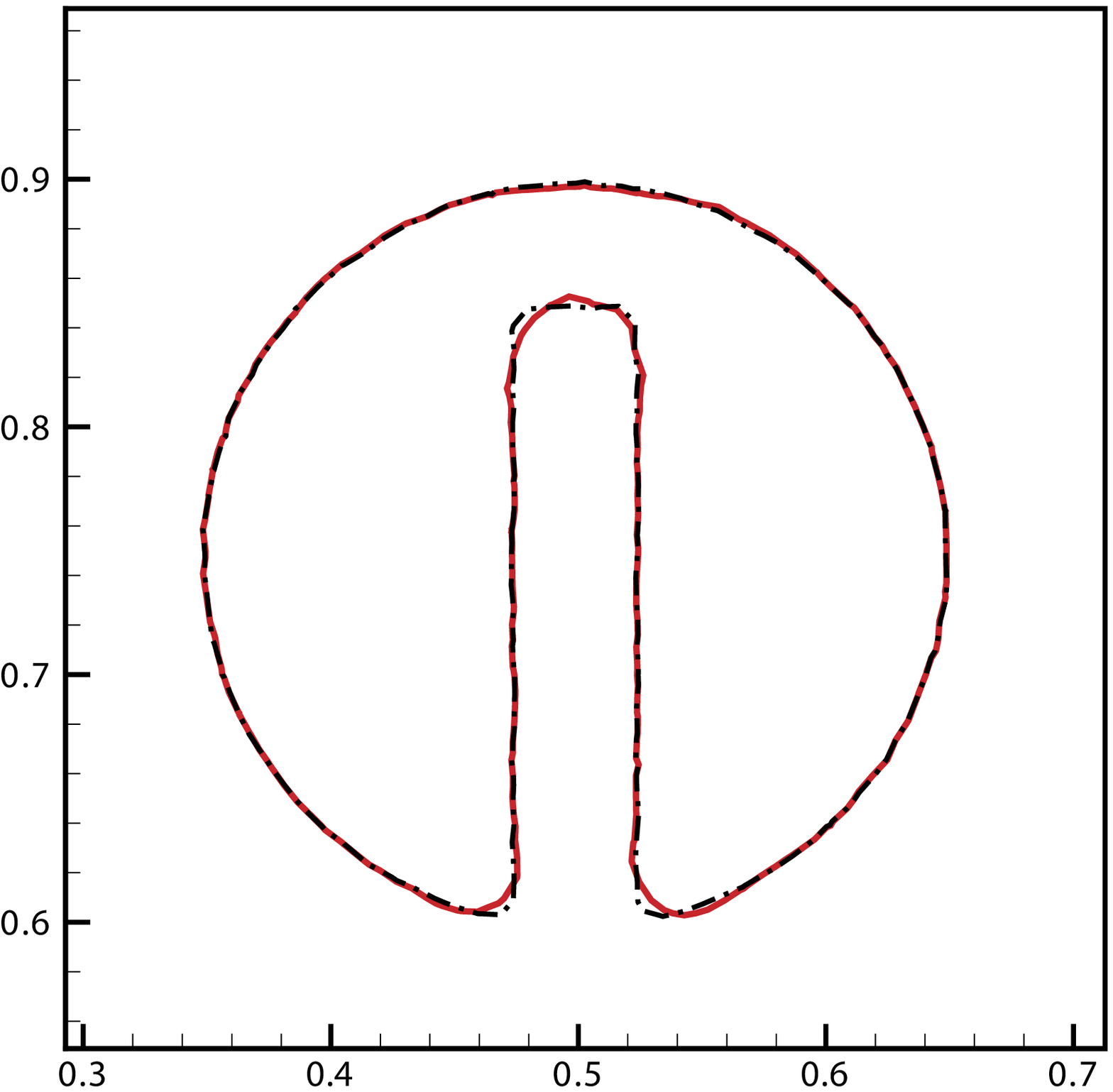}} \hspace{0.4cm}
	\subfigure[$\text{N}=200$] {
		\centering
		\includegraphics[width=0.3\textwidth]{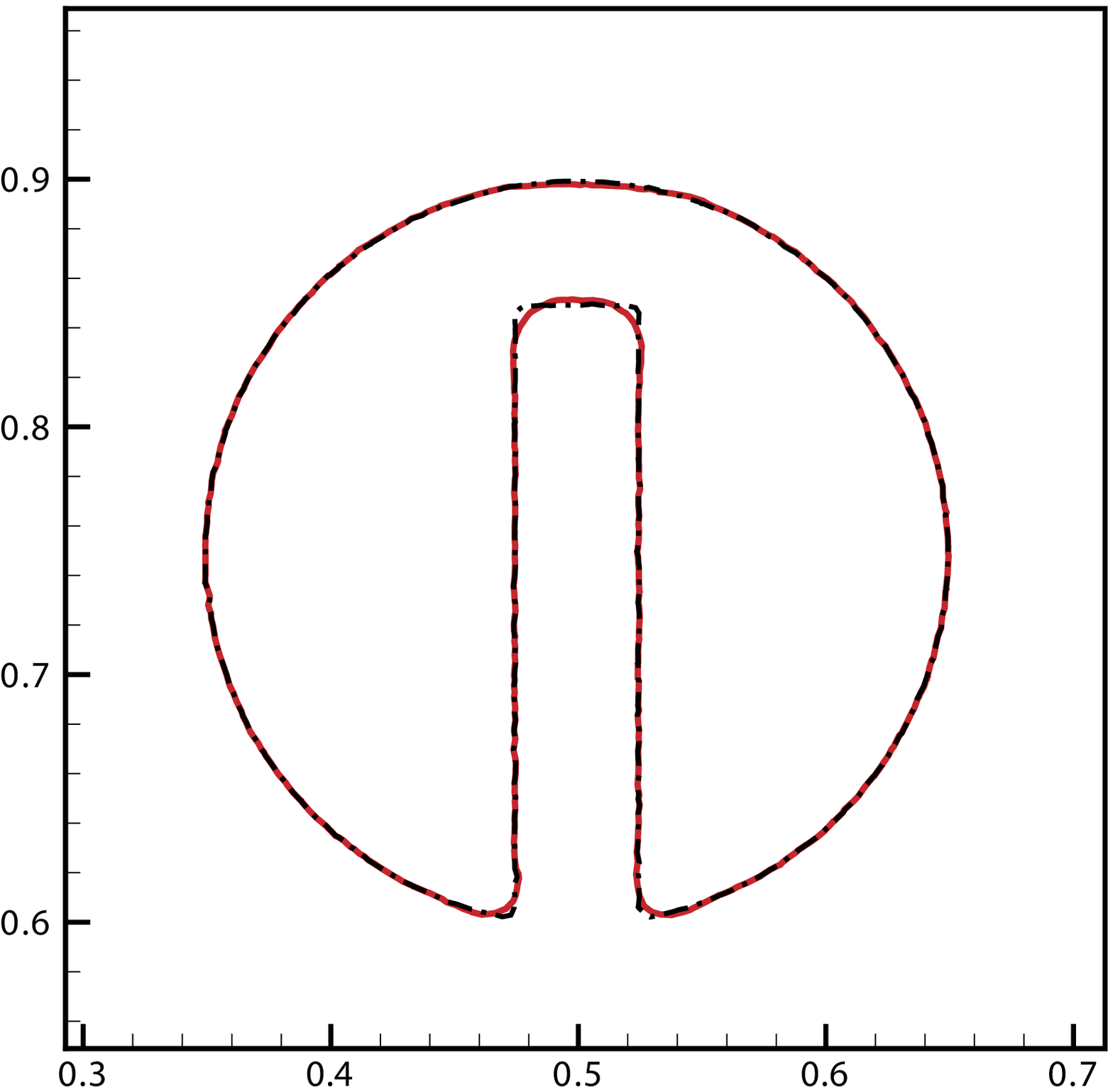}}
	\caption{\color{black}VOF 0.5 contour line in Zalesak solid body rotation test after one revolution on unstructured triangular meshes of different grid resolutions. The black dashed line stands for the exact solution, and the red solid line for the numerical solution.}
\label{zalesak_uns_results}
\end{figure}

Fig. \ref{zalesak_uns_results} presents the numerical results on unstructured triangular grid, and shows the similar solution quality to the Cartesian grid.  

The numerical errors and convergence rates are outlined in Table \ref{zalesak-uns-comparison}, which shows comparable or better numerical accuracy as compared to other THINC schemes.

\begin{table}[h]
\centering
\caption{Numerical errors ($E_r$) and convergence rates of Zalesak solid body rotation test after one revolution.} 
\begin{tabular}[t]{lSSSSS} \toprule
\multicolumn{6}{c}{Cartesian mesh}\\
\hline
{$\textbf{Methods}$}               & {$50$}              & {Order} & {$100$}               & {Order} & {$200$} \\ \midrule
{THINC-scaling}                    & {$6.87\times10^{-2}$} & {2.15}  & {$1.55\times10^{-2}$}   & {0.78}  & {$9.05\times10^{-3}$}\\
{MTHINC \cite{ii2012interface}}    & {$2.93\times10^{-2}$} & {0.86}  & {$1.61\times10^{-2}$}   & {1.03}  & {$7.91\times10^{-3}$}\\        
{UMTHINC \cite{xie2017toward}}     & {$8.12\times10^{-2}$} & {1.63}  & {$2.61\times10^{-2}$}   & {0.97}  & {$1.33\times10^{-2}$}\\  
{THINC/QQ \cite{xie2017toward}}    & {$8.96\times10^{-2}$} & {1.47}  & {$3.22\times10^{-2}$}   & {0.95}  & {$1.67\times10^{-2}$}\\\bottomrule
\multicolumn{6}{c}{Triangular mesh}\\
\hline
{$\textbf{Methods}$}               & {$50$}              & {Order} & {$100$}               & {Order} & {$200$} \\ \midrule
{THINC-scaling}                    & {\color{black}$6.81\times10^{-2}$} & {\color{black}1.40}  & {\color{black}$2.60\times10^{-2}$}   & {\color{black}1.00}  & {\color{black}$1.23\times10^{-2}$}\\
%{THINC/LS}                         & {$4.60\times10^{-2}$} & {1.29}  & {$1.88\times10^{-2}$}   & {1.01}  & {$9.35\times10^{-3}$}\\        
{THINC/QQ}                         & {$5.85\times10^{-2}$} & {1.11}  & {$2.70\times10^{-2}$}   & {0.81}  & {$1.54\times10^{-2}$}\\\bottomrule
\label{zalesak-uns-comparison}
\end{tabular}
\end{table}

\begin{figure}[htbp]
	\centering
	\subfigure[]{
		\centering
		\includegraphics[width=0.35\textwidth]{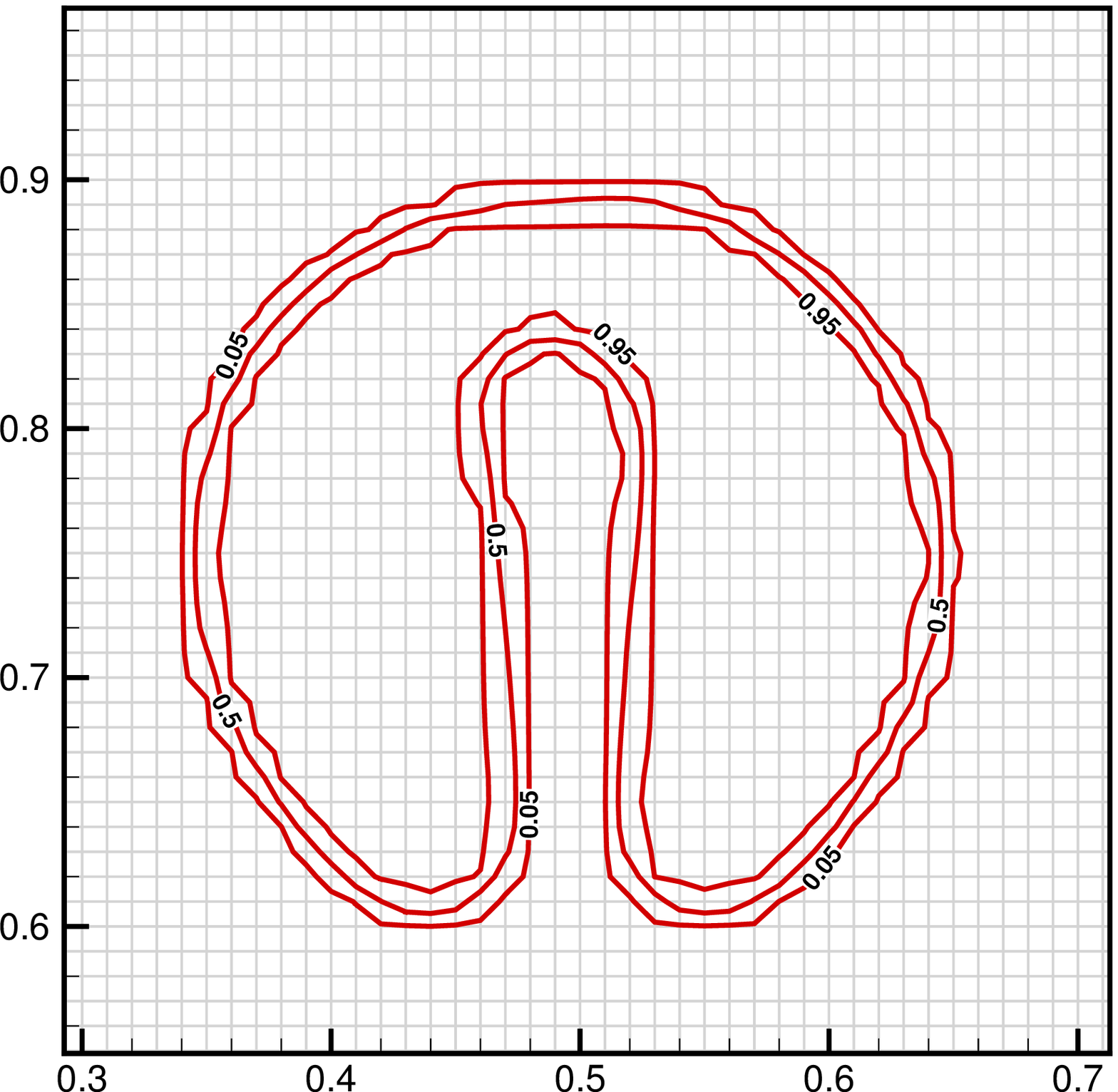}} \hspace{0.4cm}
	\subfigure[]{
		\centering
		\includegraphics[width=0.35\textwidth]{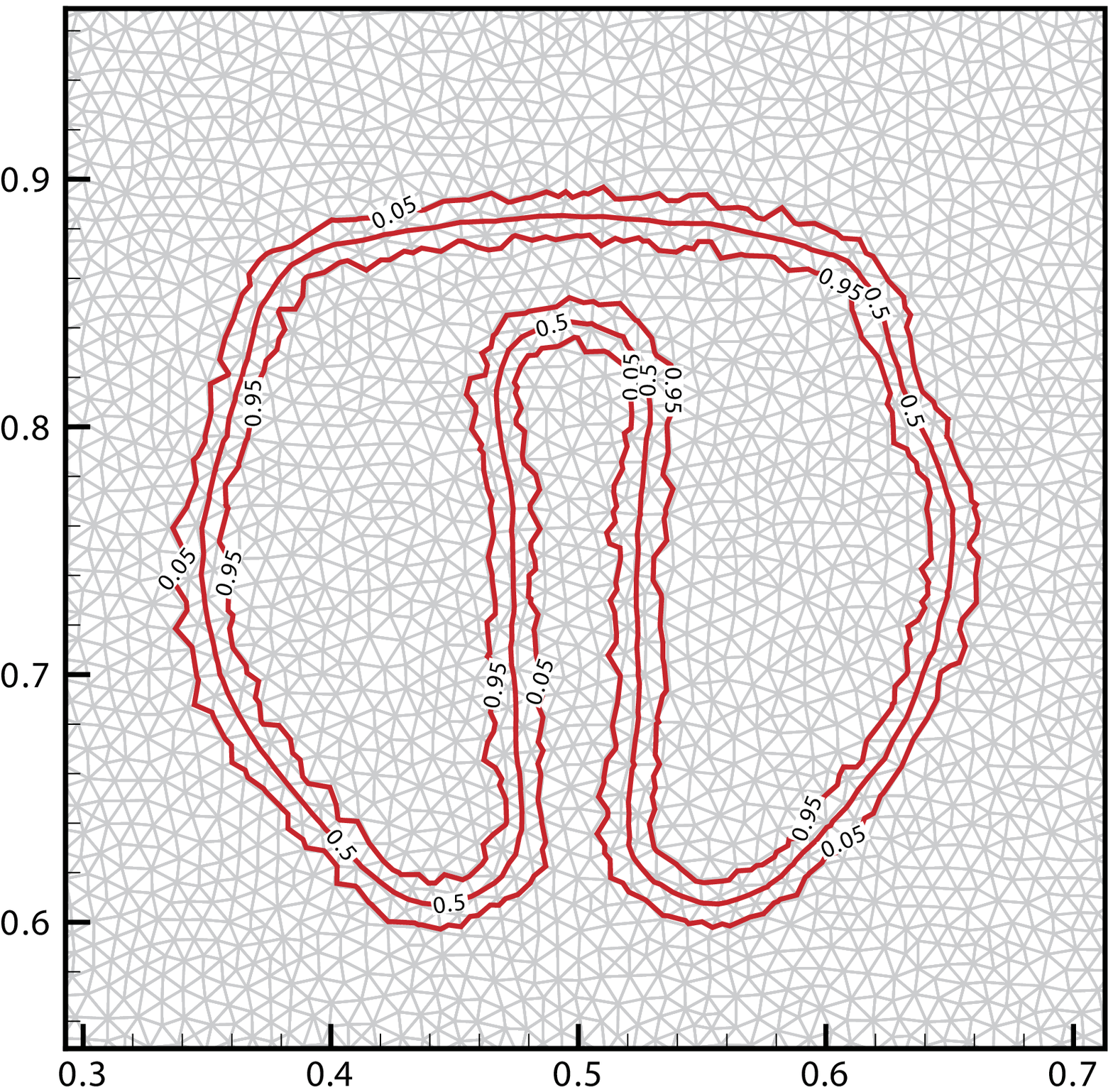}}
	\caption{\small{Numerical results showing $0.05$, $0.5$ and $0.95$ VOF contour lines in Zalesak solid body rotation test after ten revolutions on (a) structured grid, and {\color{black}(b) unstructured triangular grid with $\text{N}=100$ using THINC-scaling scheme.}}}
\label{diff-zalesak}
\end{figure}

We further verified the capability of the proposed scheme to maintain interface thickness and geometrical faithfulness in long-term computation. We plot the numerical results after ten revolutions in Fig.\ref{diff-zalesak}. It is observed that the interface thickness is preserved within 2-3 mesh cells, free from smearing-out even after ten revolutions. Moreover, the geometrical feather is well preserved in the numerical solution, which is very challenging for other existing VOF or level set methods. 

{\color{black} 
As discussed in Remark 7, the thickness of interface transition layer is  controlled by $\beta$ which is specified as $\beta=6/\Delta$ in all numerical tests for 2 and 3D presented in this paper. For a typical $100\times100$ mesh in 2D, for example, the actual value in the THINC function is 600, which is a value adequate to make the continuous THINC function mimic the Heaviside function. As $\beta\propto 1/\Delta $, refining grid resolution will reduce   
$\Delta$, thus steepen the transition jump and make the THINC function to converge the Heaviside function of VOF.  

}
%*****************************************************************************************************************
\subsubsection{Vortex deformation transport test}

The THINC-scaling scheme is further assessed by the single vortex test \cite{rider1998reconstructing}, known as Rider-Kothe shear flow test case, in which a circle initially centred at $\left(0.5,0.75\right)$ in an unit square domain is advected by time dependent velocity field given by the stream function as follows,

\begin{equation}
\Psi\left(x,y,t\right)=\frac{1}{\pi}{\sin}^2\left(\pi{x}\right){\sin}^2\left(\pi{y}\right){\cos}\left(\frac{\pi{t}}{T}\right),
\end{equation}
where $T=8$ is specified in this test. This test, as one of the most widely used benchmark tests, is more challenging to assess the capability of the scheme in capturing the heavily distorted interface with stretched tail when transported to $t=T/2$. From $t=T/2$ to $t=T$, the reverse velocity field restores the interface back to its initial shape. In case of $T=8$, the spiral tail becomes so thin that it can not be resolved by a coarse grid of finite resolution.

\begin{figure}[htbp]
    \centering
	\subfigure[$\text{N}=32$]{
		\centering
		\includegraphics[width=0.3\textwidth]{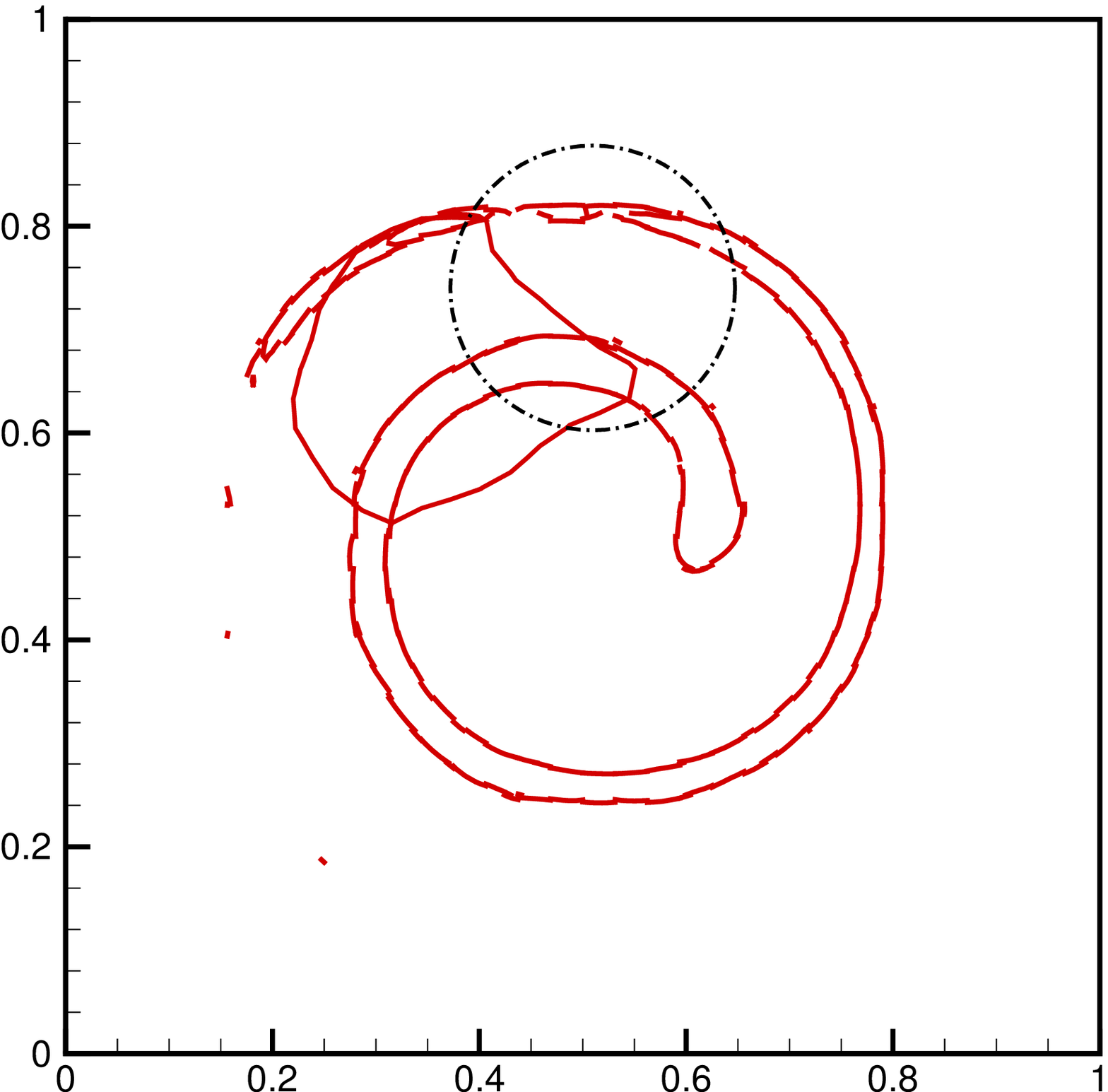}} \hspace{0.4cm}
	\centering
	\subfigure[$\text{N}=64$]{
		\centering
		\includegraphics[width=0.3\textwidth]{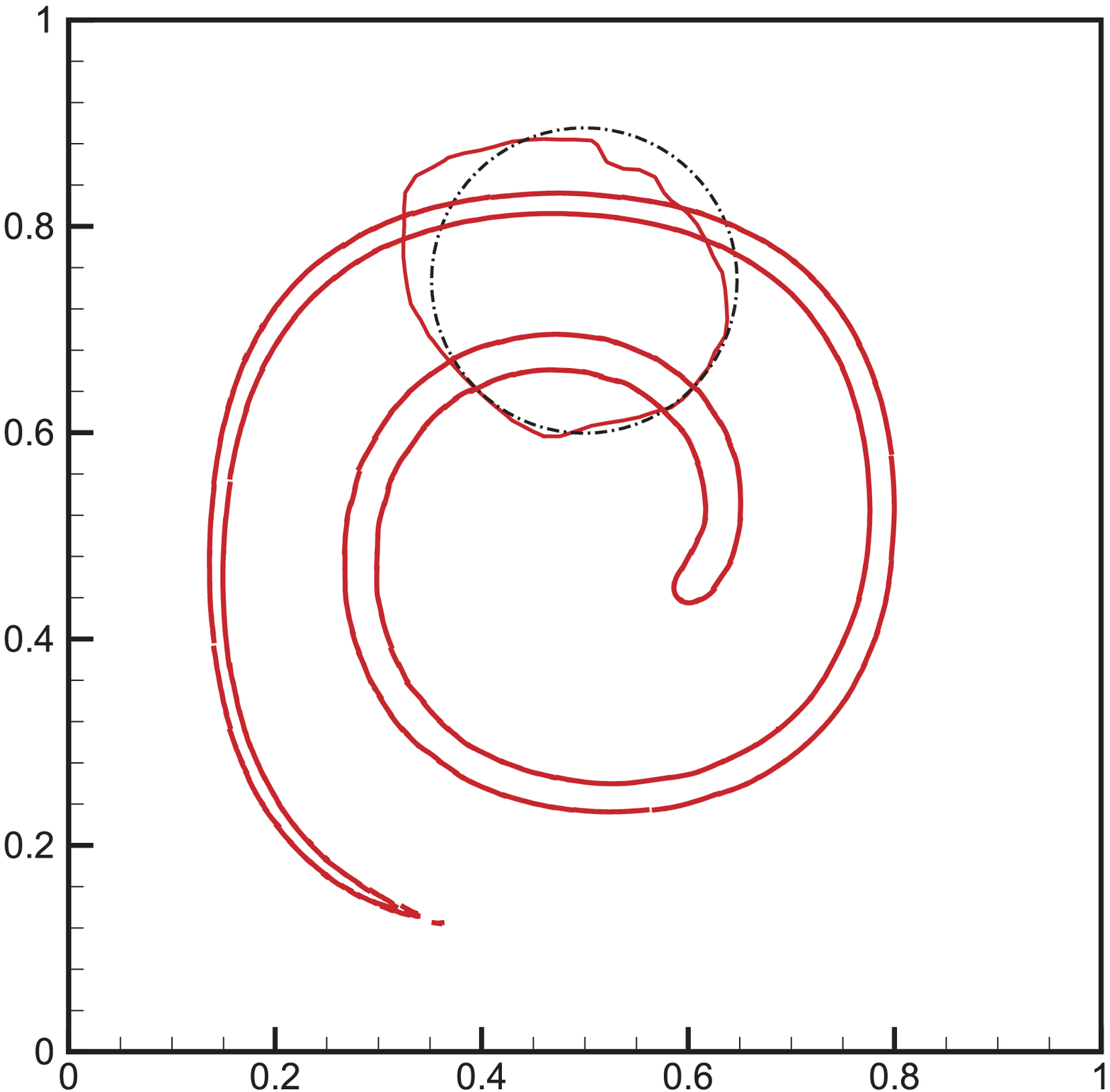}} \\
    \centering
	\subfigure[$\text{N}=128$]{
		\centering
		\includegraphics[width=0.3\textwidth]{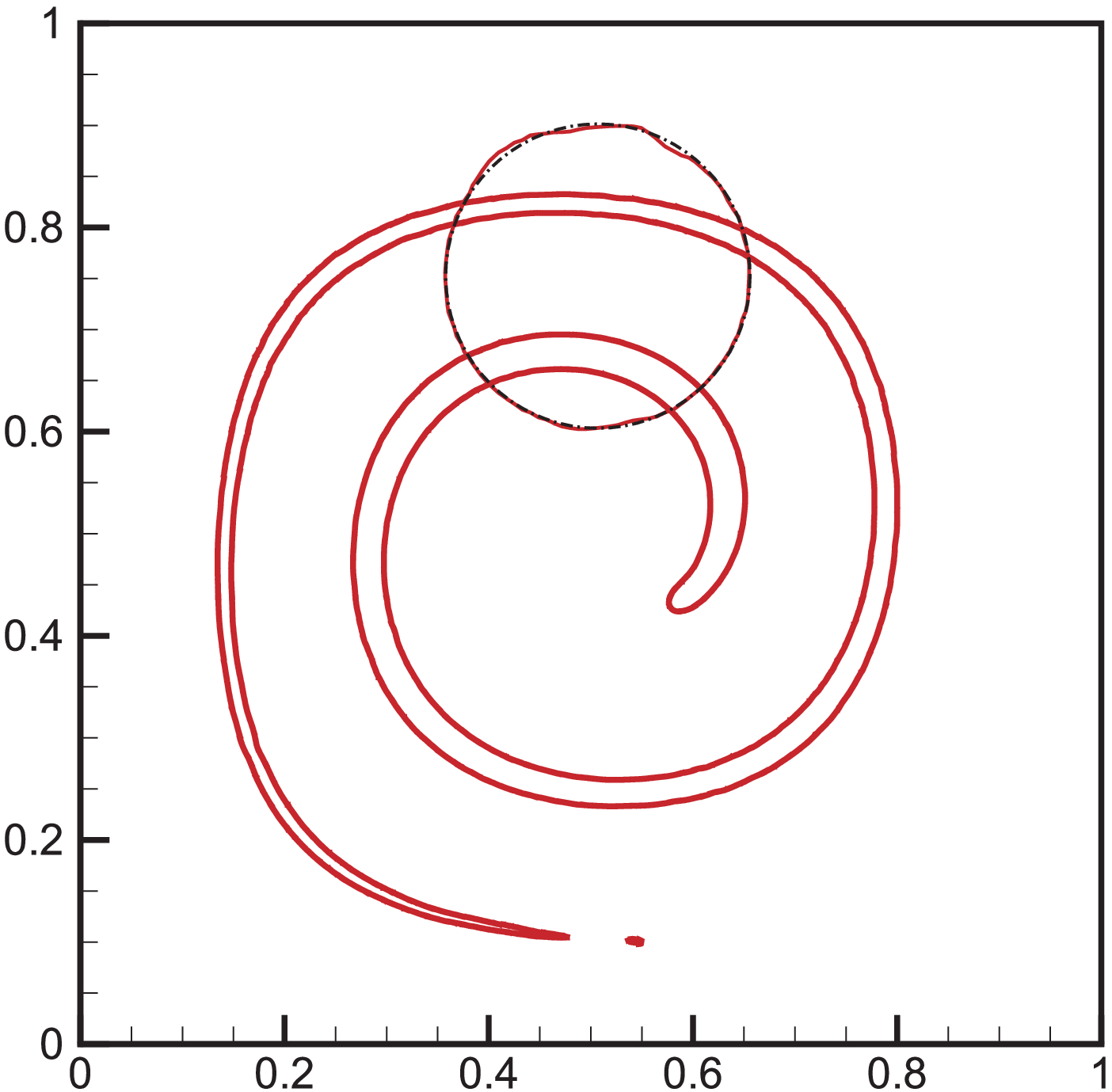}} \hspace{0.4cm}
	\subfigure[$\text{N}=256$]{
		\centering
		\includegraphics[width=0.3\textwidth]{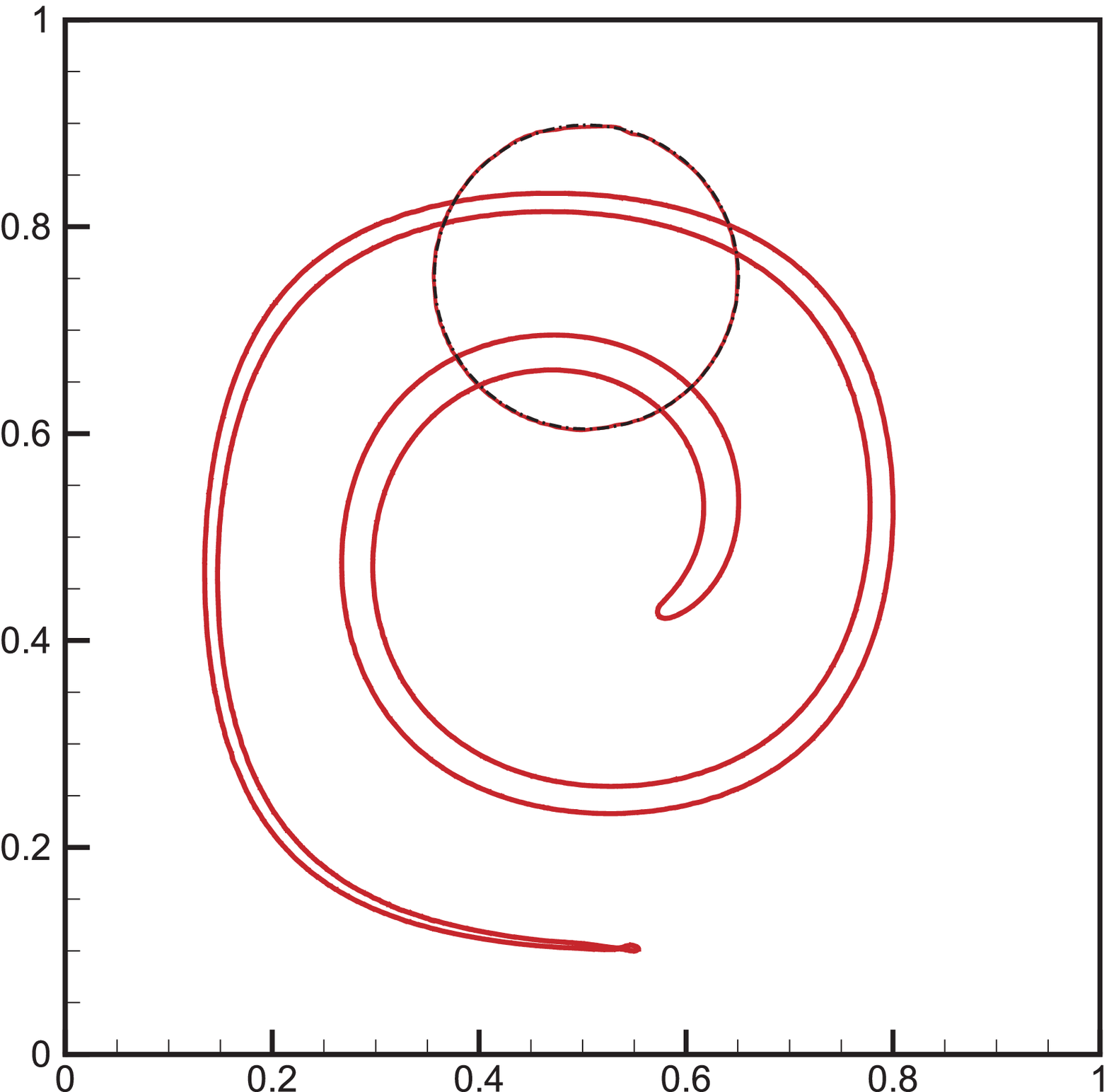}} 
	\caption{Numerical results for Rider-Kothe single vortex test showing the PSI for interface cells at $t=T/2$ and the VOF 0.5 contour lines at $t=T$ for structured grids of different resolutions.}
\label{rider_results}
\end{figure}

\par We tested this case on Cartesian grids with different resolutions. 
Numerical results on different grids at $t=T/2$ and $t=T$ are shown in Fig. \ref{rider_results}. The THINC-scaling scheme can capture the elongated tail even when the interface is under grid resolution, and can restore the initial circle with good solution quality. 

\begin{figure}[htbp]
	\centering
	\subfigure[] {
		\centering
		\includegraphics[width=0.4\textwidth]{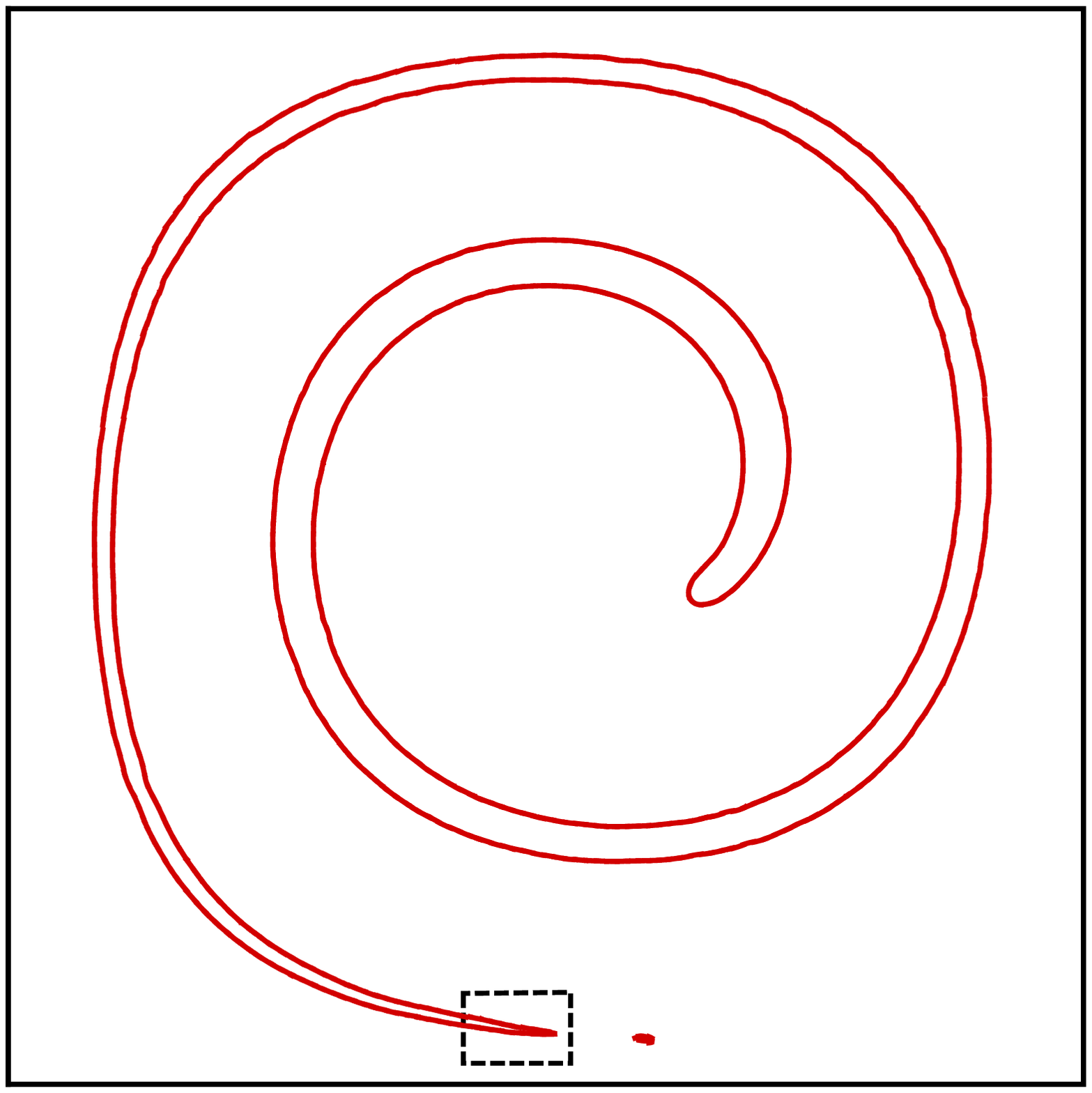}} \hspace{0.4cm}
	\centering
	\subfigure[] {
		\centering
		\includegraphics[width=0.4\textwidth]{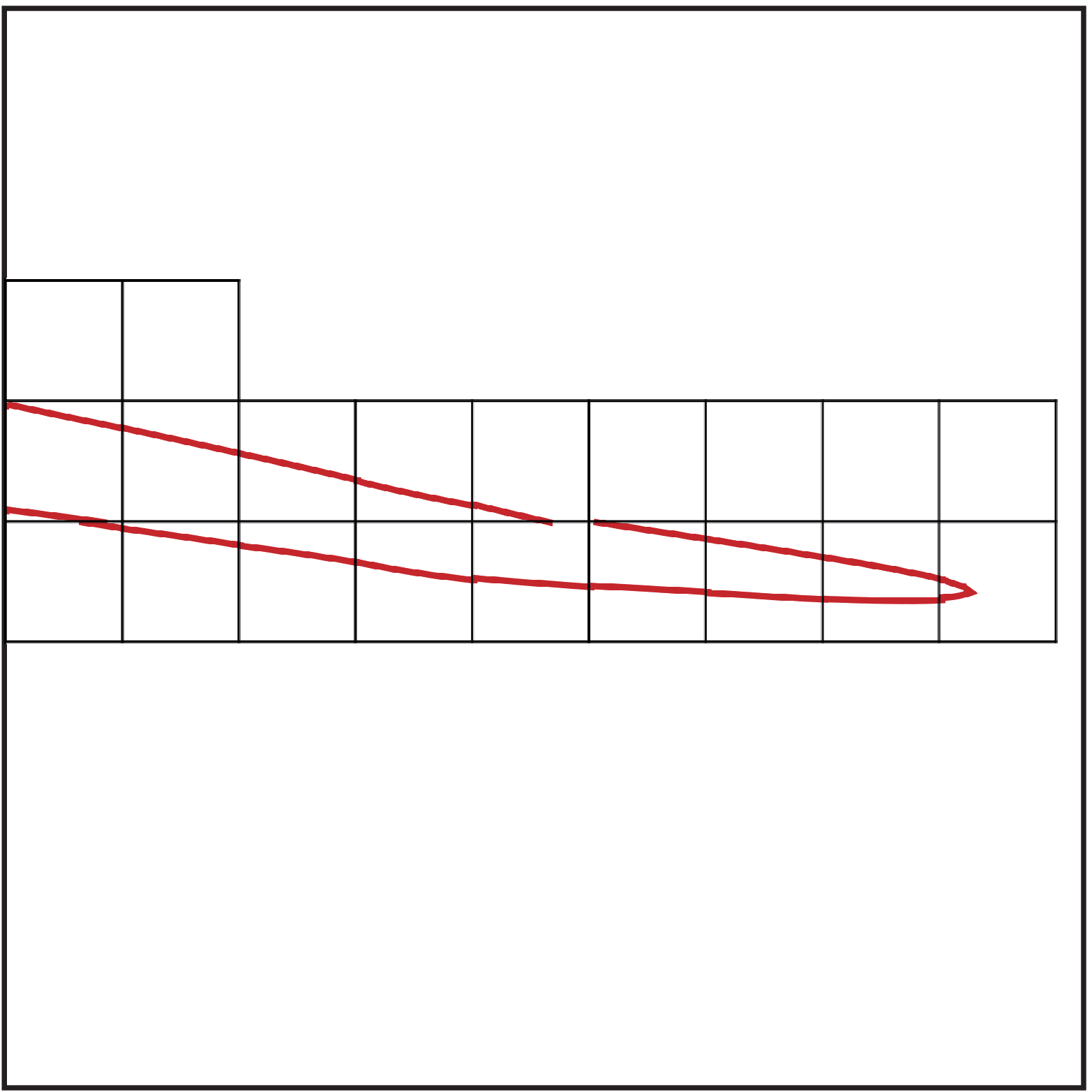}} \hspace{0.4cm}
	\caption{The distorted interface at $T/2$ on $128\times128$ grid (a) where the boxed part is enlarged as the close-up of the thin tail in panel (b) which shows that the film of tail is under the resolution of grid cell, but can be still resolved by the THINC-scaling scheme. }
\label{thin-tail}
\end{figure}

The PSI at $t=T/2$ on $128\times128$ grid is shown in Fig. \ref{thin-tail}. The tail tip is stretched into a thin film with a thickness smaller than the cell size. This sub-cell structure can still be reconstructed by the THINC-scaling scheme with quadratic or higher order polynomial representation. Consequently, the pieces of flotsam generated from the PLIC VOF methods are not observed here.  

A quantitative comparison with other VOF and hybrid methods are given in Table \ref{r-k-comparison}. It reveals the appealing accuracy of the present scheme. 

{\color{black} We should note that  among the tabulated errors the symmetric difference error used in \cite{Ahn2009sde,jemison2013MOF} provides a more rigorous metric to evaluate an explicitly reconstructed interface.   Since the THINC-scaling method simultaneously provides both the volume fraction $\bar{H}_i$ and the polynomial surface of the interface $\psi_i$, we are able to immediately assess the symmetric difference error. We define the symmetric difference error metric ($E_{sd}$)  as follows,  
\begin{equation}
E_{sd}=\sum_{c \in \mathcal{M}}\left|\left(T_{c} \cup R_{c}\right)-\left(T_{c} \cap R_{c}\right)\right|+\sum_{c \notin \mathcal{M}}\left|\bar{H}_{c}-\bar{H}_{c}^{e} \| \Omega_{c}\right|, 
\label{esd}
\end{equation}
where $\mathcal{M}$ denotes the set of interface cells identified by 
\begin{equation}
\mathcal{M}=\left\{\Omega_{i} \mid  \text { if }\left(\left(\varepsilon<\bar{H}_{i}<1-\varepsilon\right) \& \&\left(\varepsilon<\bar{H}_{i}^{e}<1-\varepsilon\right)\right)\right\}. 
\end{equation}
$T_{c}=T \cap  \Omega_{c}$ is the intersection of cell $\Omega_{c}$ with the true area ($T$) of the target material defined by $\mathcal{P}^e_{c}(\mathbf{x})\geq0$
that corresponds to the region encompassed by the zero contour of  exact level set field $\phi^e_{c}(\mathbf{x})=0$ . $R_{c}=R \cap  \Omega_{c}$ is the intersection area of cell $\Omega_{c}$ with the region ($R$) defined by  reconstructed interface polynomial
\begin{equation}
\psi_{c}(\mathbf{x}) \equiv \mathcal{P}_{c}(\mathbf{x})+\phi_{c}^{\Delta}\geq0. 
\end{equation}
Given $\phi^e_{c}(\mathbf{x})$ and $\psi_{c}(\mathbf{x})$, the cell-wise intersection area is computed by checking if the values of  $\phi^e_{c}(\mathbf{x})$ or $\psi_{c}(\mathbf{x})$ are positive over a set of finer control volumes sub-divided within each interface cell.   

We include the  symmetric difference error of THINC-scaling scheme in Table \ref{r-k-comparison} as well. It is observed that the  symmetric difference errors are slightly larger than $E\left(L_1\right)$, but shows a similar convergence behavior over different grid resolutions.  
}

\begin{table}[ht]
\centering
\footnotesize
\caption{Numerical errors $E\left(L_1\right)$ and convergence rates for Rider-Kothe test on Cartesian grid. For those with superscript ``$*$", the errors are measured by the symmetric difference error metric \eqref{esd}, which is equivalent to that used in \cite{jemison2013MOF}.} 
\resizebox{\columnwidth}{!}{
\begin{tabular}[t]{lSSSSSSS} \toprule
{$\textbf{Methods}$}               & {$32$} & {Order} &  {$64$} & {Order} & {$128$} & {Order} & {$256$} \\ \midrule
{THINC-scaling}                    & {$9.25\times10^{-2}$} & {2.86} & {$1.27\times10^{-2}$} & {3.22}  & {$1.36\times10^{-3}$}   & {1.77}  & {$3.98\times10^{-4}$}\\
{THINC-scaling$^*$}                & {$9.29\times10^{-2}$} & {2.86} & {$1.28\times10^{-2}$} & {3.16}  & {$1.43\times10^{-3}$}   & {1.76}  & {$4.23\times10^{-4}$}\\
{isoAdvector-plicRDF \cite{Henning2019VOF}} & {-} & {-} & {$1.26\times10^{-2}$} & {2.27}  & {$2.61\times10^{-3}$}   & {2.19}  & {$5.71\times10^{-4}$}\\
{UFVFC-Swartz \cite{Maric2018VOF}} & {-} & {-} & {$5.74\times10^{-3}$} & {1.98}  & {$1.45\times10^{-3}$}   & {1.94}  & {$3.77\times10^{-4}$}\\   
{Owkes and Desjardins \cite{Owkes2014VOF}}  & {-} & {-} & {$7.58\times10^{-3}$} & {2.01}  & {$1.88\times10^{-3}$}   & {2.21}  & {$4.04\times10^{-4}$}\\
{Rider-Kothe/Puckett \cite{rider1998reconstructing}}  & {$4.78\times10^{-2}$} & {2.78} & {$6.96\times10^{-3}$} & {2.27}  & {$1.44\times10^{-3}$}   & {-}  & {-}\\
{Stream/Puckett \cite{harvie2000new}}  & {$3.72\times10^{-2}$} & {2.45} & {$6.79\times10^{-3}$} & {2.52}  & {$1.18\times10^{-3}$}   & {-}  & {-}\\
{Stream/Youngs \cite{harvie2000new}}  & {$3.61\times10^{-2}$} & {1.85} & {$1.00\times10^{-2}$} & {2.21}  & {$2.16\times10^{-3}$}   & {-}  & {-}\\
{EMFPA/Puckett \cite{lopez2004volume}}  & {$3.77\times10^{-2}$} & {2.52} & {$6.58\times10^{-3}$} & {2.62}  & {$1.07\times10^{-3}$}   & {-}  & {-}\\
{CVTNA-PCFSC \cite{Liovic2006PLIC}}  & {$2.34\times10^{-3}$} & {2.12} & {$5.38\times10^{-4}$} & {2.03}  & {$1.31\times10^{-4}$}   & {-}  & {-}\\
{Hybrid markers-VOF \cite{aulisa2003mixed}}  & {$2.53\times10^{-2}$} & {3.19} & {$2.78\times10^{-3}$} & {2.54}  & {$4.78\times10^{-4}$}   & {-}  & {-}\\
{Markers-VOF \cite{lopez2005improved}}  & {$7.41\times10^{-3}$} & {1.83} & {$2.12\times10^{-3}$} & {2.31}  & {$4.27\times10^{-4}$}   & {-}  & {-}\\
{DS-MOF$$ \cite{jemison2013MOF}}  & {$3.45\times10^{-2}$} & {1.78} & {$1.00\times10^{-2}$} & {3.17}  & {$1.11\times10^{-3}$}   & {-}  & {-}\\
{DS-CLSVOF$$ \cite{jemison2013MOF}}  & {$5.45\times10^{-2}$} & {2.37} & {$1.05\times10^{-2}$} & {2.59}  & {$1.74\times10^{-3}$}   & {-}  & {-}\\
{DS-CLSMOF$$ \cite{jemison2013MOF}}  & {$2.92\times10^{-2}$} & {2.40} & {$5.51\times10^{-3}$} & {2.00}  & {$1.37\times10^{-3}$}   & {-}  & {-}\\
{THINC/QQ} \cite{xie2017toward}   & {$6.70\times10^{-2}$} & {1.98} & {$1.52\times10^{-2}$} & {2.33}  & {$3.06\times10^{-3}$}   & {-}  & {-}\\
{THINC/SW} \cite{xiao2011revisit}   & {$3.90\times10^{-2}$} & {1.36} & {$1.52\times10^{-2}$} & {1.94}  & {$3.96\times10^{-3}$}   & {-}  & {-}\\
{THINC/WLIC} \cite{yokoi2007efficient}   & {$4.16\times10^{-2}$} & {1.37} & {$1.61\times10^{-2}$} & {2.18}  & {$3.56\times10^{-3}$}   & {-}  & {-}
\\\bottomrule
\label{r-k-comparison}
\end{tabular}
}
\end{table}

 We solved the Rider-Kothe shear flow test case on unstructured grid with triangular cells. In order to capture fine interface structures, the choice of $\beta$ is same as discussed in the previous test case. We experimented with two cases using respectively 6 (GP=6) and 12 (GP=12) quadrature points to calculate integration \eqref{mass_constraint} or \eqref{THINCf-integration}. 
 
\begin{figure}[htbp]
	\centering
	\subfigure[$\text{N}=32$] {
		\centering
		\includegraphics[width=0.4\textwidth]{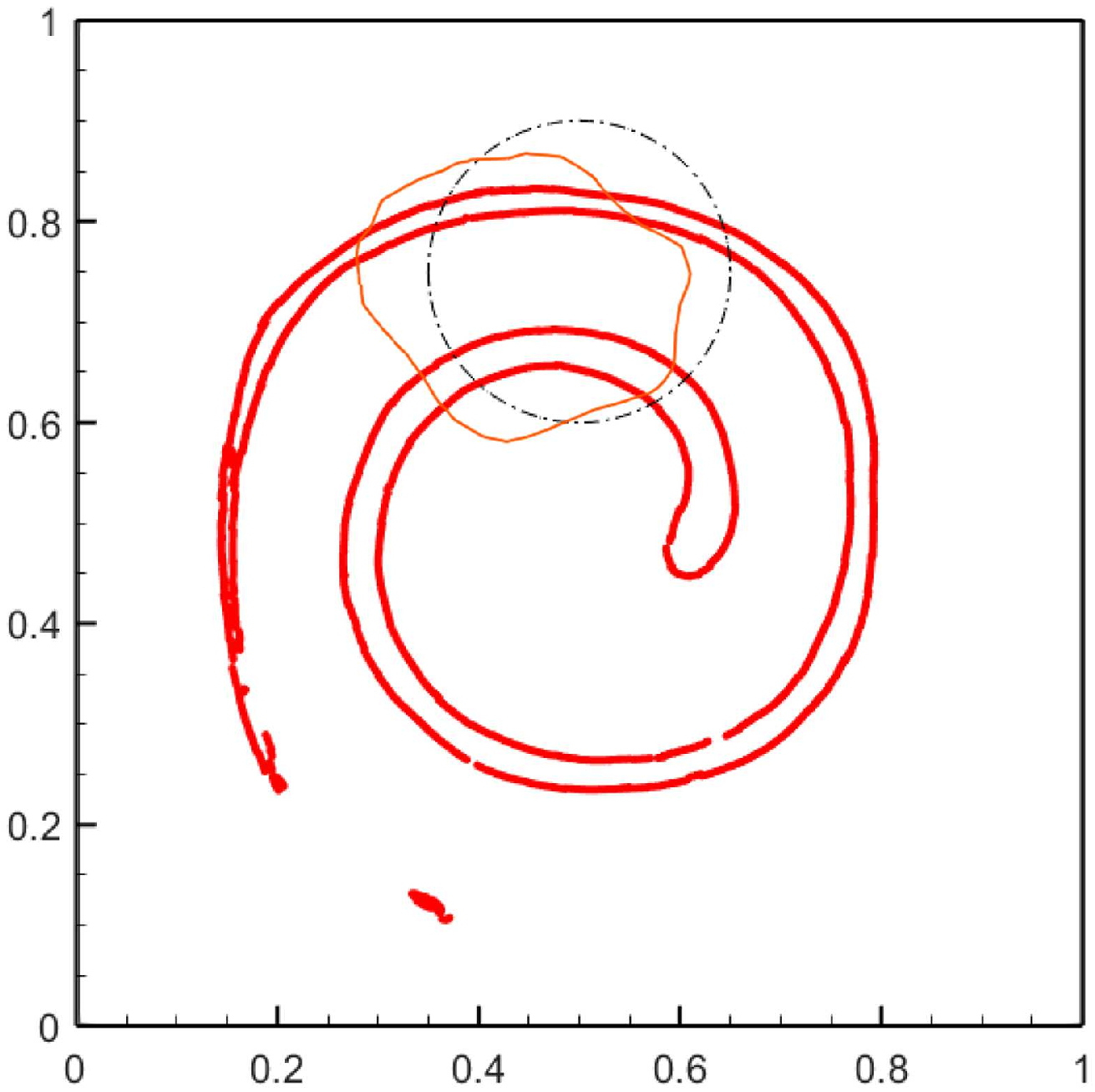}}\hspace{0.4cm}
	\centering
	\subfigure[$\text{N}=64$] {
		\centering
		\includegraphics[width=0.4\textwidth]{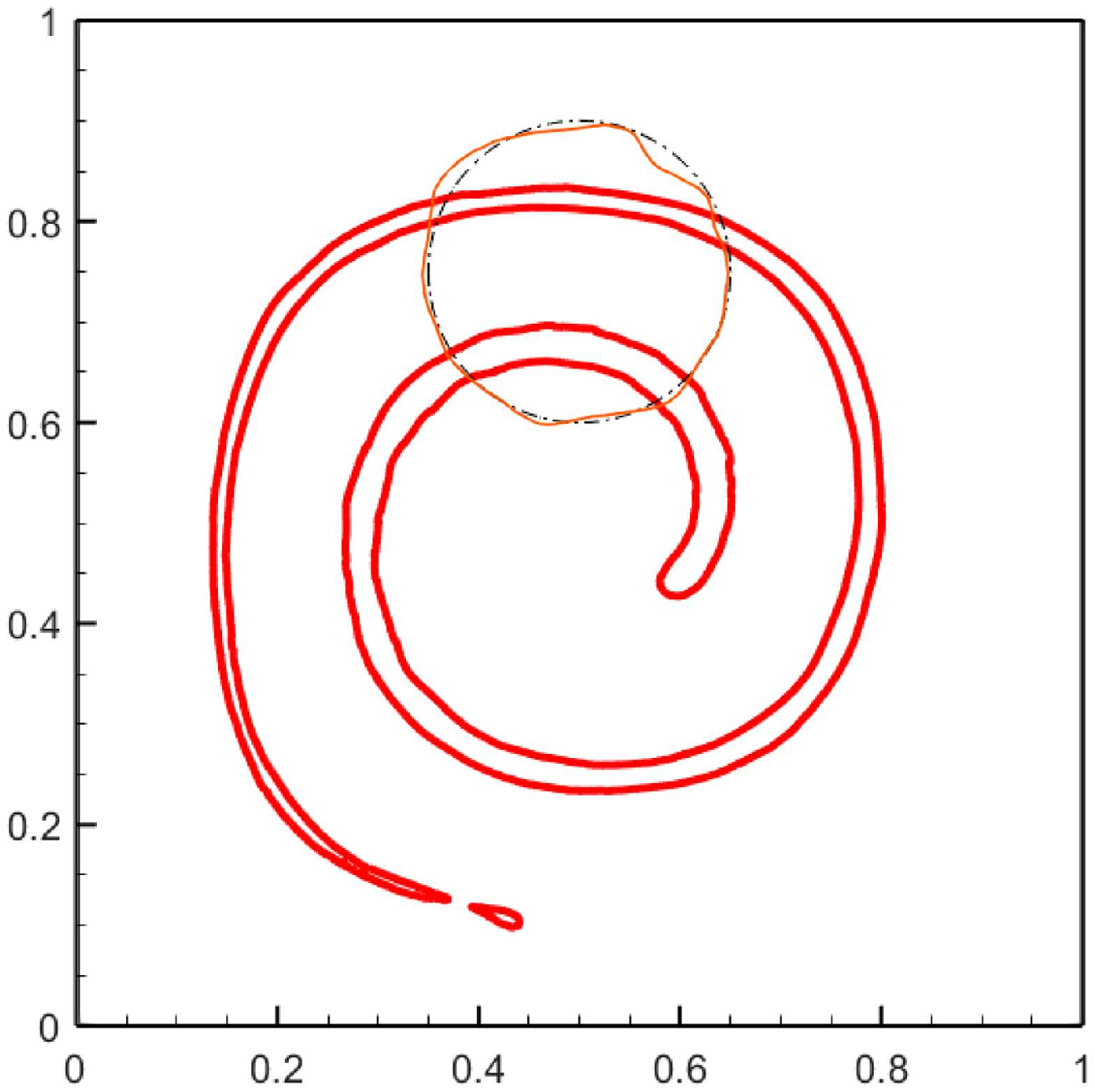}}
	\centering
	\subfigure[$\text{N}=128$] {
		\centering
		\includegraphics[width=0.4\textwidth]{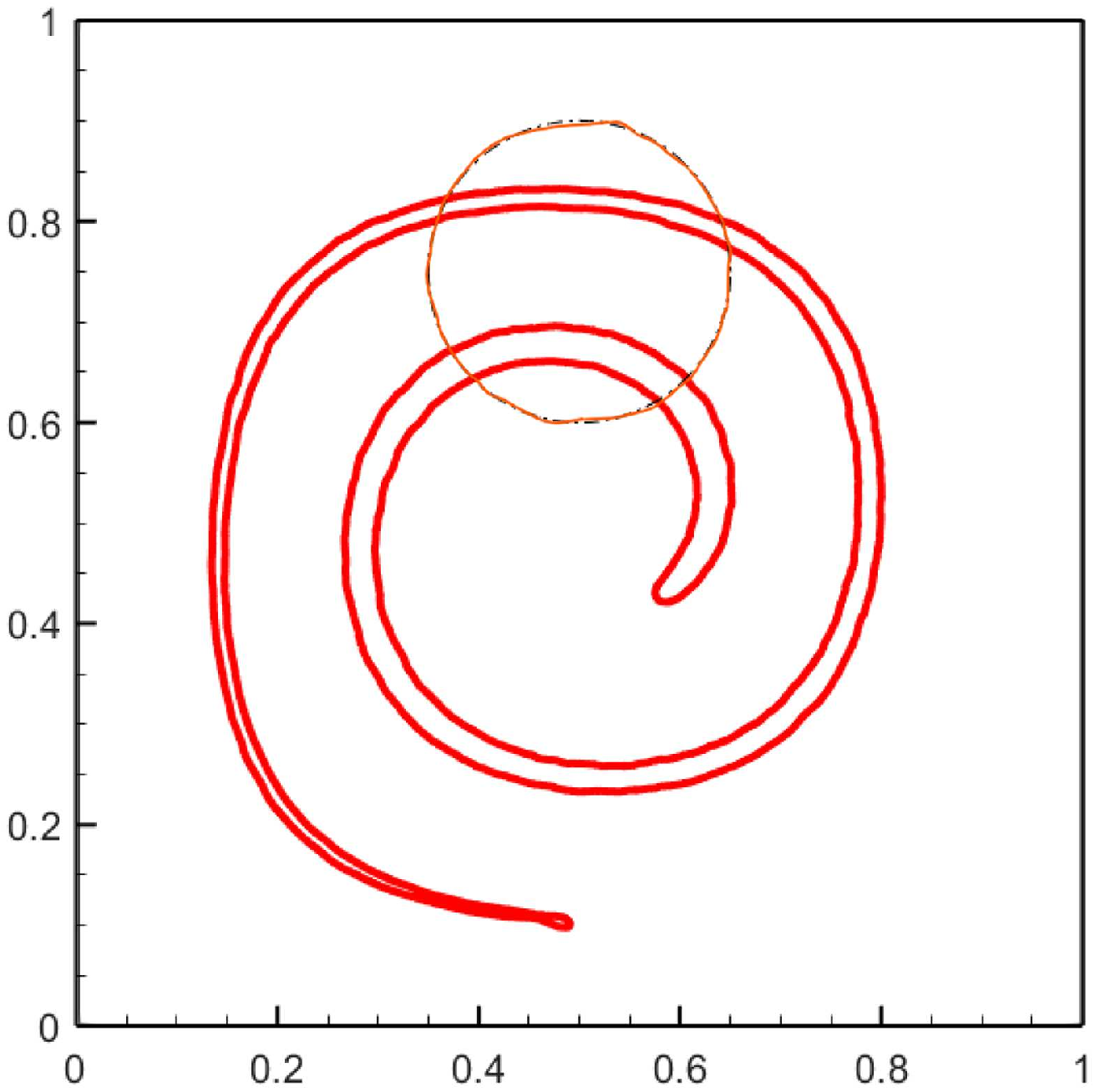}}\hspace{0.4cm}
	\centering
	\subfigure[$\text{N}=256$] {
		\centering
		\includegraphics[width=0.4\textwidth]{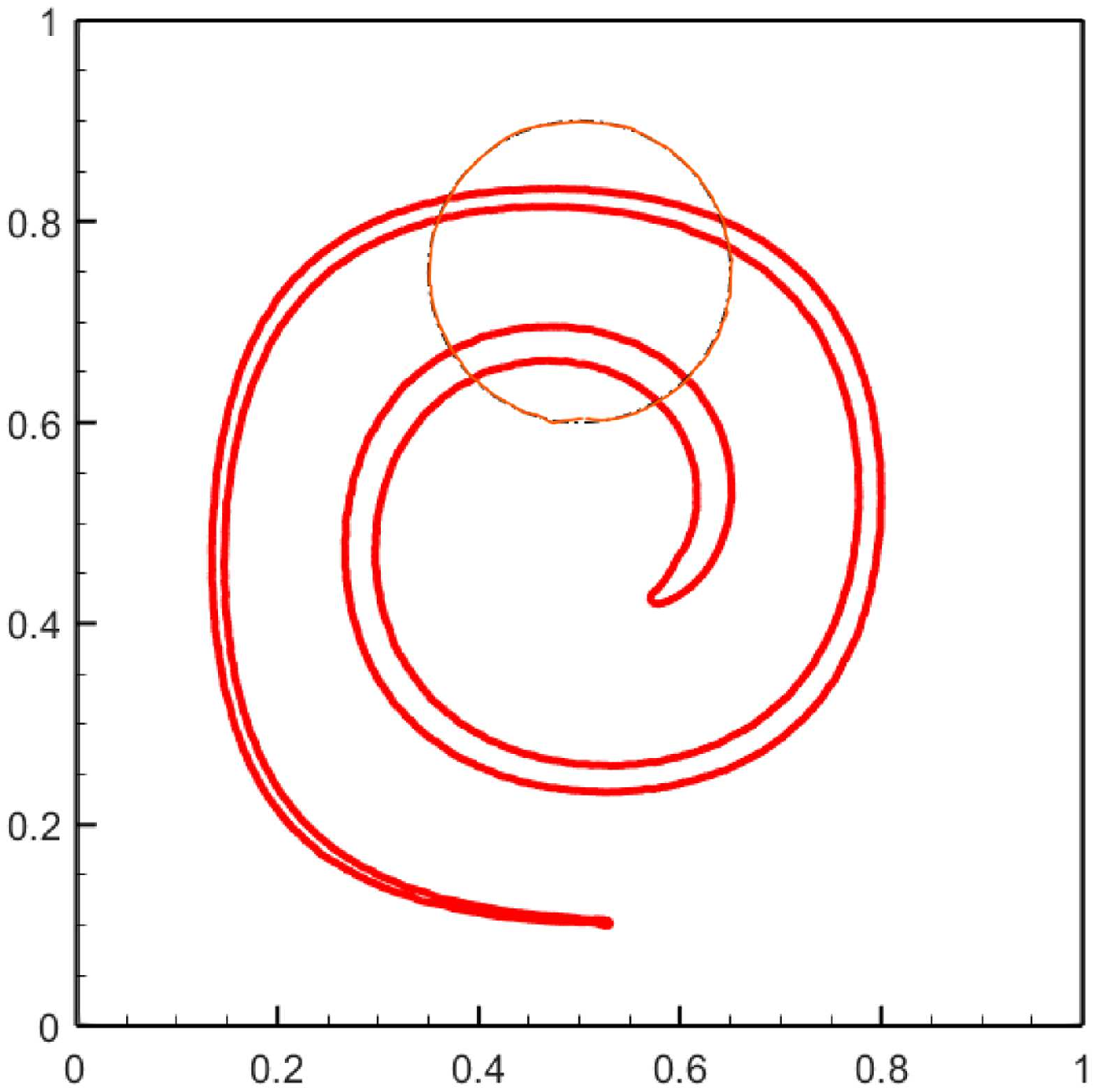}}
	\caption{{\color{black}Numerical results for Rider-Kothe single vortex test showing PSI at $t=T/2$ and VOF 0.5-contour at $t=T$ for triangular unstructured grids of different  resolutions. Area integration is done using $\text{GP}=6$ (6 Gaussian quadrature points).}}
\label{RK-uns-2}
\end{figure}

\begin{figure}[htbp]{
\begin{center}
\includegraphics[width=0.7\textwidth]{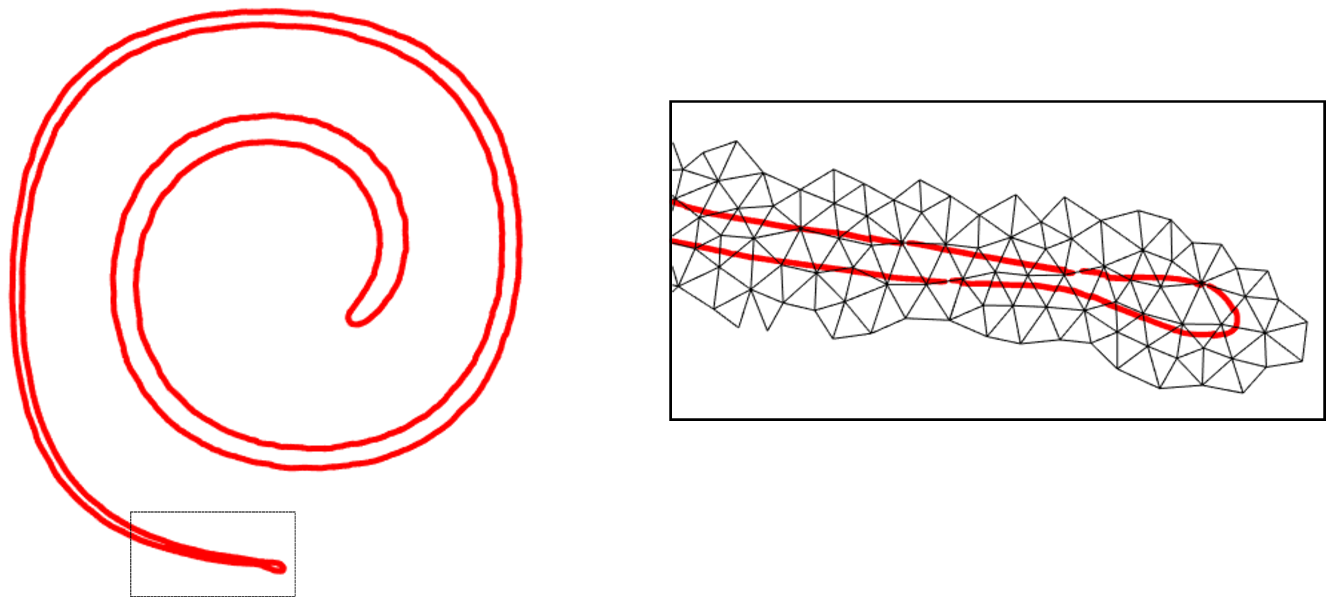}
\end{center}} 
\caption{\color{black}Enlarged view of PSI in interface cells highlighted in dashed box in left panel for $\text{N}=128$ and $\text{GP}=6$.}
\label{uns_PSI_tail2}
\end{figure}

 We plot the numerical results in Fig.\ref{RK-uns-2}. The PSIs in the interface cells accurately present the reconstructed interface. The flotsams in the PLIC reconstruction \cite{Henning2019VOF,jofre2014} are avoided in the present results. As observed from the enlarged view of thin tail part in Fig.\ref{uns_PSI_tail2}, our scheme is capable of  retrieving interface structures which are even under grid resolution. Similar to the results on structured grid, we can observe that there are few cells in which two interface segments are identified by using quadratic PSIs.  
 
\begin{figure}[htbp]
	\centering
	\subfigure[$\text{N}=32$] {
		\centering
		\includegraphics[width=0.4\textwidth]{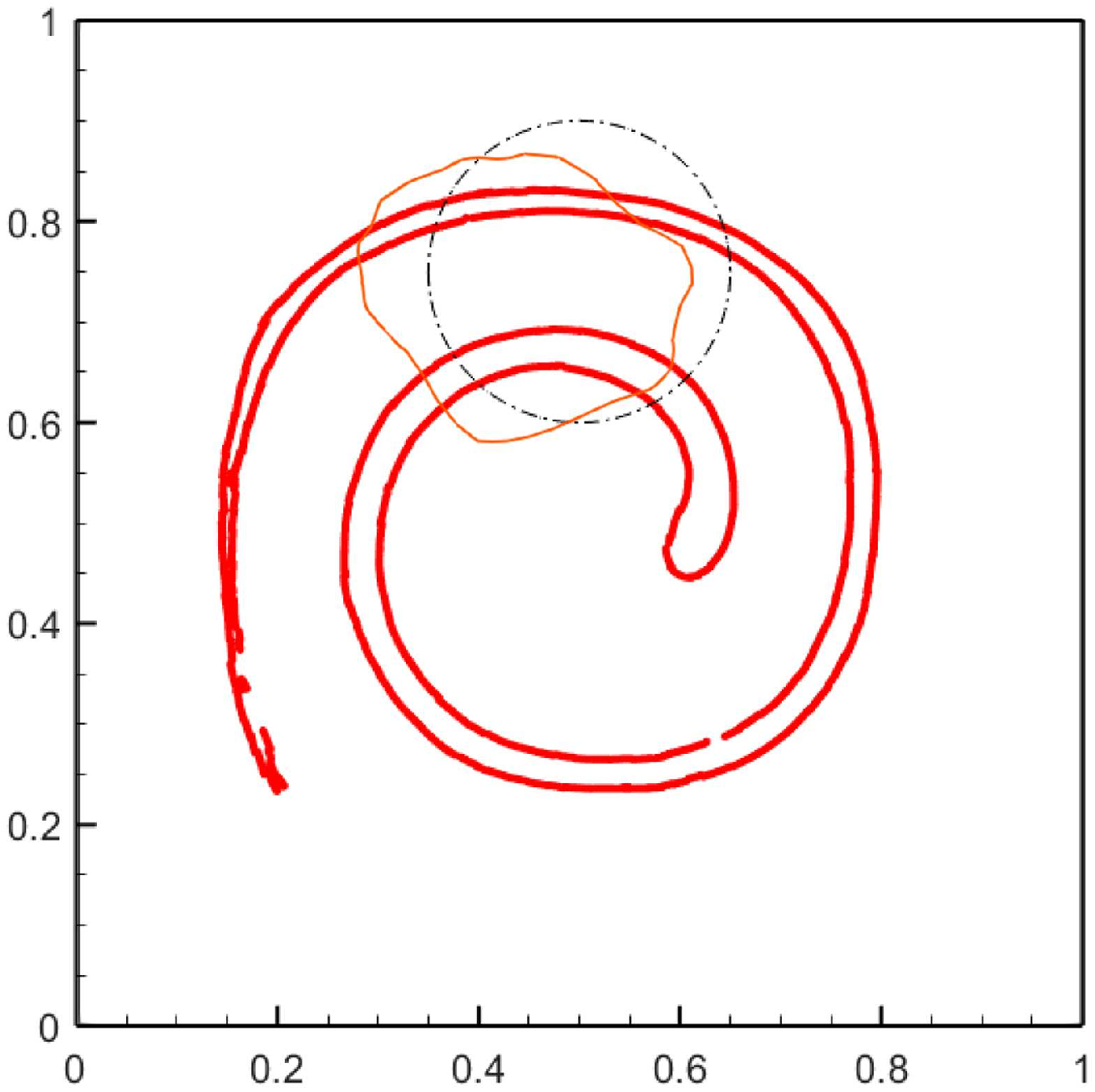}}\hspace{0.4cm}
	\centering
	\subfigure[$\text{N}=64$] {
		\centering
		\includegraphics[width=0.4\textwidth]{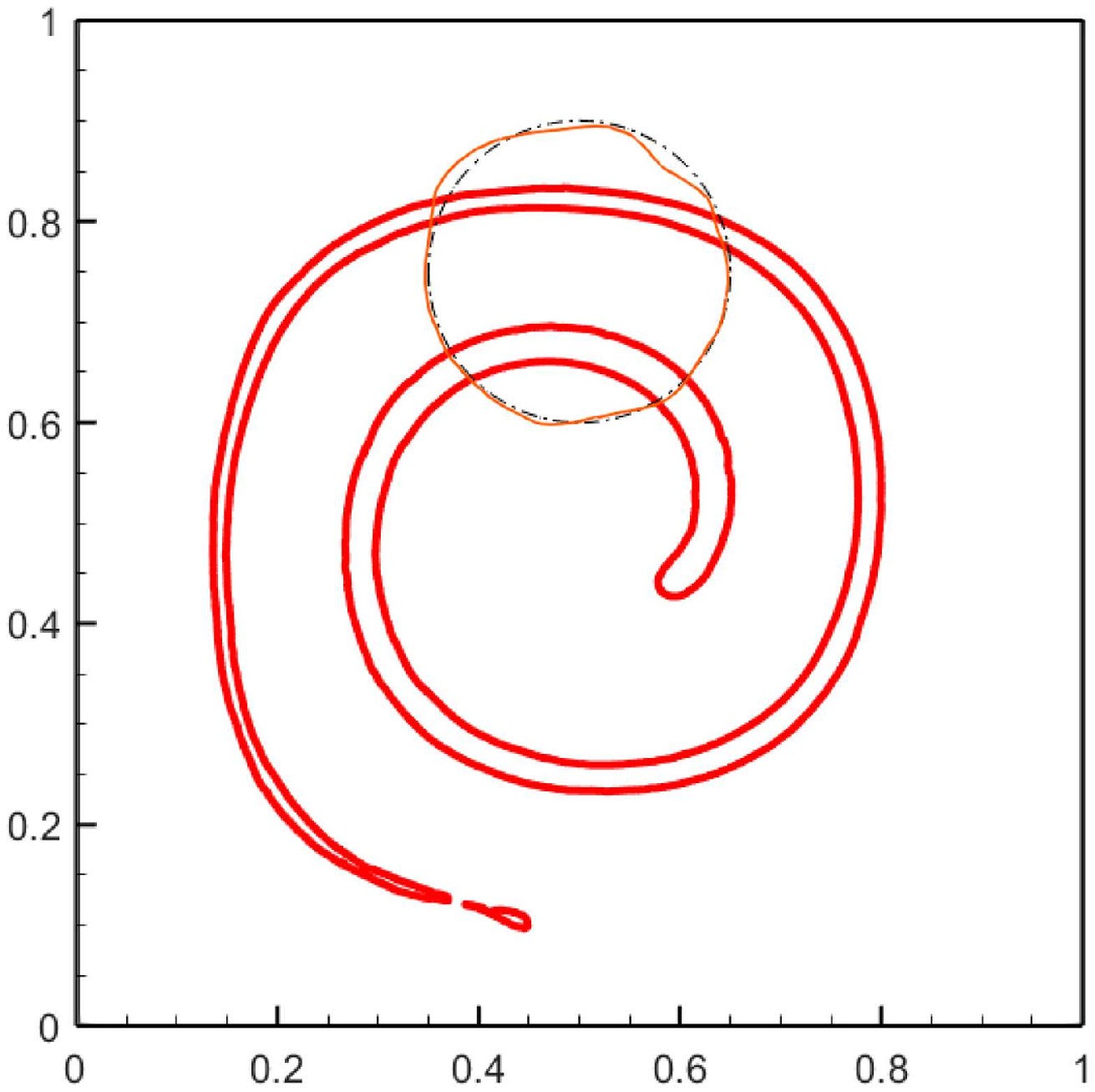}}
	\centering
	\subfigure[$\text{N}=128$] {
		\centering
		\includegraphics[width=0.4\textwidth]{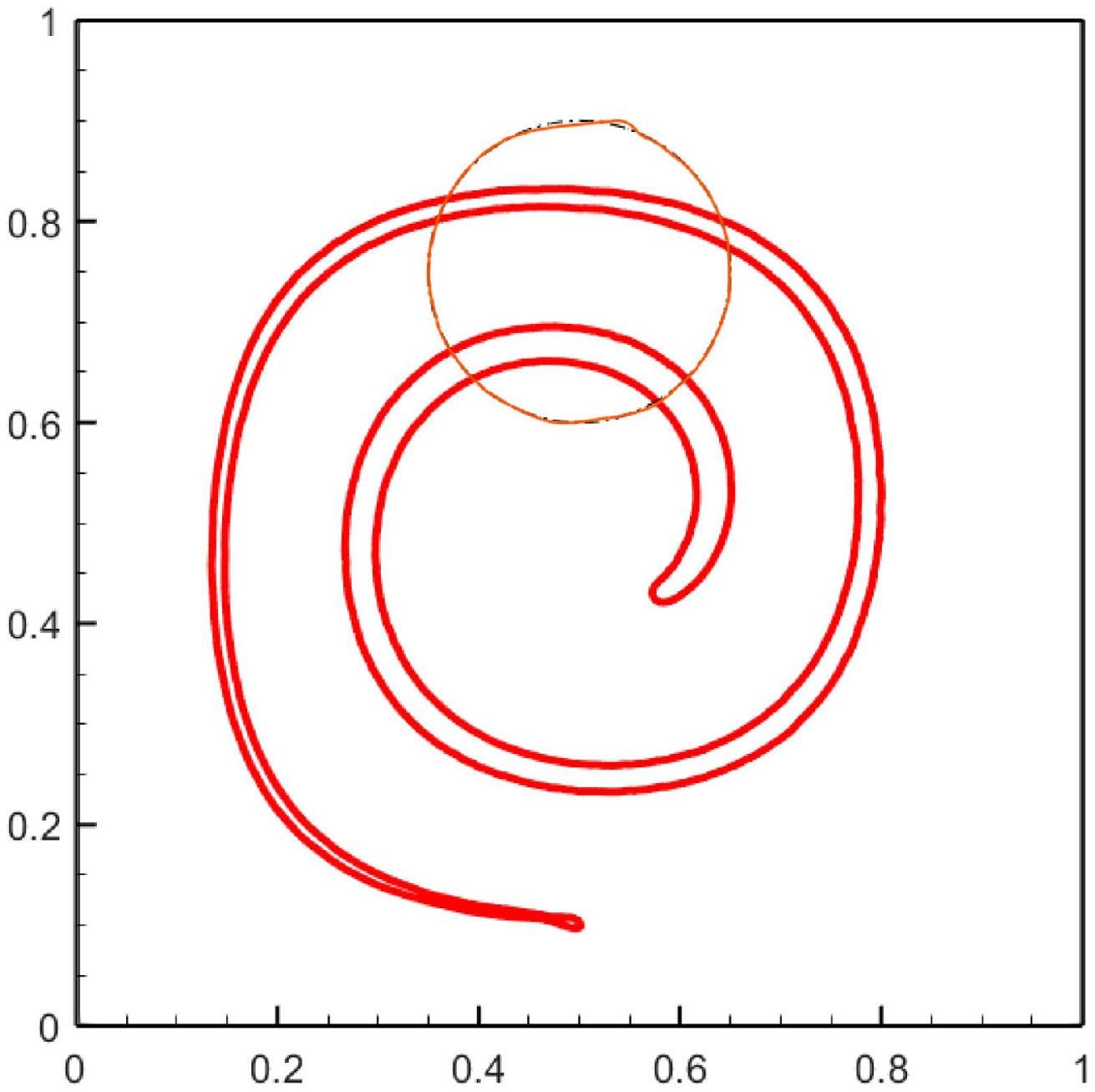}}\hspace{0.4cm}
	\centering
	\subfigure[$\text{N}=256$] {
		\centering
		\includegraphics[width=0.4\textwidth]{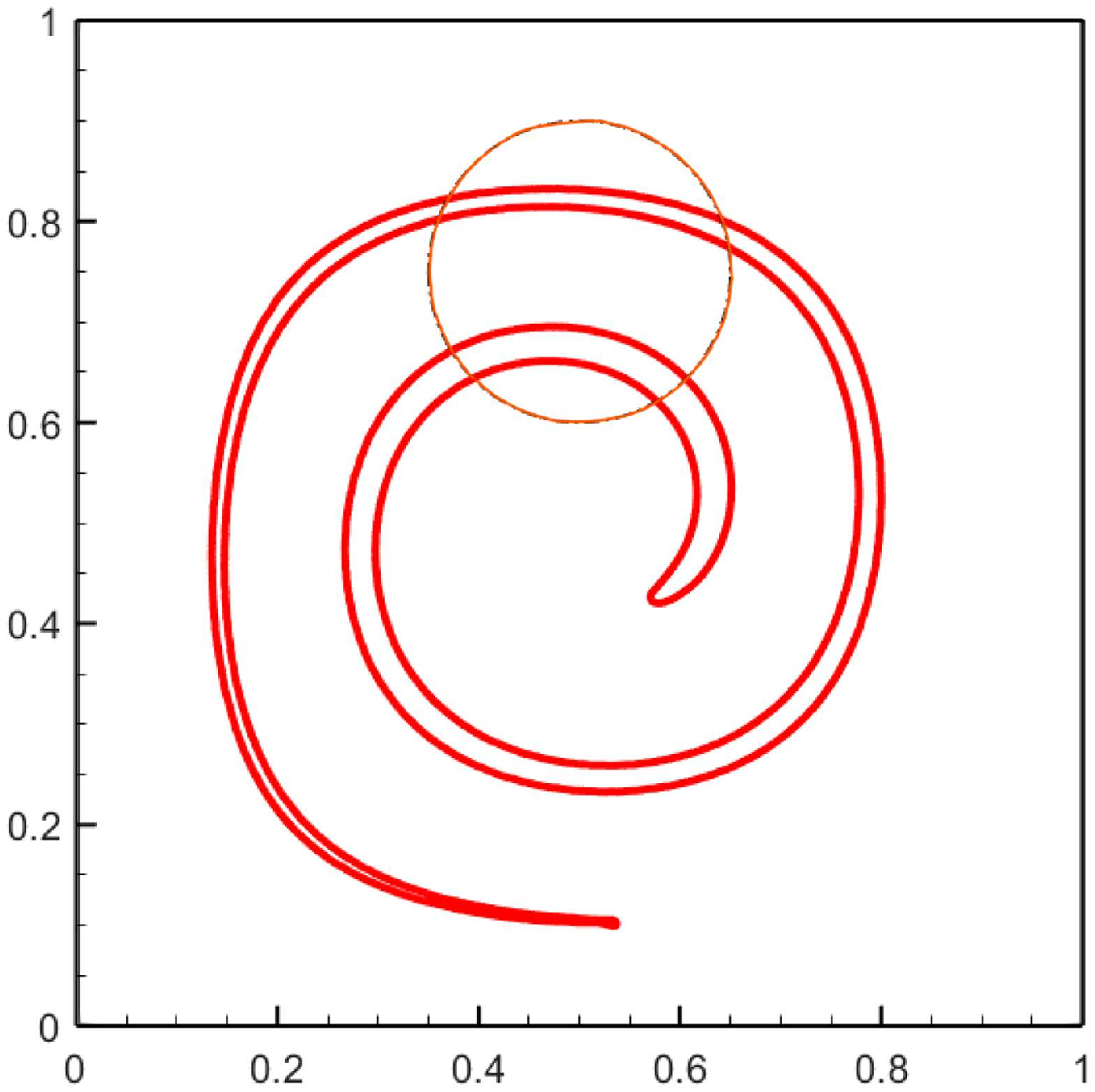}}
	\caption{\color{black}Same as Fig.\ref{RK-uns-2}, but using 12 Gaussian quadrature points for area integration.}
\label{RK-uns}
\end{figure}

\begin{figure}[htbp]{
\begin{center}
\includegraphics[width=0.7\textwidth]{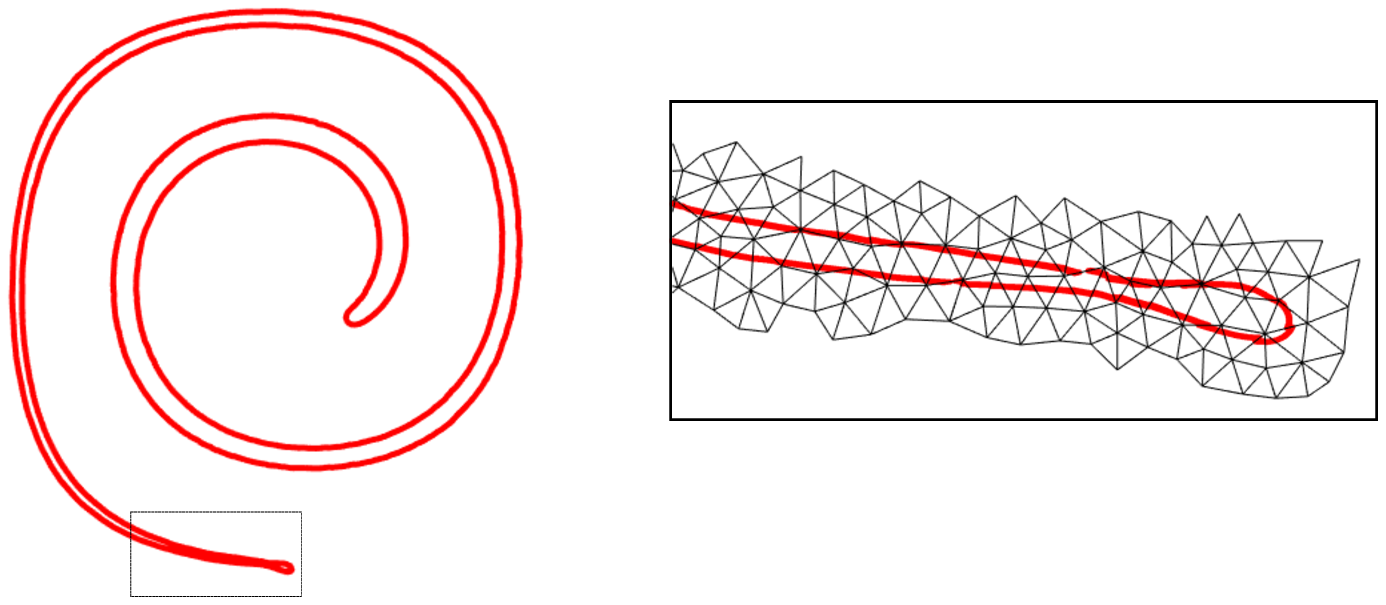}
\end{center}} 
\caption{\color{black}Same as Fig.\ref{uns_PSI_tail2}, but using 12 Gaussian quadrature points for area integration.}
\label{uns_PSI_tail}
\end{figure}

 In order to examine the sensitiveness of number of Gaussian quadrature points for area integration on 2D triangular grid used for solving \eqref{THINCf-integration}, we also show the numerical results using 12 Gaussian points in Fig.\ref{RK-uns}  and Fig.\ref{uns_PSI_tail}. Although using more Gaussian quadrature points might improve numerical accuracy, the improvement looks less significant. Thus, we suggest using 6 Gaussian quadrature points for area integration on triangular grid element in practice as a better trade-off between solution quality and computational cost. 
 
 A quantitative comparison is given in Table \ref{r-k-unstructured-comparison}. The THINC-scaling scheme shows superiority in accuracy on unstructured grid as well. {\color{black} We also include the  symmetric difference error of THINC-scaling scheme in Table \ref{r-k-unstructured-comparison}. Similar to the structured grid case, the  symmetric difference errors are slight larger than $E\left(L_1\right)$, but shows a similar convergence behavior over different grid resolutions.  
}

\begin{table}[ht]
\centering
\footnotesize
\caption{Numerical errors $E\left(L_1\right)$ and convergence rates for Rider-Kothe test on unstructured triangular grid.  For those with superscript ``$*$", the symmetric difference error metric  \eqref{esd} is used.} 
\resizebox{\columnwidth}{!}{
\begin{tabular}[t]{lSSSSSSS} \toprule
{$\textbf{Methods}$}               & {$32$}                & {Order} & {$64$}                  & {Order} & {$128$}  &  {Order}  & {256}\\ \midrule
{THINC-scaling (GP=6)}             & {\color{black}$3.48\times10^{-2}$} & {\color{black}2.80}  & {\color{black}$4.99\times10^{-3}$}   & {\color{black}1.66}  & {\color{black}$1.58\times10^{-3}$} &  {\color{black}1.12}  & {\color{black}$7.29\times10^{-4}$}\\
{THINC-scaling$^*$ (GP=6)}             & {\color{black}$3.50\times10^{-2}$} & {\color{black}2.56}  & {\color{black}$5.92\times10^{-3}$}   & {\color{black}1.57}  & {\color{black}$2.00\times10^{-3}$} &  {\color{black}1.02}  & {\color{black}$9.85\times10^{-4}$}\\
{THINC-scaling (GP=12)}            & {\color{black}$3.51\times10^{-2}$} & {\color{black}2.91}  & {\color{black}$4.68\times10^{-3}$}   & {\color{black}1.88}  & {\color{black}$1.27\times10^{-3}$} &  {\color{black}1.44}  &  {\color{black}$4.69\times10^{-4}$}\\
{THINC/QQ (GP=6)}                          & {$4.02\times10^{-2}$} & {2.17}  & {$8.93\times10^{-3}$}   & {1.69}  & {$2.76\times10^{-3}$} & {1.05}      & {$1.33\times10^{-3}$}\\
{THINC/QQ (GP=12)}                         & {$4.01\times10^{-2}$} & {1.90}  & {$1.07\times10^{-2}$}   & {1.99}  & {$2.68\times10^{-3}$} & {1.01}   & {$1.33\times10^{-3}$}\\
{isoAdvector-plicRDF \cite{Henning2019VOF}} & {-} & {-}  & {$2.21\times10^{-2}$}   & {2.62}  & {$3.58\times10^{-3}$} & {2.25} & {$7.51\times10^{-4}$}\\
{Shahbazi/Paraschivoiu \cite{shahbazi2003VOF}} & {$3.55\times10^{-2}$} & {2.78}  & {$7.17\times10^{-3}$}   & {2.27}  & {$1.44\times10^{-3}$} & {-} & {-}
\\\bottomrule
\label{r-k-unstructured-comparison}
\end{tabular}
}
\end{table}

\begin{figure}[htbp]
	\centering
	\subfigure[$\text{N}=32$] {
		\centering
		\includegraphics[width=0.3\textwidth]{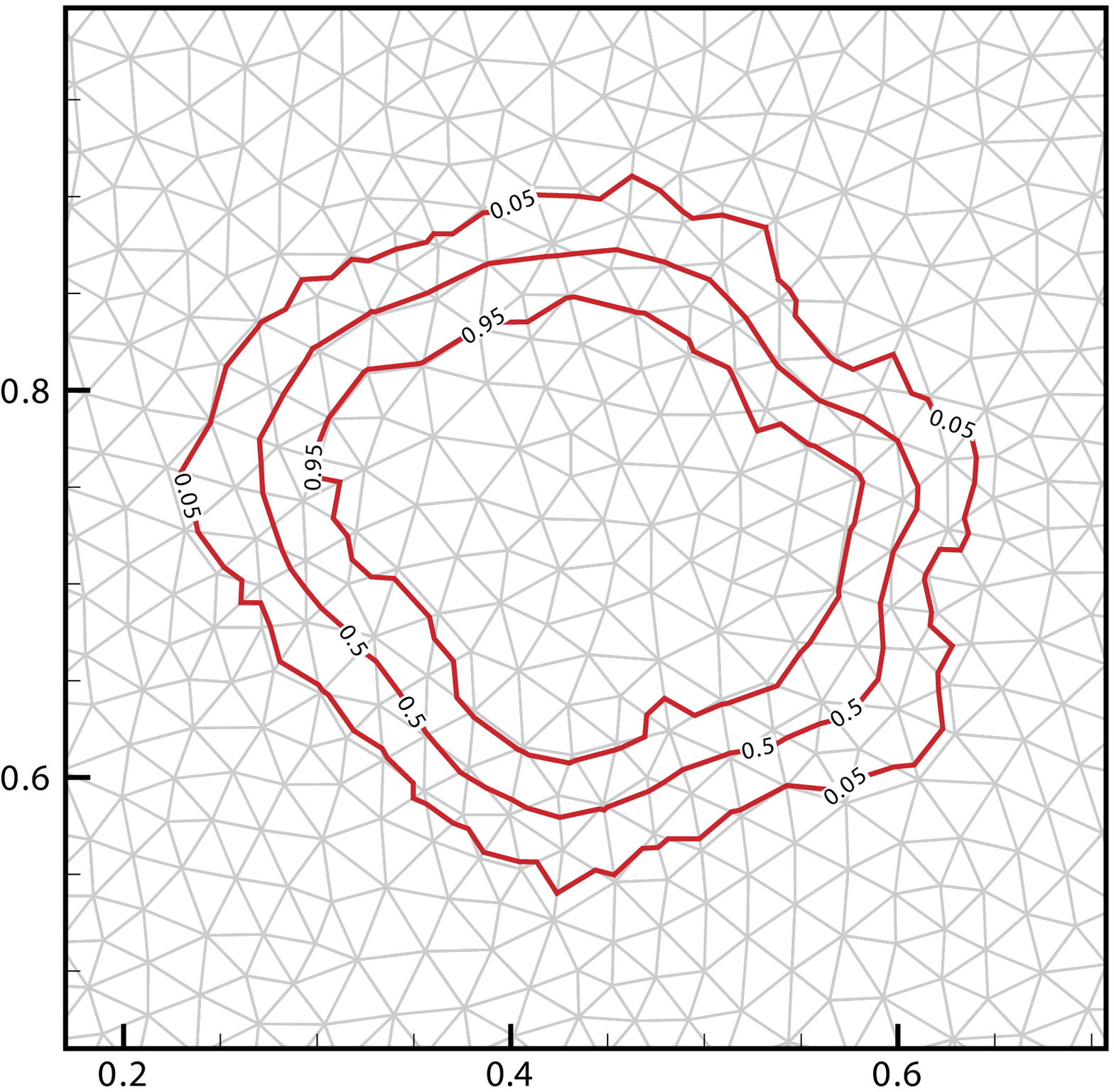}}
		\hspace{0.4cm}
	\centering
	\subfigure[$\text{N}=64$] {
		\centering
		\includegraphics[width=0.3\textwidth]{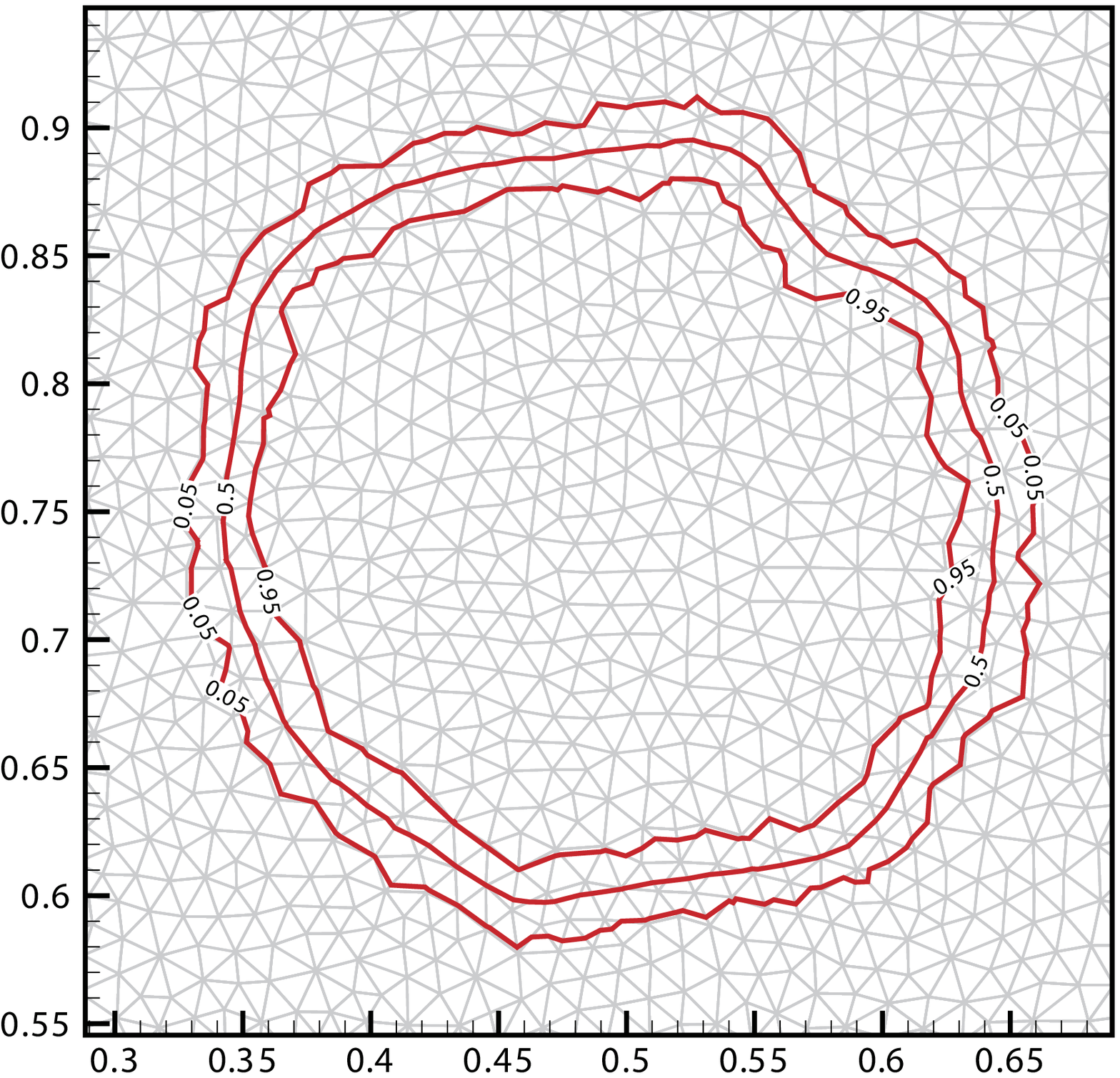}}
		\hspace{0.4cm}
	\centering
	\subfigure[$\text{N}=128$] {
		\centering
		\includegraphics[width=0.3\textwidth]{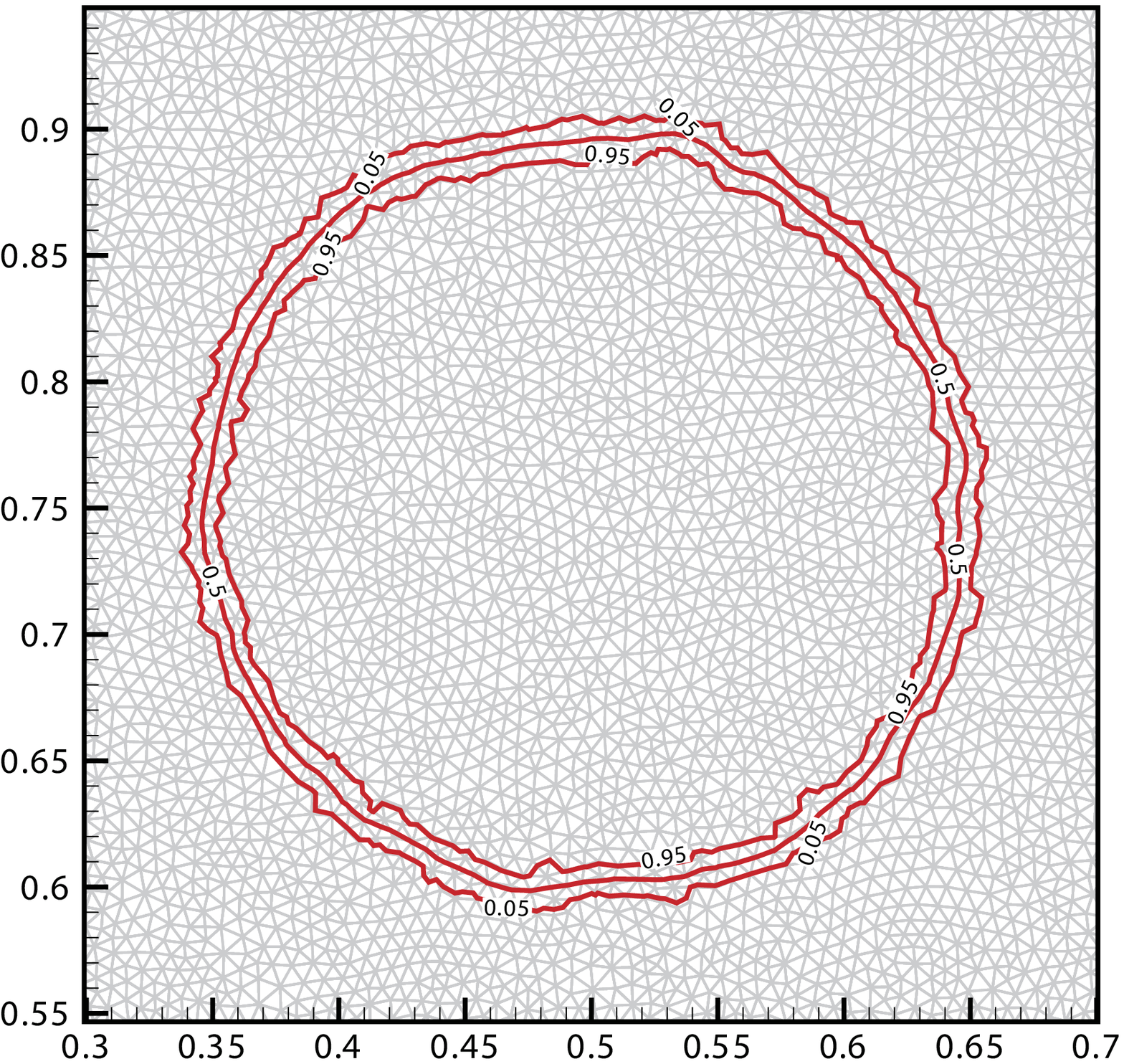}}
	\caption{\color{black}VOF $0.05$, $0.5$ and $0.95$ contour lines in Rider-Kothe single vortex test at $t=T$ using THINC-scaling scheme.}
\label{diff-uns}
\end{figure}

%~~~~~~~~~~~~~~~~~~~~~~~~~~~~~~~~ Numerical diffusion~~~~~~~~~~~~~~~~~~~~~~~~~~~~~~~~
We also show the capability of the proposed scheme in preserving the compactness of the interface thickness. From Fig.\ref{diff-uns}, showing the zoomed view of VOF $0.05$, $0.5$ and $0.95$ contour lines, we observe that the interface thickness remains intact in 3 cells after one rotation at $t=T$. It illustrates that the scheme is capable of preventing numerical diffusion and maintaining interface sharpness.

%~~~~~~~~~~~~~~~~~~~~~~~~~~~~~~~~~~~~~~~~~~~~~~~~~~~~~~~~~~~~~~~~~~~~~~~~~~~~~~~~~~~~~~~~~~~~~~~~~~~~~~~~~~~~~
\subsection{Three dimensional benchmark tests}

\subsubsection{Deformation flow}
We further test the proposed scheme in 3D using the vortex deformation flow test proposed in \cite{leveque1996high}. In this test, a sphere of radius $0.15$ initially centered at $\left(0.35,0.35,0.35\right)$ in an unit cube domain is advected by time dependent velocity field given by
\begin{equation}
\textbf{u}\left(\textbf{x},t\right)=\cos\left(\pi t/T\right)
\begin{pmatrix}
2\sin^2\left(\pi x\right)\sin\left(2\pi y\right)\sin\left(2\pi z\right)\\
-\sin\left(2\pi x\right)\sin^2\left(\pi y\right)\sin\left(2\pi z\right)\\
-\sin\left(2\pi x\right)\sin\left(2\pi y\right)\sin^2\left(\pi z\right)
\end{pmatrix}
\label{3D-vortex-vel}
\end{equation}
where time period $T=3$, and maximum CFL number is set to be $0.1$ and $0.5$ for Cartesian and unstructured tetrahedral grid respectively. Same as in 2D case, we use $\beta=6.0/\Delta$, where $\Delta=\min\left(\Delta x,\Delta y,\Delta z\right)$ in case of Cartesian grid. For tetrahedral grid, we set $\Delta=\min\left(l_i\right)$ with $i=1,\cdots,6$, where the $\min$ operator is over every length $l_i$ of the edges of the tetrahedral element.

\begin{figure}[htbp]
	\centering
	\subfigure[$\text{N}=32$ $(t=T/2)$] {
		\centering
		\includegraphics[width=0.4\textwidth]{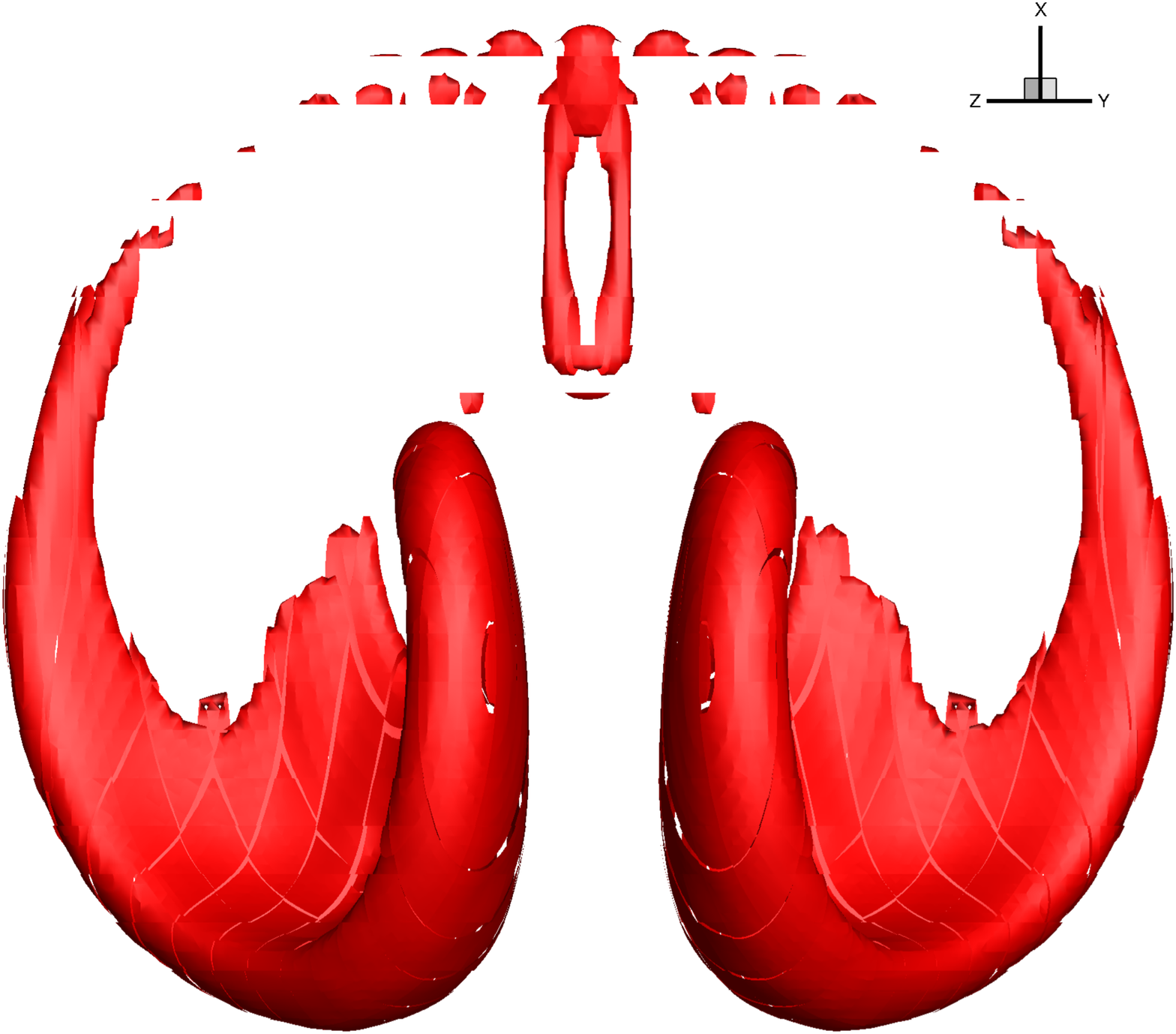}} \hspace{0.4cm}
	\centering
	\subfigure[$\text{N}=32$ $(t=T)$] {
		\centering
		\includegraphics[width=0.2\textwidth]{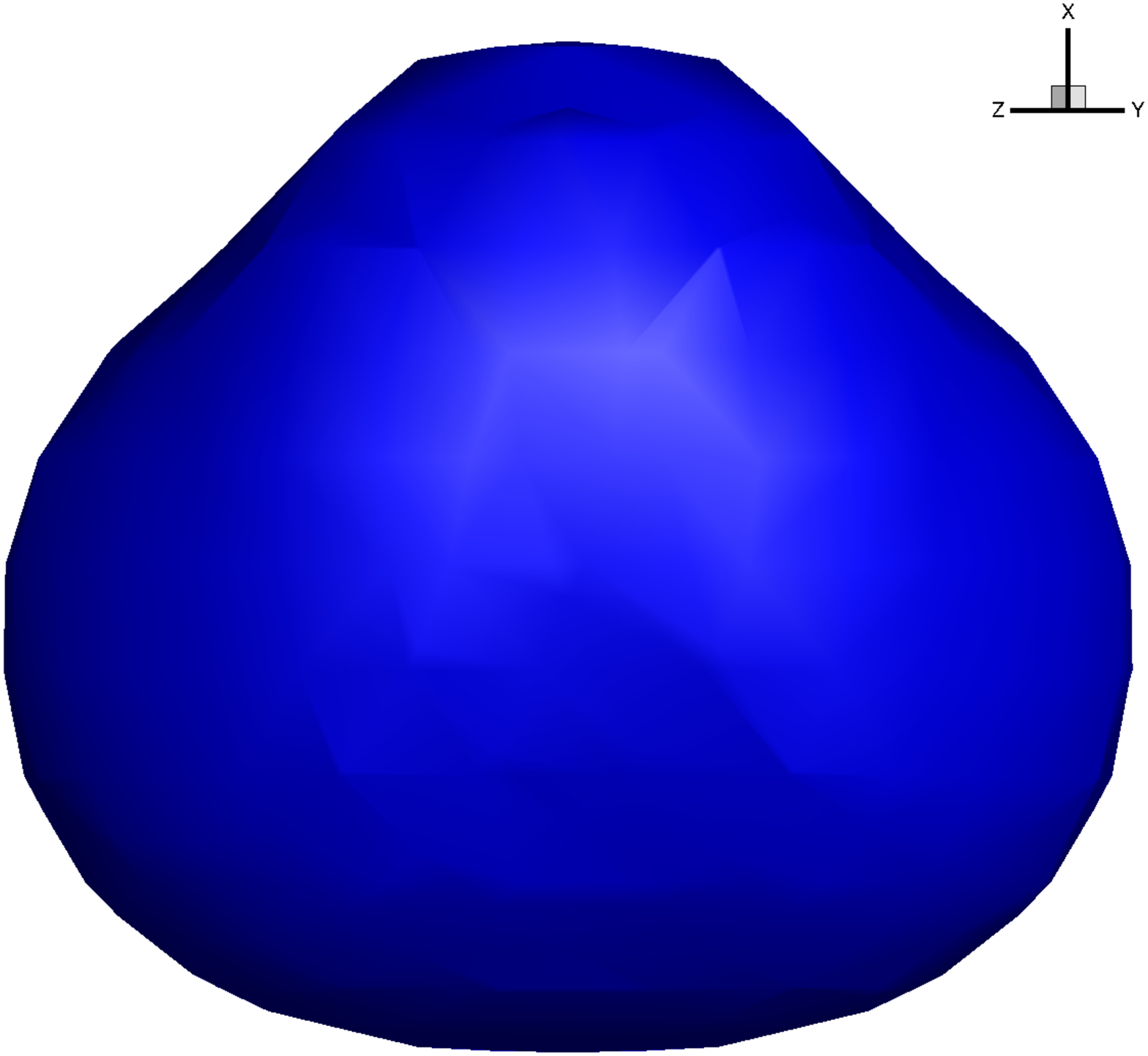}}
	\centering
	\subfigure[$\text{N}=64$ $(t=T/2)$] {
		\centering
		\includegraphics[width=0.4\textwidth]{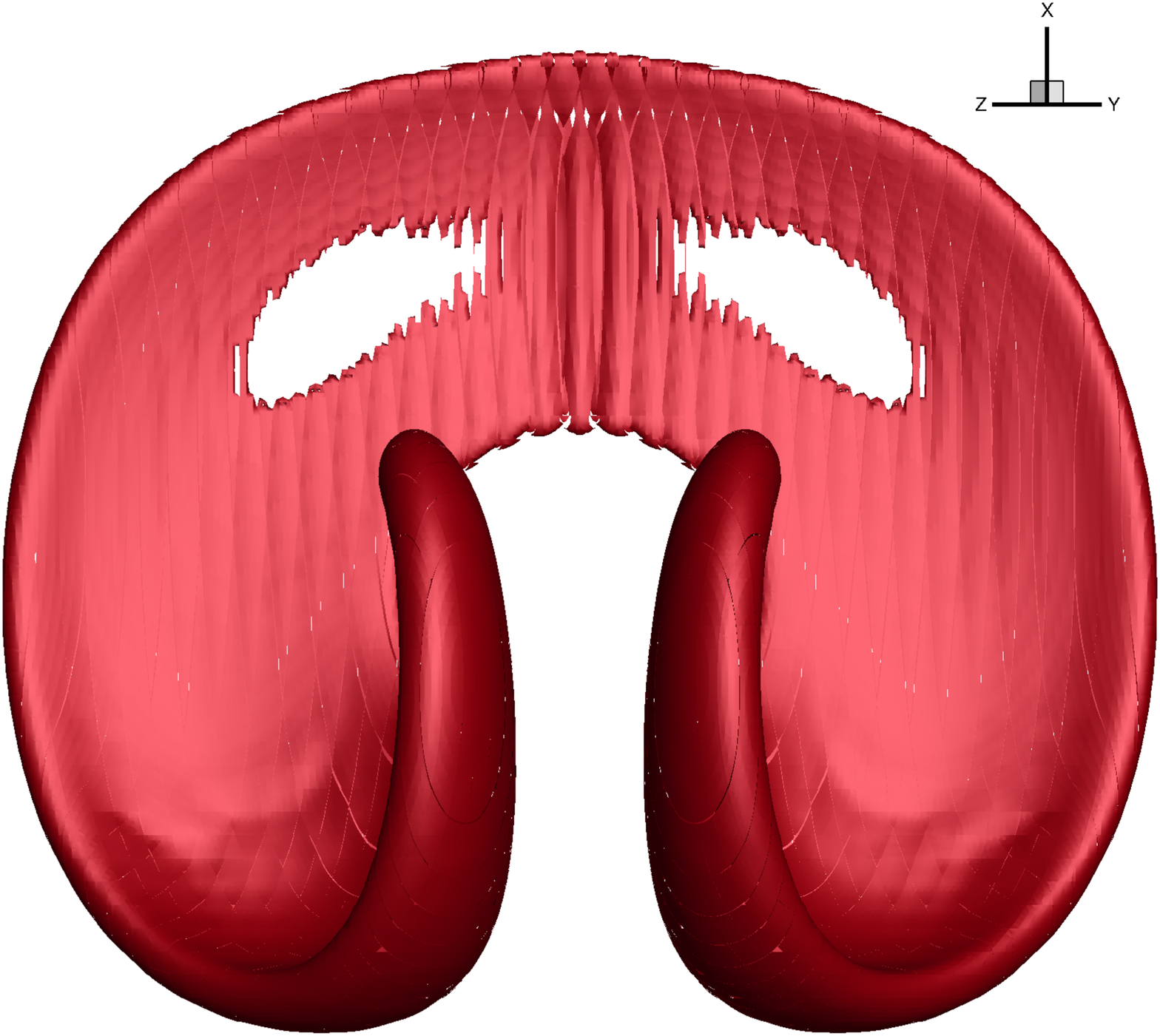}} \hspace{0.4cm}
	\centering
	\subfigure[$\text{N}=64$ $(t=T)$] {
		\centering
		\includegraphics[width=0.2\textwidth]{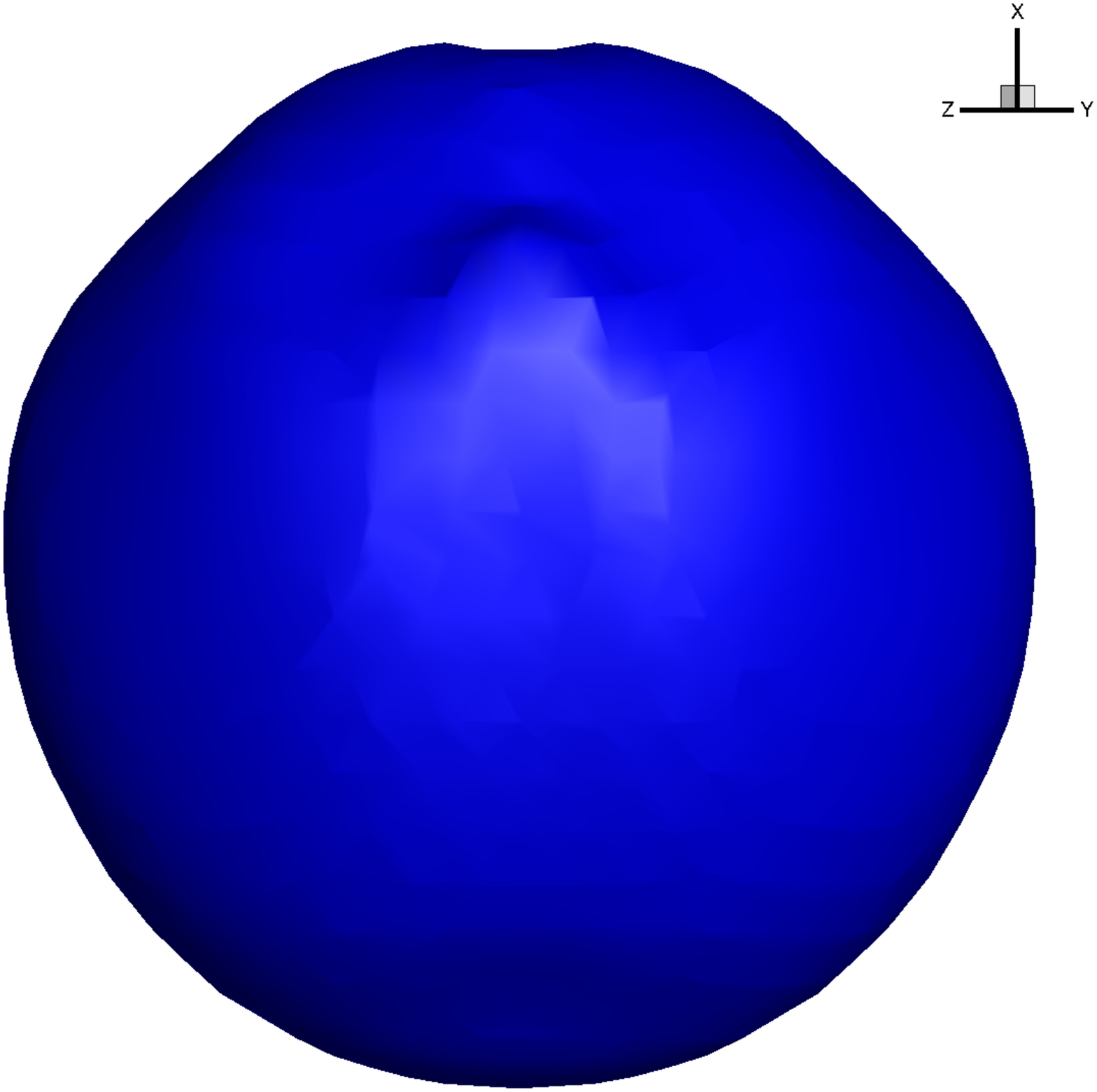}}
	\centering
	\subfigure[$\text{N}=128$ $(t=T/2)$] {
		\centering
		\includegraphics[width=0.4\textwidth]{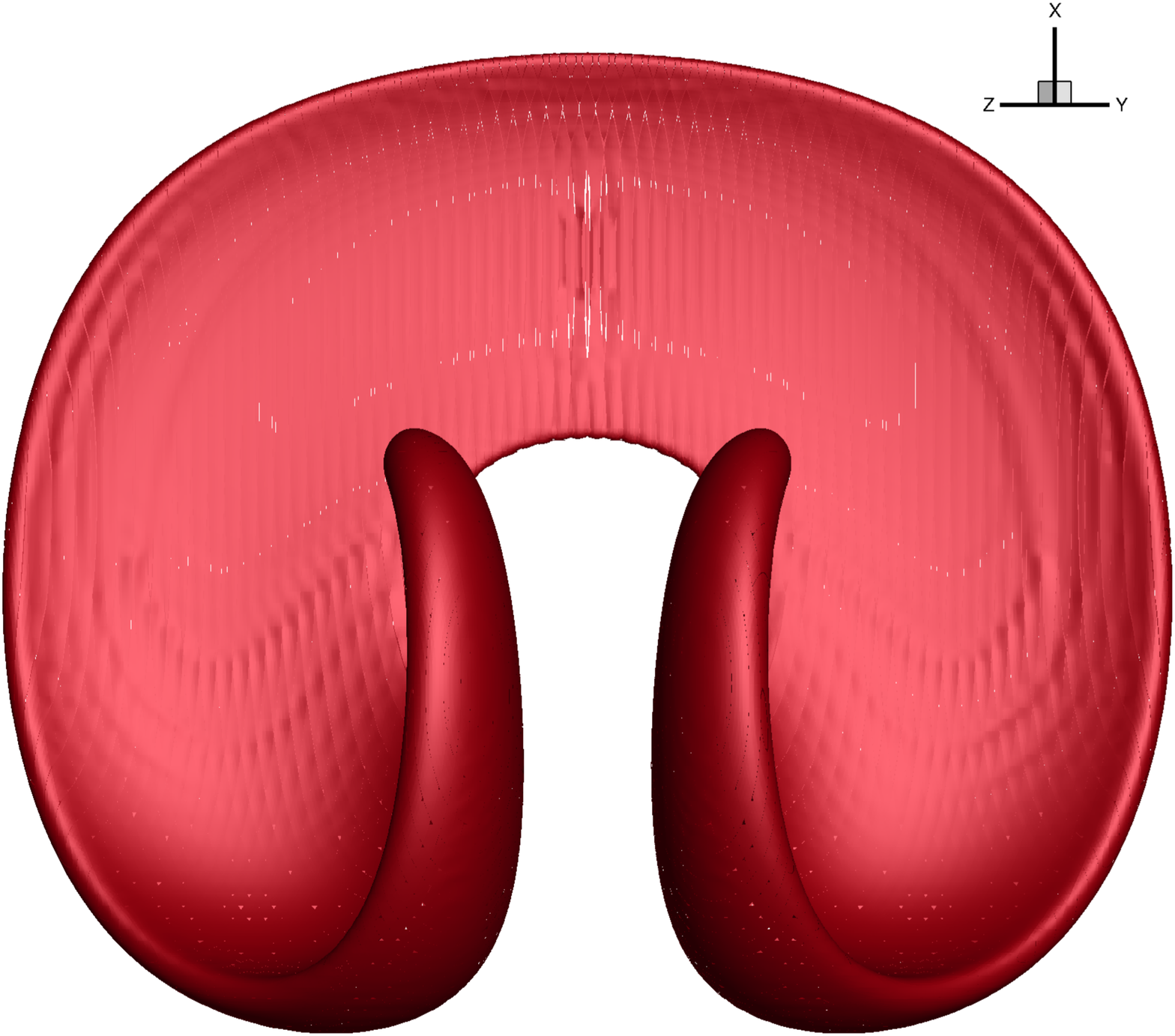}} \hspace{0.4cm}
	\centering
	\subfigure[$\text{N}=128$ $(t=T)$] {
		\centering
		\includegraphics[width=0.2\textwidth]{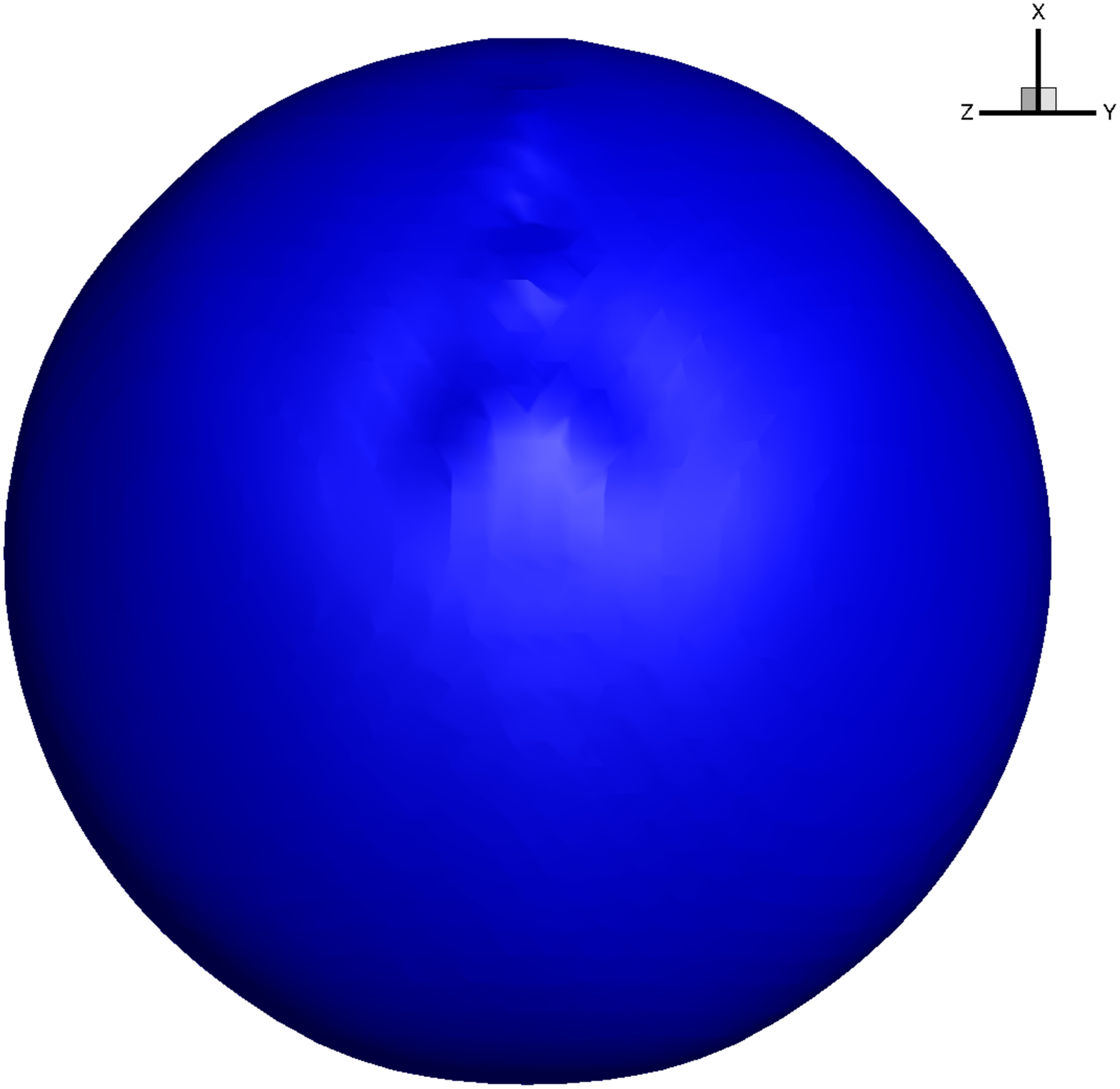}}
	\caption{Numerical results for 3D deformation flow on Cartesian mesh showing PSI at $t=T/2$ (left column) and VOF 0.5-isosurface at $t=T$ (right column).}
\label{3Ddeform-struct}
\end{figure}

Numerical results of THINC-scaling scheme on Cartesian grids with different mesh resolutions are shown in Fig.\ref{3Ddeform-struct}, where at $t=T/2$ the interface is represented by PSI of interface cells, reconstructed using cell-wise uniformly distributed 5 sample points in each direction, and at $t=T$ the restored sphere interface is represented by VOF 0.5-isosurface. As observed from the visual comparison, our results are more geometrically faithful in comparison with other latest variants of PLIC VOF methods reported in \cite{Maric2018VOF,Henning2019VOF}.

\begin{figure}[htbp]{
\begin{center}
\includegraphics[width=0.8\textwidth]{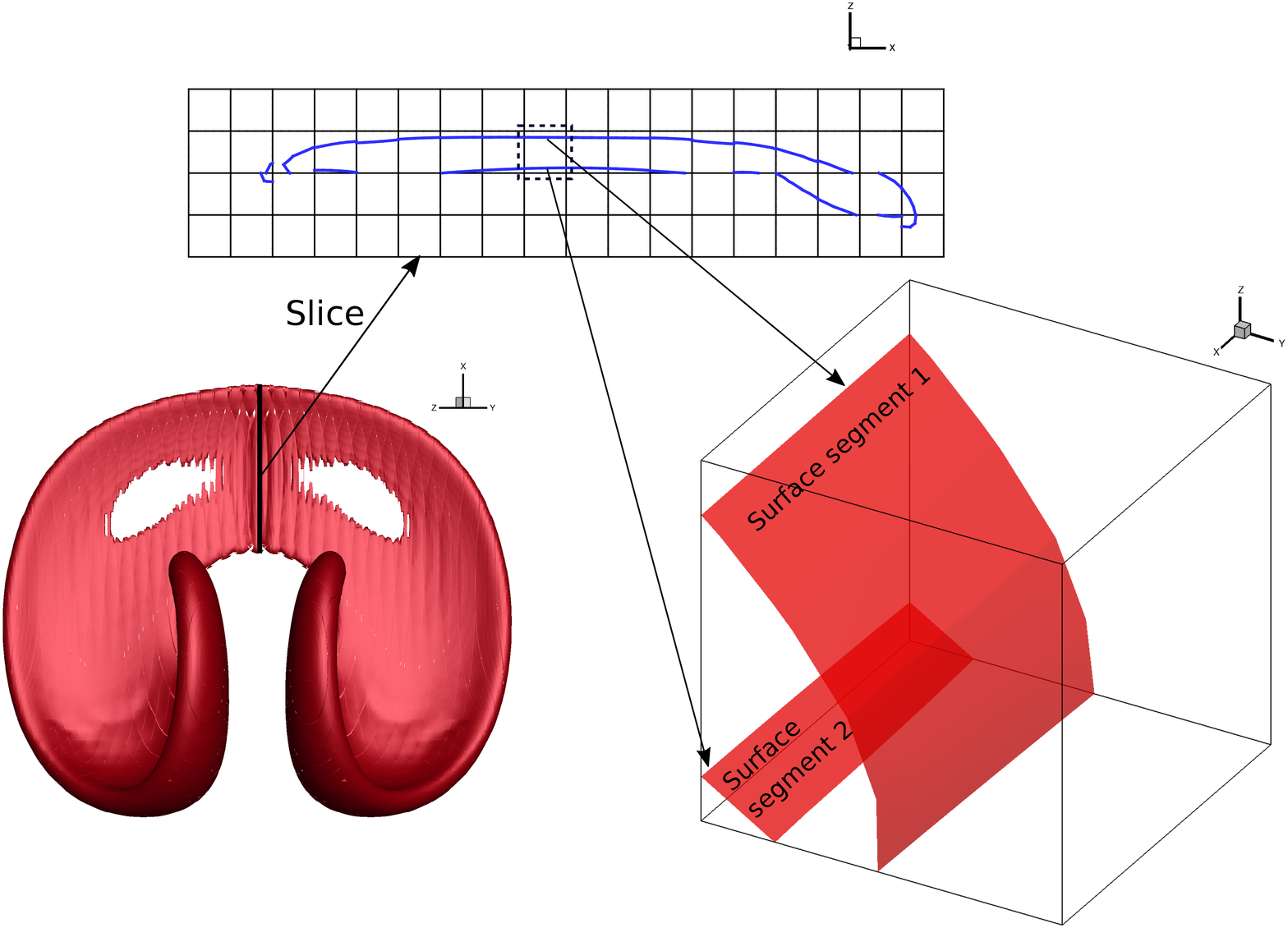}
\end{center}} 
\caption{2D slice of PSI obtained in 3D deformation flow and surface segments retrieved for thin film on Cartesian mesh with $\text{N}=64$.}
\label{3Dsurface_segs}
\end{figure}

Similar to the 2D case, we can also retrieve in 3D the two interface segments of a thin film structure under the grid resolution by use of the quadratic PSI. In Fig.\ref{3Dsurface_segs}, we show the 2D slice extracted from the thinned part of the deformed interface at $t=T/2$ for Cartesian mesh with $\text{N}=64$. It is observed that there are cells where sub-cell interface structures under grid resolution can be retrieved by two reconstructed surface segments as the PSI of quadratic polynomial. As these sub-grid-size structures can be resolved by using high-order surface polynomial, the THINC-scaling scheme is able to represent these thin structures more accurately as shown in the left column of Fig.\ref{3Ddeform-struct}
 which become otherwise the flotsams in the results of PLIC VOF methods. 

%~~~~~~~~~~~~~~~~~~~~~~~~~~~~~~~~Newton iteration convergence~~~~~~~~~~~~~~~~~~~~~~~~~~~~~~~~
We analyzed the Newton-Raphson convergence for solving \eqref{mass_constraint}, by calculating the average number of Newton iterations, $\bar{N}_{\rm iter}$, at each time step given by
\begin{equation}
\bar{N}_{\rm iter}=\frac{\sum_{i=1}^{N_d}N_{si}}{N_d},    
\end{equation}
where $N_d$ is the number of interface cells and $N_{si}$ denotes the number of Newton iterations in cell $\Omega_i$. The iteration process is terminated when tolerance of $10^{-11}$ is satisfied. {\color{black} Our numerical experiments on both Cartesian and unstructured grids, including all benchmark tests presented this paper,  show that the Newton-Raphson method takes few iterations with  $N_{\rm iter}$ ranging between 2 and 4 to converge throughout the computations due to its quadratic convergence property proved in the appendix A.}

%~~~~~~~~~~~~~~~~~~~~~~~~~~~~~Mass error~~~~~~~~~~~~~~~~~~~~~~~~~~~~~~
 {\color{black} In order to evaluate numerical conservativeness, we examined quantitatively the variation of  the total VOF summed up over the whole computational domain during numerical experiments.  
 It is observed that the total VOF values  for all numerical tests on different meshes remain unchanged up to the machine precision. Numerical conservativeness is rigorously guaranteed as the VOF field is computed by a finite volume formulation in flux form. }

%~~~~~~~~~~~~~~~~~~~~~~~~~~~~~Numerical errors and convergence rates~~~~~~~~~~~~~~~~~~~~~~~~~~~~~~
A quantitative comparison with other VOF methods on Cartesian grid and unstructured tetrahedral grid are outlined in Table \ref{3Ddef-struct-comparison} and Table \ref{3Ddeform-uns-comparison} respectively. The THINC-scaling scheme shows competitive results as compared with other VOF and hybrid methods.

\begin{table}[ht]
\centering
\caption{Numerical errors $E\left(L_1\right)$ and convergence rates for 3D deformation flow test on Cartesian grid.}
\begin{tabular}[t]{lSSSSS} \toprule
{$\textbf{Methods}$}               & {$32$}              & {Order} & {$64$}               & {Order} & {$128$} \\ \midrule
{THINC-scaling}                    & {$7.57\times10^{-3}$} & {1.53}  & {$2.62\times10^{-3}$}   & {1.96}  & {$6.72\times10^{-4}$}\\
{THINC/QQ} \cite{xie2017toward}    & {$7.96\times10^{-3}$} & {1.46}  & {$2.89\times10^{-3}$}   & {1.68}  & {$9.05\times10^{-4}$}\\
{UMTHINC} \cite{xie2017unstructured} & {$8.06\times10^{-3}$} & {1.41}  & {$3.04\times10^{-3}$}   & {1.69}  & {$9.40\times10^{-4}$}\\
{UFVFC-Swartz \cite{Maric2018VOF}} & {$5.86\times10^{-3}$} & {1.91}  & {$1.56\times10^{-3}$}   & {2.34}  & {$3.08\times10^{-4}$}\\  {Owkes and Desjardins \cite{Owkes2014VOF}}  & {$6.98\times10^{-3}$} & {1.73}  & {$2.10\times10^{-3}$}   & {1.89}  & {$5.62\times10^{-4}$}\\  
{isoAdvector-plicRDF \cite{Henning2019VOF}} & {$8.36\times10^{-3}$} & {1.36}  & {$3.25\times10^{-3}$}   & {2.31}  & {$6.57\times10^{-4}$}\\
{Youngs \cite{jofre2014}} & {$7.47\times10^{-3}$} & {1.43}  & {$2.77\times10^{-3}$}   & {1.77}  & {$8.14\times10^{-4}$}\\
{LVIRA \cite{jofre2014}} & {$6.92\times10^{-3}$} & {1.51}  & {$2.43\times10^{-3}$}   & {1.93}  & {$6.37\times10^{-4}$}\\
{RK-3D \cite{hernandez2008new}} & {$7.85\times10^{-3}$} & {1.51}  & {$2.75\times10^{-3}$}   & {1.89}  & {$7.41\times10^{-4}$}\\
{FMFPA-3D \cite{hernandez2008new}} & {$7.44\times10^{-3}$} & {1.42}  & {$2.79\times10^{-3}$}   & {1.97}  & {$7.14\times10^{-4}$}\\
{Improved ELVIRA-3D \cite{lopez2008new}} & {$7.35\times10^{-3}$} & {1.45}  & {$2.69\times10^{-3}$}   & {2.05}  & {$6.51\times10^{-4}$}\\
{DS-MOF \cite{jemison2013MOF}} & {$5.72\times10^{-3}$} & {1.50}  & {$2.02\times10^{-3}$}   & {-}  & {-}\\
{DS-CLSVOF \cite{jemison2013MOF}} & {$6.92\times10^{-3}$} & {1.70}  & {$2.13\times10^{-3}$}   & {-}  & {-}\\
{DS-CLSMOF \cite{jemison2013MOF}} & {$4.81\times10^{-3}$} & {1.27}  & {$1.99\times10^{-3}$}   & {-}  & {-}\\ \bottomrule
\label{3Ddef-struct-comparison}
\end{tabular}
\end{table}

\begin{table}[ht]
\centering
\caption{Numerical errors $E\left(L_1\right)$ and convergence rates for 3D deformation flow test on unstructured tetrahedral grid.} 
\begin{tabular}[t]{lSSS} \toprule
{$\textbf{Methods}$}          & {$32$}                & {Order} & {$64$} \\ \midrule
{THINC-scaling}               & {\color{black}$6.64\times10^{-3}$} & {1.69}  & {$2.06\times10^{-3}$} \\
%{THINC/LS}                    & {$6.62\times10^{-3}$} & {1.69}  & {$2.05\times10^{-3}$} \\
{THINC/QQ}                    & {$8.72\times10^{-3}$} & {1.71}  & {$2.67\times10^{-3}$} \\
{isoAdvector-plicRDF \cite{Henning2019VOF}} & {$1.31\times10^{-2}$} & {1.06}  & {$6.34\times10^{-3}$} \\
{Youngs \cite{jofre2014}} & {$1.02\times10^{-2}$} & {1.20}  & {$4.45\times10^{-3}$} \\
{LVIRA \cite{jofre2014}} & {$1.02\times10^{-2}$} & {1.53}  & {$3.54\times10^{-3}$} \\\bottomrule
\label{3Ddeform-uns-comparison}
\end{tabular}
\end{table}

%~~~~~~~~~~~~~~~~~~~~~~~~~~~~~~~~~~~~~~~~~~~~~~~~~~~~~~~~~~~~~~~~~~~~~~~~~~~~~~~~~~~~~~~~
\subsubsection{Shear flow}
We carried out the 3D shear flow benchmark test introduced in \cite{Liovic2006PLIC}, and later used to evaluate other geometric VOF methods \cite{Maric2018VOF,Henning2019VOF,jofre2014}. In this test, a cuboid computational domain of  $\left[1,1,2\right]$ was partitioned with a mesh generated from uniformly distributed nodes with $\text{N}$ in $x$ and $y$ directions and $2\text{N}$ in $z$ direction. A sphere initially centered at $\left(0.5,0.75,0.5\right)$ is transported in $+z$ direction, and then returned to its initial position by time dependent velocity field given by
\begin{equation}
\textbf{u}\left(\textbf{x},t\right)=\cos\left(\pi t/T\right)
\begin{pmatrix}
\sin^2\left(\pi x\right)\sin\left(2\pi y\right)\\
-\sin\left(2\pi x\right)\sin^2\left(\pi y\right)\\
\left(1-2r\right)^2
\end{pmatrix}
\label{3D-shear-vel}
\end{equation}
where $r=\sqrt{\left(x-0.5\right)^2+\left(y-0.5\right)^2}$, time period $T=3$, and maximum CFL number is set to be $0.1$ for Cartesian grid and $0.5$ for tetrahedral grid in our simulations.

\begin{figure}[htbp]
	\centering
	\subfigure[$\text{N}=32$ $(t=T/2)$] {
		\centering
		\includegraphics[width=0.4\textwidth]{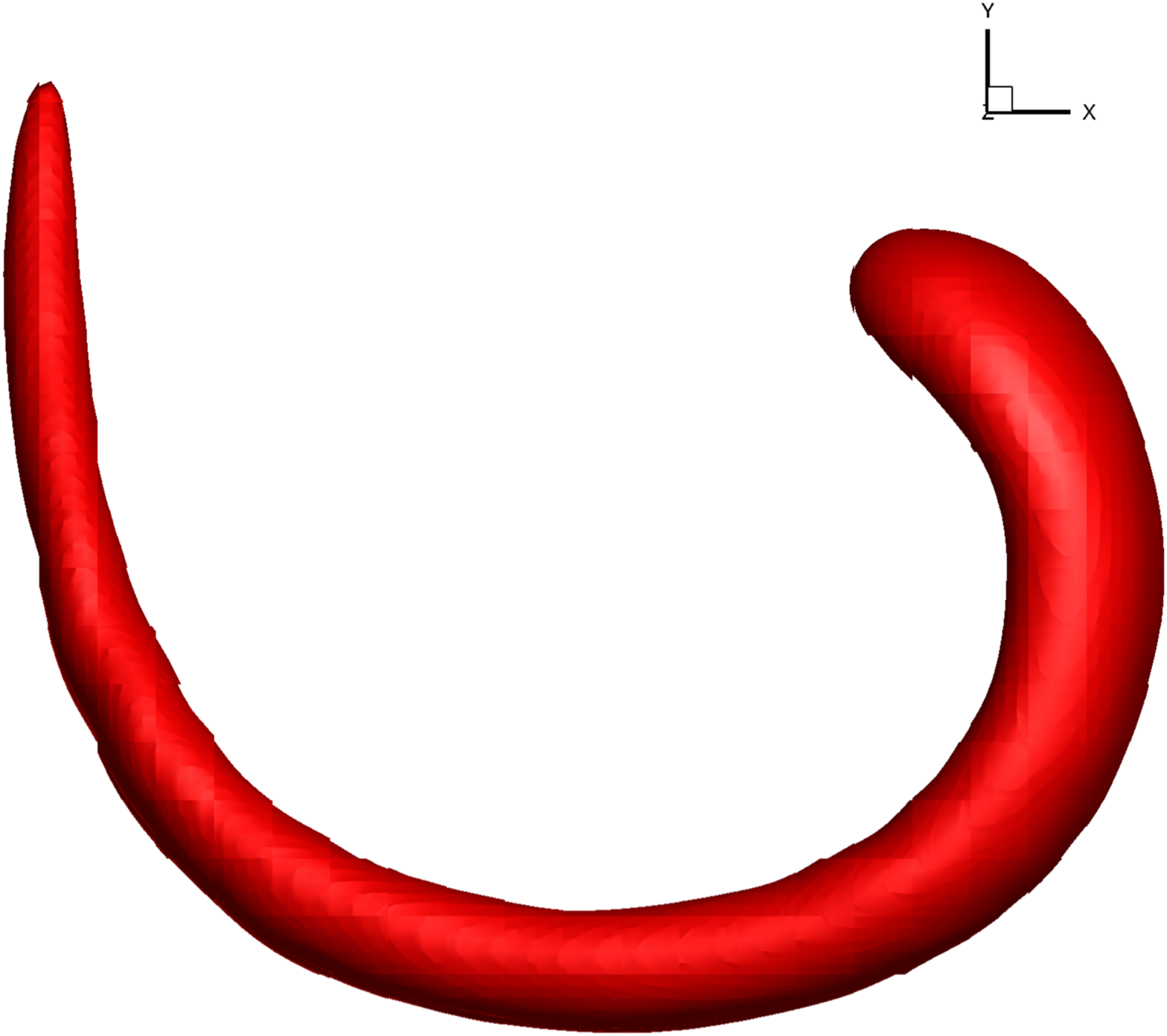}} \hspace{0.3cm}
	\centering
	\subfigure[$\text{N}=32$ $(t=T)$] {
		\centering
		\includegraphics[width=0.2\textwidth]{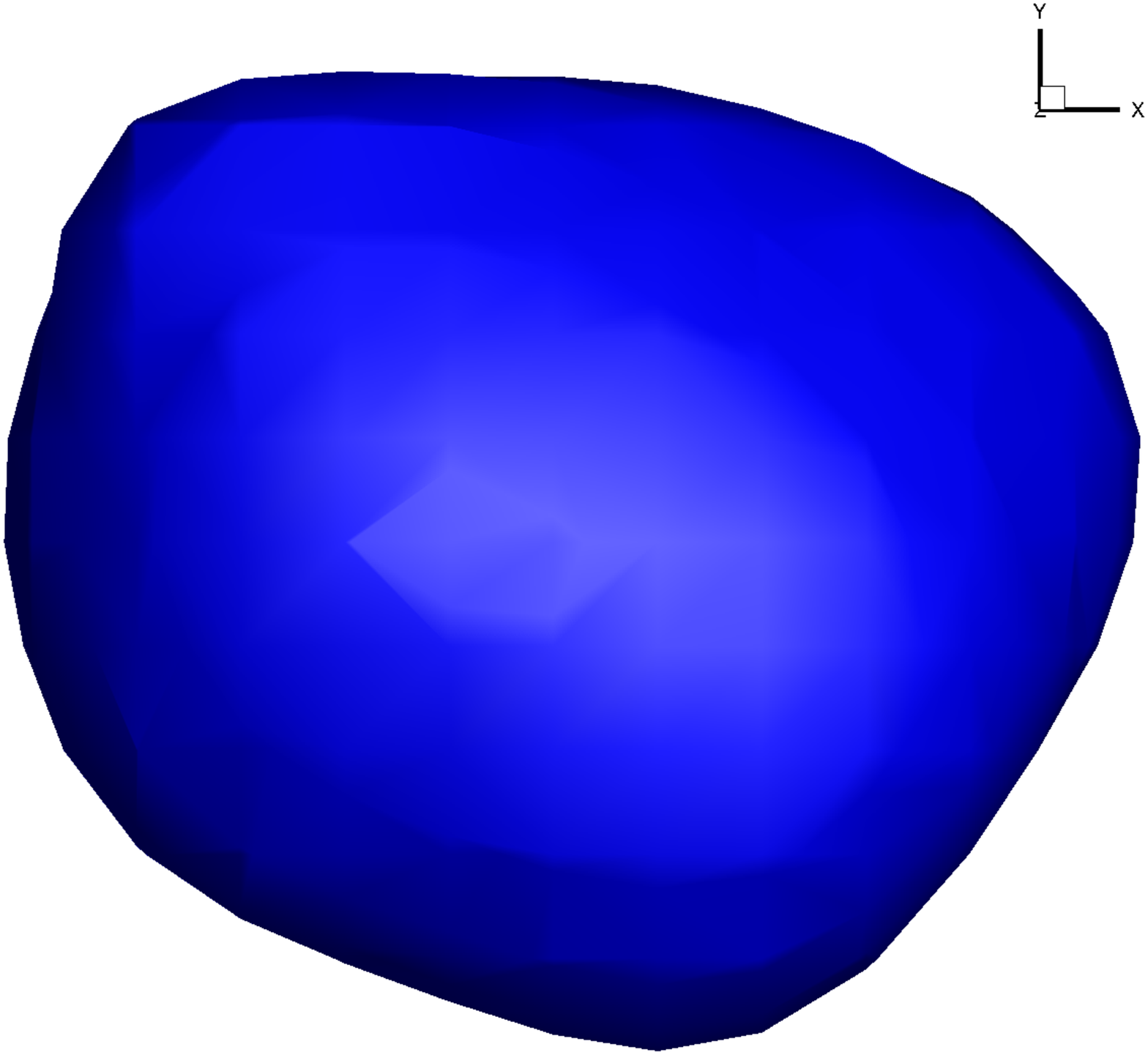}}
	\centering
	\subfigure[$\text{N}=64$ $(t=T/2)$] {
		\centering
		\includegraphics[width=0.4\textwidth]{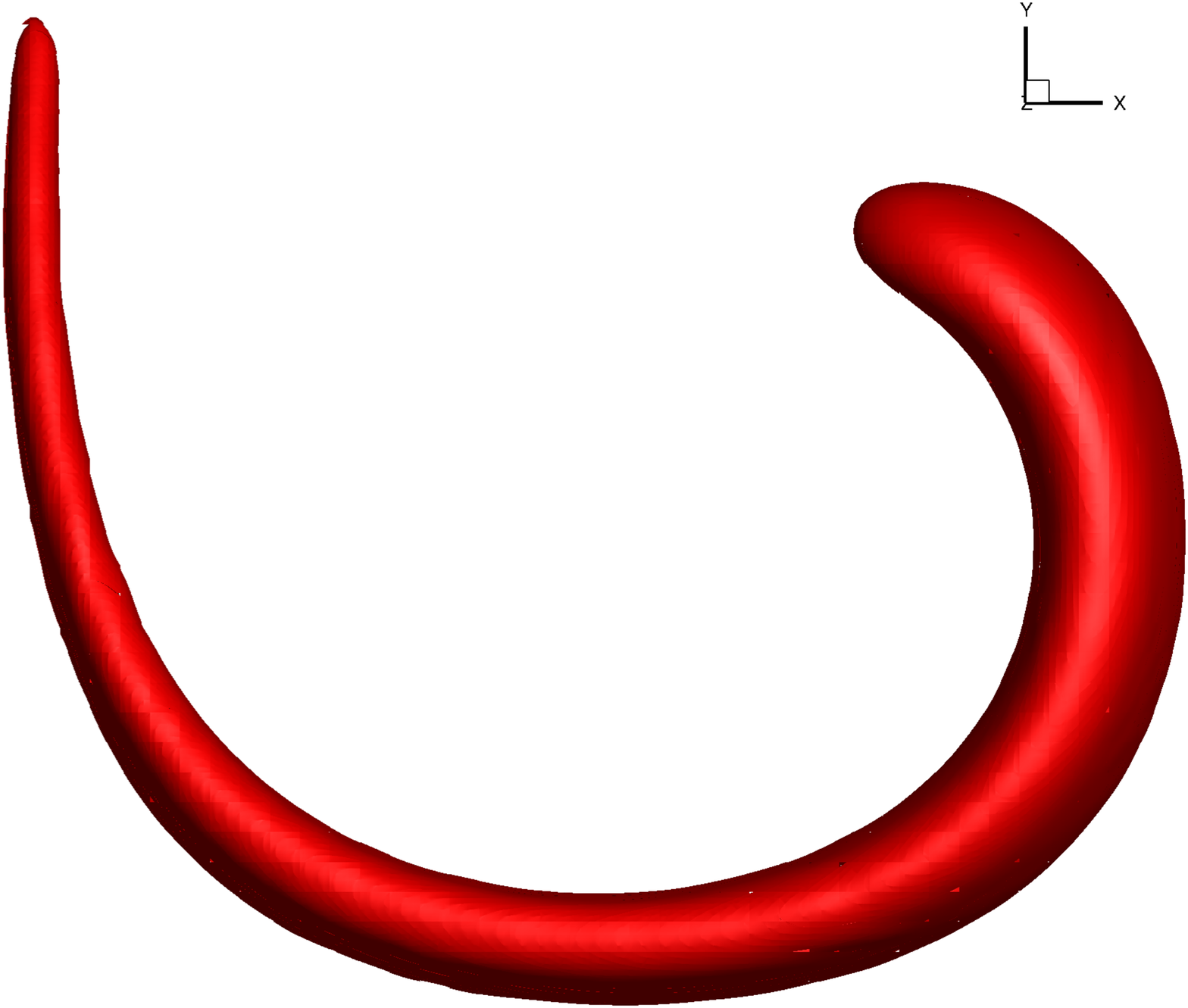}} \hspace{0.3cm}
	\centering
	\subfigure[$\text{N}=64$ $(t=T)$] {
		\centering
		\includegraphics[width=0.2\textwidth]{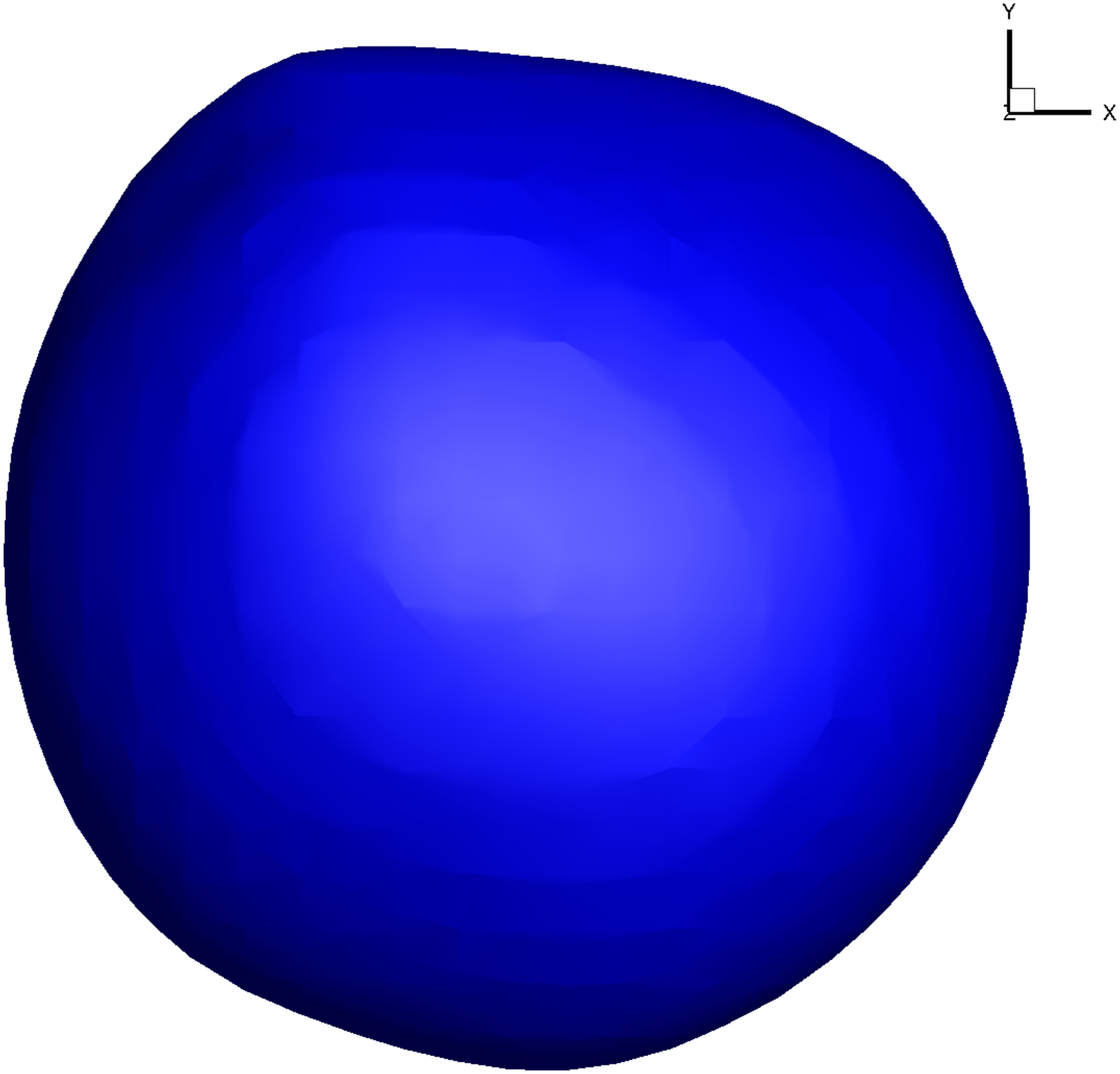}}
	\centering
	\subfigure[$\text{N}=128$ $(t=T/2)$] {
		\centering
		\includegraphics[width=0.4\textwidth]{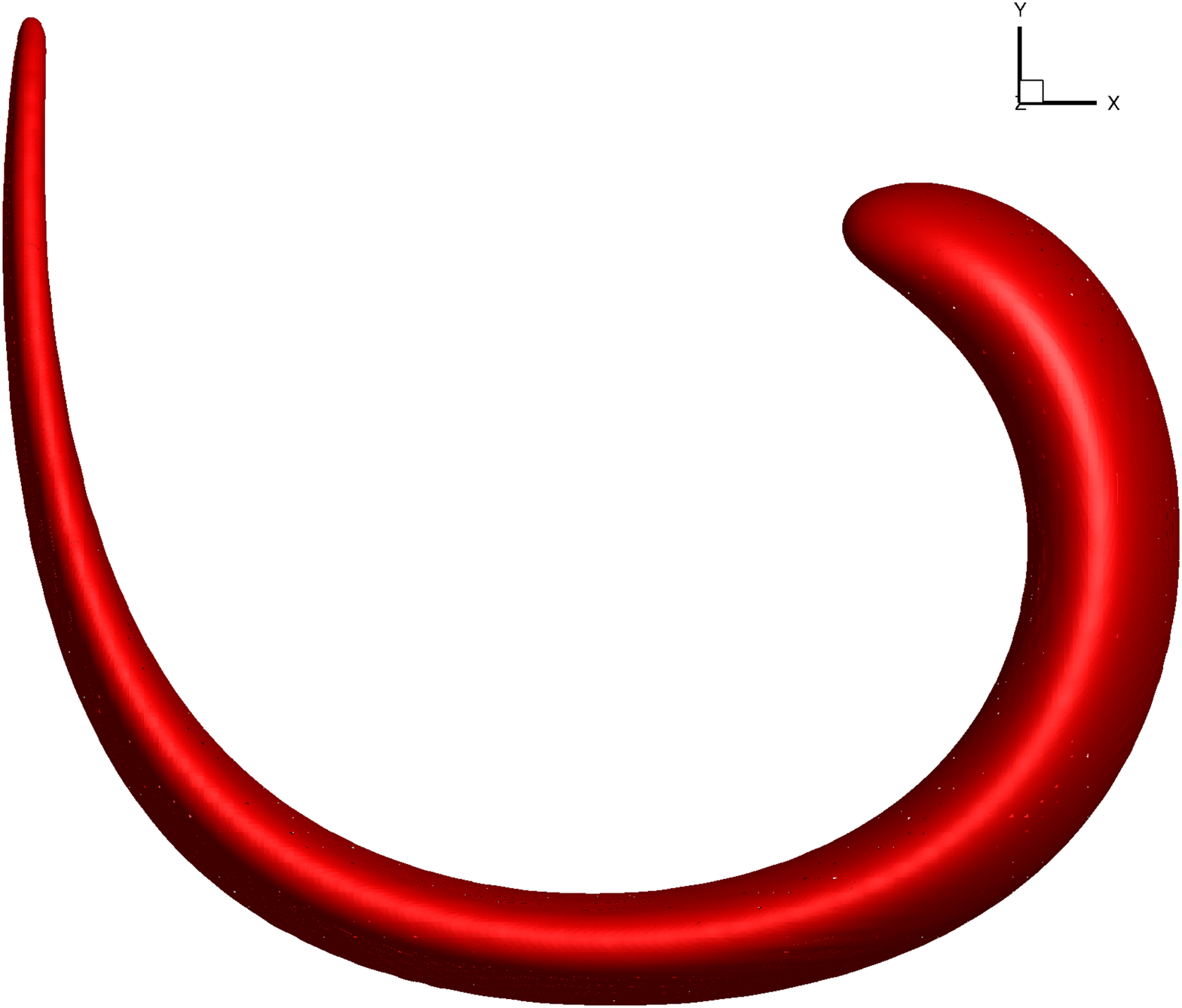}} \hspace{0.3cm}
	\centering
	\subfigure[$\text{N}=128$ $(t=T)$] {
		\centering
		\includegraphics[width=0.2\textwidth]{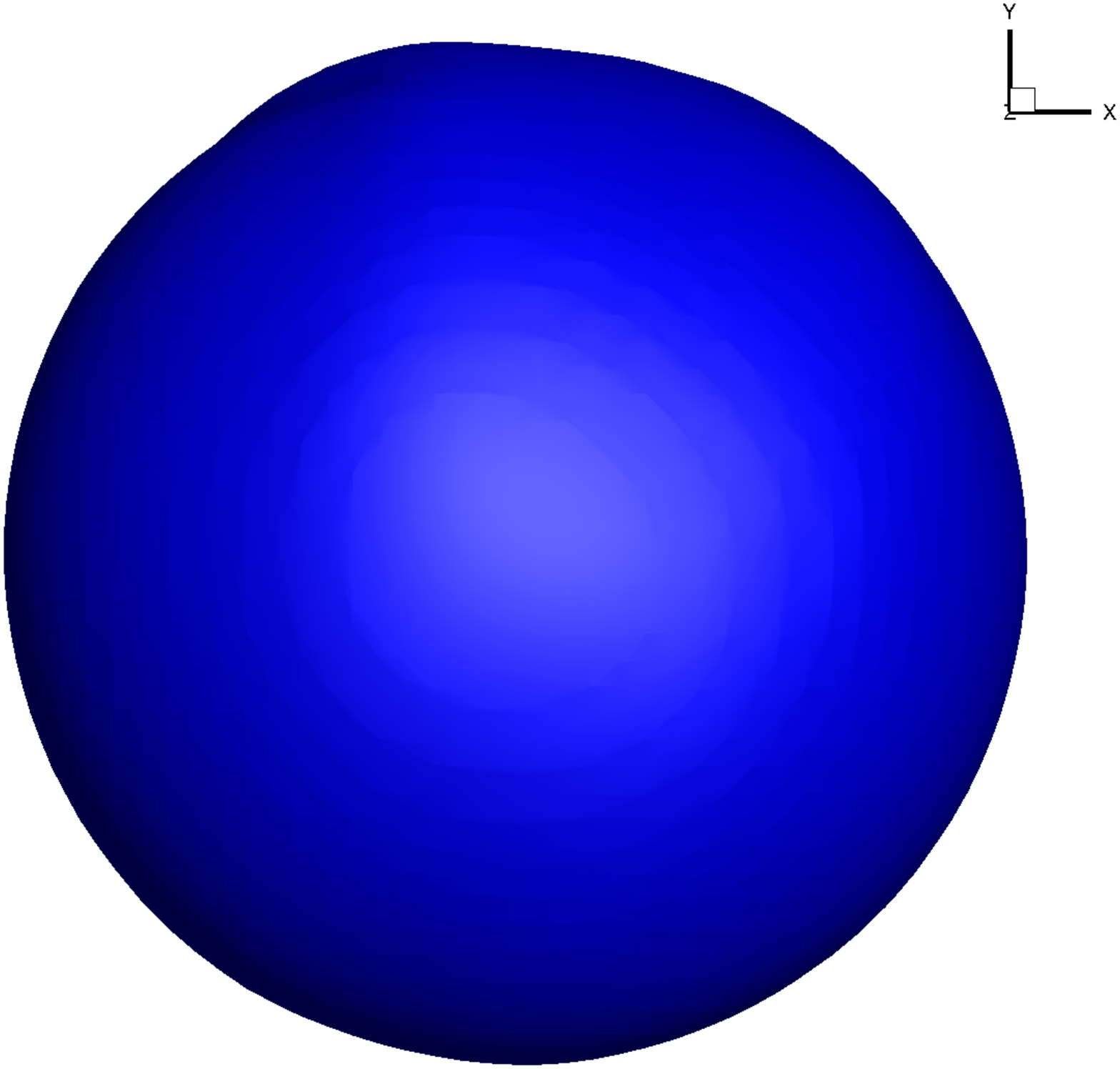}}
	\caption{Numerical results for 3D shear flow on Cartesian mesh. Shown are PSI at $t=T/2$ (left column), and  VOF 0.5-isosurfaces at $t=T$ (right column).}
\label{3Dshear-struct}
\end{figure}

Numerical results of THINC-scaling scheme on Cartesian grids with different resolutions are shown in Fig.\ref{3Dshear-struct}. As can be seen from PSI at $t=T/2$ (left column in Fig.\ref{3Dshear-struct}), a good resolution of the elongated thin tail part is obtained. Also, at $t=T$ the restored sphere achieves geometrically faithful results with more refined mesh. On visual comparison with the PLIC VOF methods results shown in \cite{Maric2018VOF,Henning2019VOF}, we notice more accurate interface representation in our results, especially as observed from the elongated tail part result at $t=T/2$.   

The quantitative comparison of THINC-scaling scheme with other VOF methods, including some of the latest variants of PLIC VOF methods on Cartesian grids are shown in Table \ref{3Dshear-struct-comparison}.

\begin{table}[ht]
\centering
\caption{Numerical errors $E\left(L_1\right)$ and convergence rates for 3D shear flow test on Cartesian grid.} 
\begin{tabular}[t]{lSSSSS} \toprule
{$\textbf{Methods}$}               & {$32$}              & {Order} & {$64$}               & {Order} & {$128$} \\ \midrule
{THINC-scaling}                    & {$3.59\times10^{-3}$} & {1.76}  & {$1.06\times10^{-3}$}   & {1.75}  & {$3.16\times10^{-4}$}\\
{UFVFC-Swartz \cite{Maric2018VOF}} & {$1.97\times10^{-3}$} & {2.21}  & {$4.25\times10^{-4}$}   & {1.77}  & {$1.24\times10^{-4}$}\\  
{isoAdvector-plicRDF \cite{Henning2019VOF}} & {$4.06\times10^{-3}$} & {1.91}  & {$1.08\times10^{-3}$}   & {2.14}  & {$2.44\times10^{-4}$}\\
{Youngs \cite{jofre2014}} & {$4.06\times10^{-3}$} & {1.66}  & {$1.29\times10^{-3}$}   & {1.24}  & {$5.45\times10^{-4}$}\\
{CVTNA-PCFSC \cite{Liovic2006PLIC}} & {$2.86\times10^{-3}$} & {2.00}  & {$7.14\times10^{-4}$}   & {2.19}  & {$1.56\times10^{-4}$}\\\bottomrule
\label{3Dshear-struct-comparison}
\end{tabular}
\end{table}

\begin{figure}[htbp]
    \centering
	\subfigure[$\text{N}=32$ $(t=T/2)$] {
		\centering
		\includegraphics[width=0.4\textwidth]{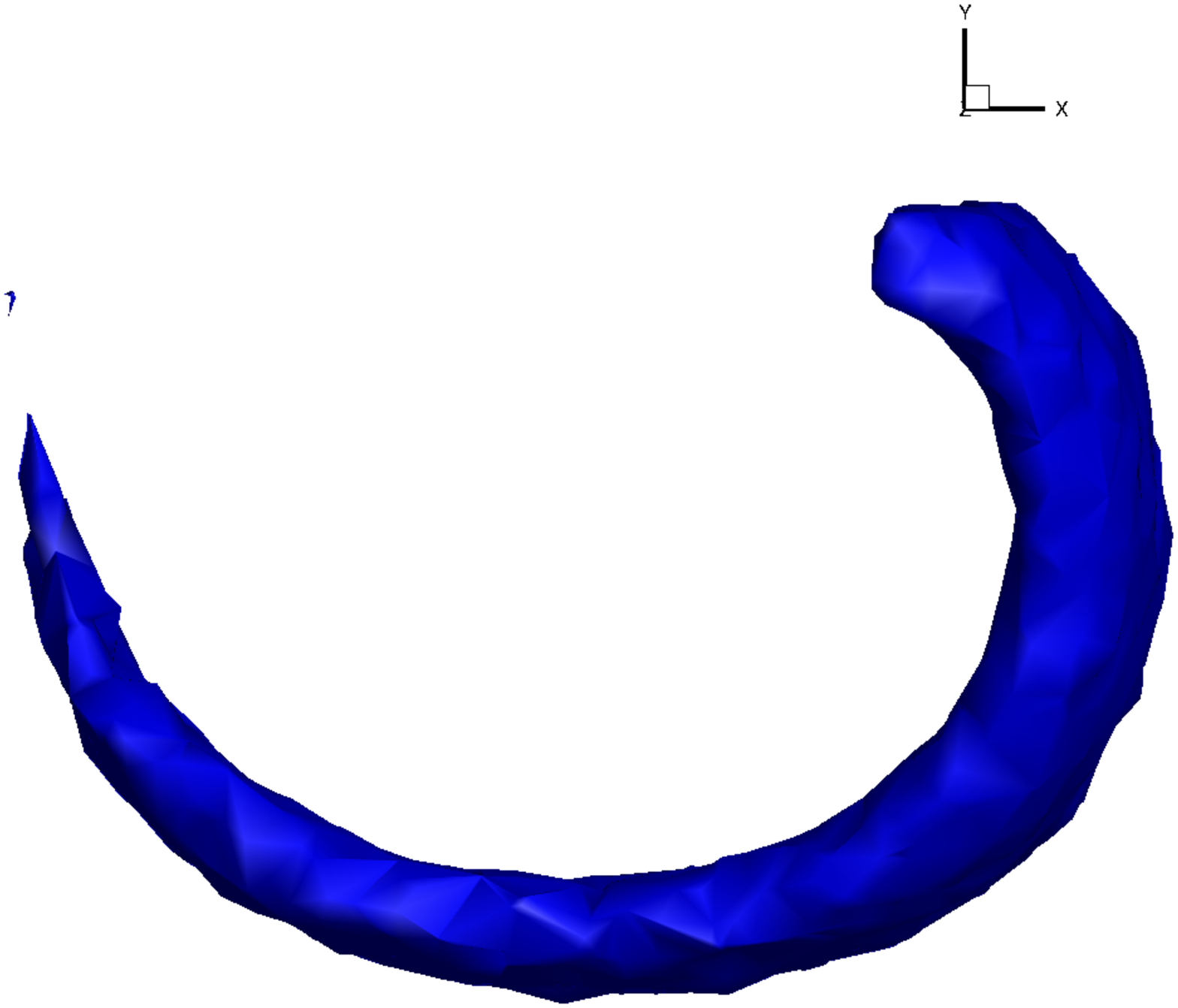}} \hspace{0.4cm}
	\centering
	\subfigure[$\text{N}=32$ $(t=T)$] {
		\centering
		\includegraphics[width=0.2\textwidth]{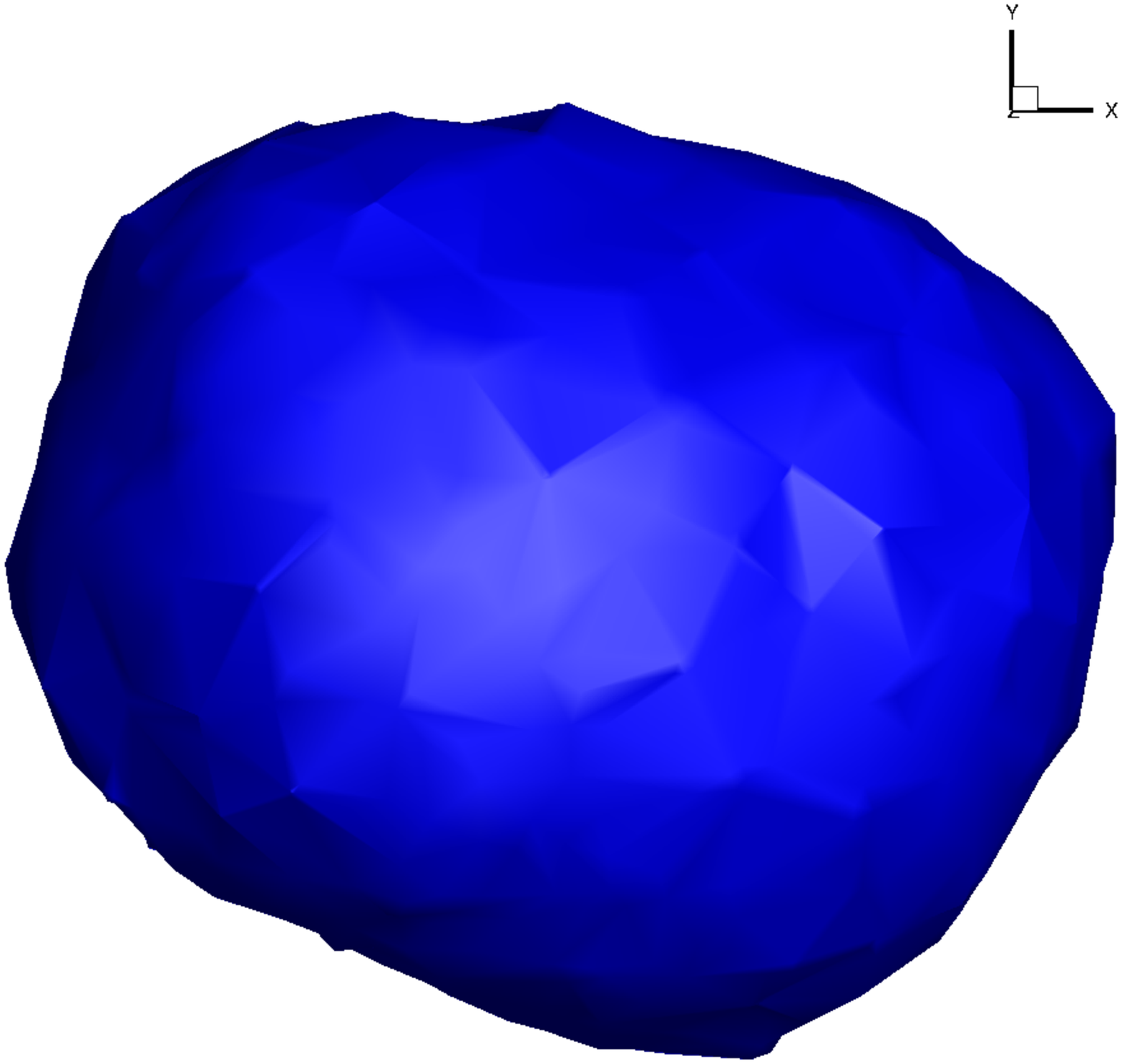}}
	\centering
	\subfigure[$\text{N}=64$ $(t=T/2)$] {
		\centering
		\includegraphics[width=0.4\textwidth]{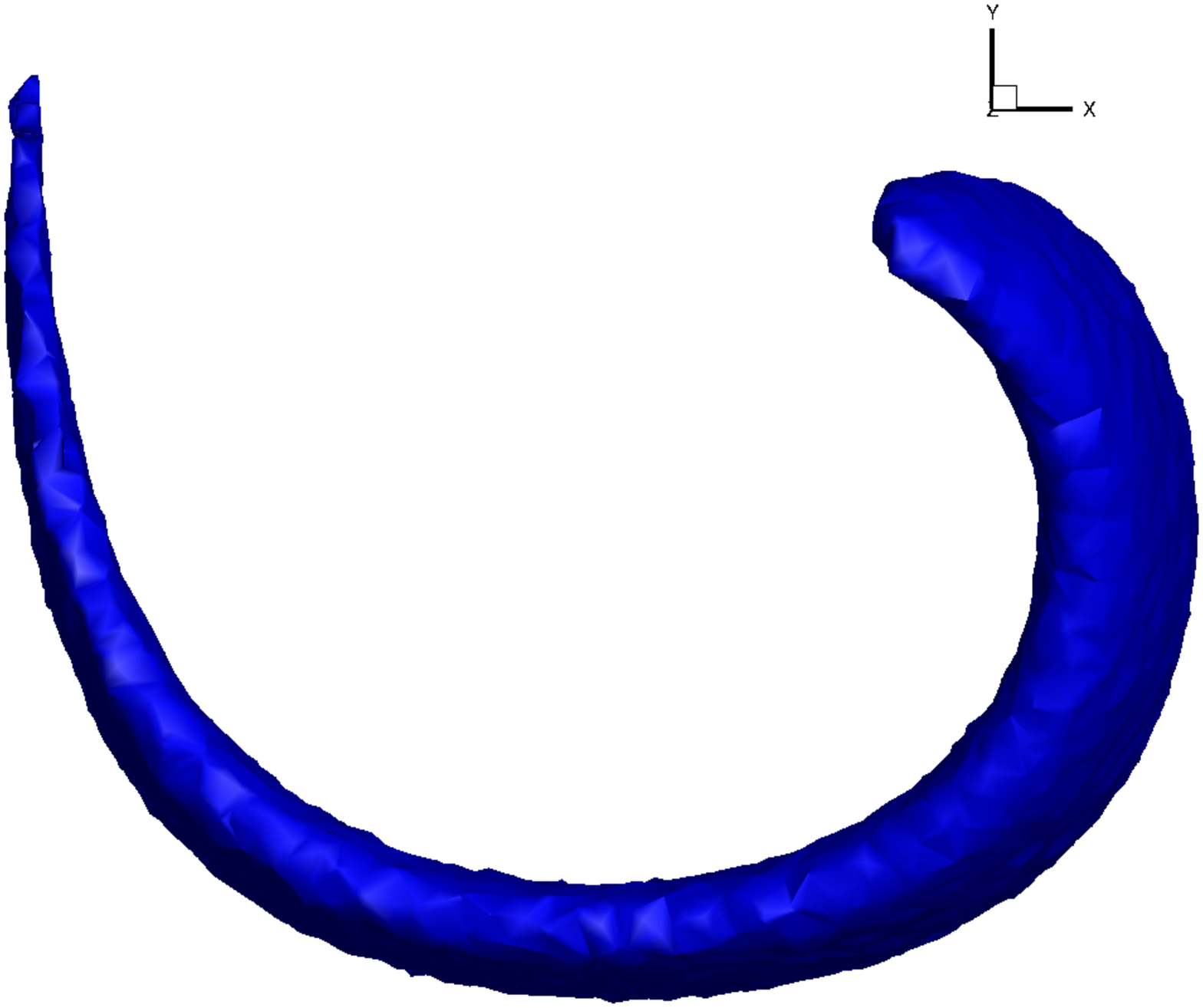}} \hspace{0.4cm}
	\centering
	\subfigure[$\text{N}=64$ $(t=T)$] {
		\centering
		\includegraphics[width=0.2\textwidth]{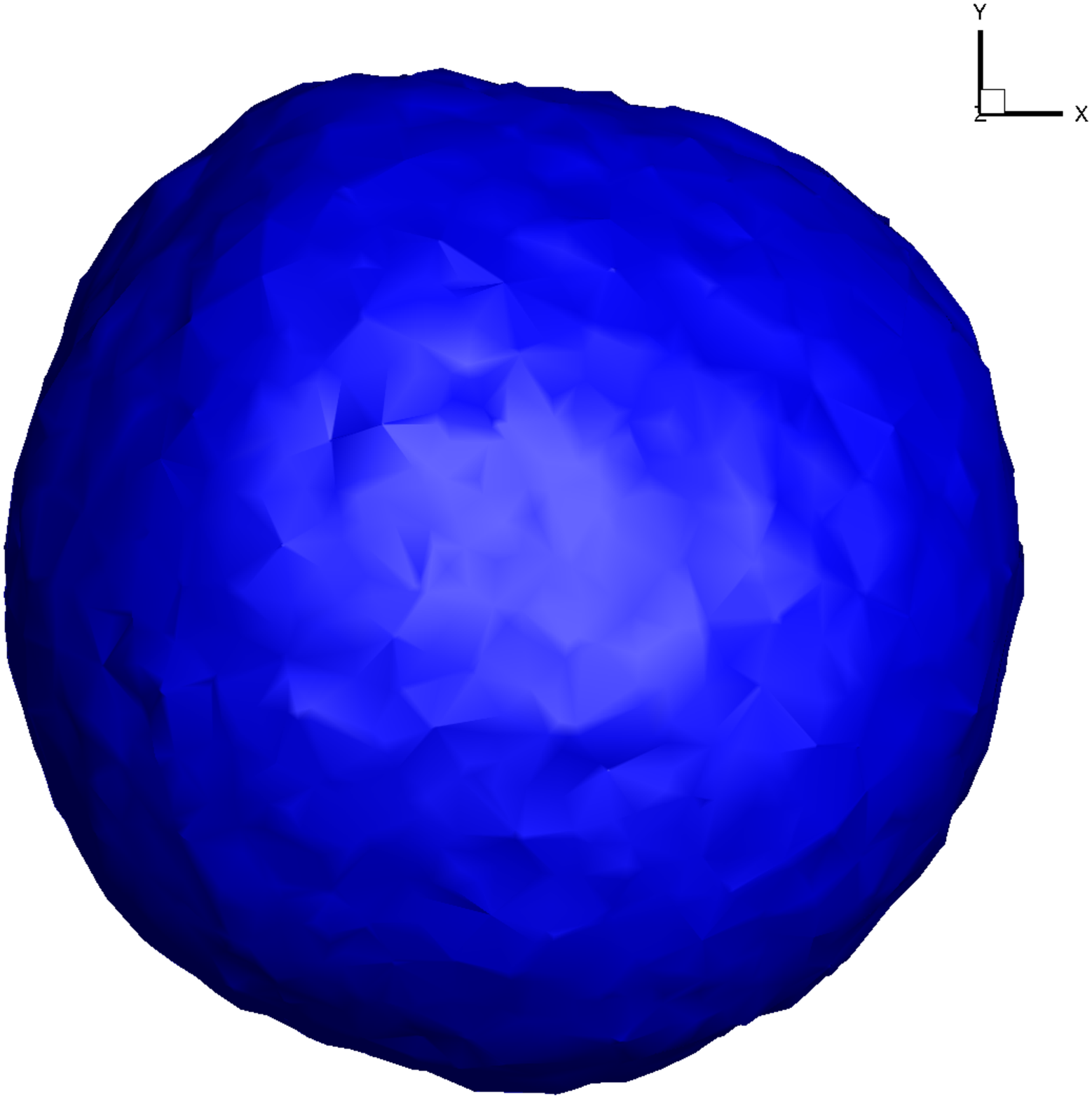}}
	\caption{\color{black}Numerical results of THINC-scaling scheme for 3D shear flow on unstructured tetrahedral mesh. Shown are the VOF 0.5-isosurfaces at $t=T/2$ and $t=T$.}
\label{3Dshear-unstruct-TSc}
\end{figure}

%~~~~~~~~~~~~~~~~~~~~~~~~~~~~~~~~~~~~~~~~~~~~~~~~~~~~~~~~~~~~~~~~~~~~~~~~~~~~~~~~~~~~~~~~~~~~~~~~~~~~~~~~~~~~~~~~~~~
We also show the numerical results of THINC-scaling scheme on tetrahedral grids of different resolutions in Fig.\ref{3Dshear-unstruct-TSc}. Here, the interface is represented using VOF 0.5-isosurface  at both $t=T/2$ and $t=T$. As compared to other methods described in \cite{jofre2014,Henning2019VOF}, the numerical results of our proposed scheme look more superior regarding the geometrical fidelity. As shown in Table \ref{3Dshear-uns-comparison}
 for quantitative evaluation, the numerical errors of the present scheme are smaller than other methods available for comparison.  

\begin{table}[ht]
\centering
\caption{Numerical errors $E\left(L_1\right)$ and convergence rates for 3D shear flow test on unstructured tetrahedral grid.} 
\begin{tabular}[t]{lSSS} \toprule
{$\textbf{Methods}$}          & {$32$}                & {Order} & {$64$} \\ \midrule
{THINC-scaling}               & {\color{black}$3.38\times10^{-3}$} & {1.67}  & {$1.05\times10^{-3}$} \\
%{THINC/LS}                    & {$3.37\times10^{-3}$} & {1.70}  & {$1.04\times10^{-3}$} \\
{THINC/QQ}                    & {$4.25\times10^{-3}$} & {1.65}  & {$1.35\times10^{-3}$} \\
{isoAdvector-plicRDF \cite{Henning2019VOF}} & {$8.43\times10^{-3}$} & {1.55}  & {$2.88\times10^{-3}$} \\
{Youngs \cite{jofre2014}} & {$6.15\times10^{-3}$} & {1.60}  & {$2.03\times10^{-3}$} \\
{LVIRA \cite{jofre2014}} & {$5.97\times10^{-3}$} & {1.87}  & {$1.64\times10^{-3}$} \\\bottomrule
\label{3Dshear-uns-comparison}
\end{tabular}
\end{table}

\section{Conclusions}
We propose a novel scheme to unify the solution procedures of level set and VOF methods that are two interface-capturing methods based on completely different concepts and numerical methodologies. The underlying idea of the scheme, THINC-scaling scheme, is to use the THINC function to scale/convert between the level set function and a continuous Heaviside function which mimics the VOF field. 

The high-quality THINC function is constructed by using: (1) high-order polynomials computed from the level set field to accurately retrieve the geometrical information of the interface, and (2) the constraint of VOF value to guarantee rigorous numerical conservativeness. Being a function handleable by conventional calculus tools, the THINC function facilitates efficient and accurate computations for interface-capturing, which takes the advantages from both VOF and level set methods. 

Being an interface-capturing method of great practical significance, THINC-scaling scheme doesn't involve explicit geometrical reconstruction and is algorithmically simple, which allows representing interface with high-order polynomials and implementing on unstructured grids straightforwardly without substantial difficulty. Even without geometrical reconstruction, an interface can be retrieved and well defined by the PSI (Polynomial Surface of the Interface) equation in THINC-scaling scheme. Using high-order polynomial to represent the moving interface, THINC-scaling scheme is able to resolve sub-grid structures.    

We verified the THINC-scaling scheme with widely used benchmark tests for moving interfaces on both structured and unstructured grids in comparison with other existing methods, which demonstrate the superior solution quality and the great potential of the proposed scheme as a moving interface-capturing scheme for practical utility. {\color{black}  Some efforts to make it available for applications, such as parallelizing the code and merging it to fluid solvers, are in progress. }

\section*{Acknowledgment}

This work was supported in part by the fund from JSPS (Japan Society for the Promotion of Science) under Grant Nos. 18H01366 and 19H05613. RA was supported in part by SNF project 200020\_175784.

\appendix
\setcounter{table}{0}
\setcounter{figure}{0}
\setcounter{equation}{0}
\renewcommand{\thetable}{\Alph{section}1}
\renewcommand{\thefigure}{\Alph{section}1}
%\appendix

\section{Numerical quadrature of THINC function}
\numberwithin{equation}{section}
\makeatletter 

The integration in space \eqref{mass_constraint} can be approximated using Gaussian quadrature as follows. 
\begin{equation}
\sum_{g=1}^G \omega_{g}\frac{1}{2}\left(1+\tanh\left(\beta\left(\mathcal{P}_{i}\left({\textbf{x}_{g}}\right)+\phi_i^{\Delta}\right)\right)\right)={\bar H}_{i}^{n},
\label{THINCf-integration}    
\end{equation}
where $\textbf{x}_{g}$ denotes a Gaussian quadrature point, and $\omega_g$ the corresponding weight satisfying $\sum_{g=1}^G \omega_{g}=1$. 
For the quadratic surface reconstruction in the present work, we follow \cite{xie2017toward} and use 6 and 9  Gaussian points for  triangular and rectangular elements respectively in 2D, while 11 and 27 points are used for tetrahedral and cubic elements in 3D. 

We recast \eqref{THINCf-integration} into 
\begin{equation}
\sum_{g=1}^G \omega_g \frac{\tanh\left(\beta\mathcal{P}_i\left(\textbf{x}_g\right)\right)+\tanh\left(\beta\phi_i^\Delta\right)}{1+\tanh\left(\beta\mathcal{P}_i\left(\textbf{x}_g\right)\right)\cdot\tanh\left(\beta\phi_i^\Delta\right)}=2\left({\bar H}_{i}^{n}-\frac{1}{2}\right), 
\label{THINCf-int2}    
\end{equation}
which is further simplified as 
\begin{equation}
f\left(D\right):=\sum_{g=1}^G\omega_g\frac{A_g+D}{1+A_g D}-C=0,
\label{THINCf-int3}    
\end{equation}
with
\begin{equation}
A_g=\tanh\left(\beta \mathcal{P}_i\left(\textbf{x}_g\right)\right), \quad D=\tanh\left(\beta\phi_i^{\Delta}\right), \quad\text{and}\quad 
C=2\left({\bar H}_{i}^{n}-\frac{1}{2}\right).
\end{equation}
{We note that $A_g$ and $C$ belong to $[-1,1]$, and that we look for a solution $D$ in $[-1,1]$.}

Given VOF value ${\bar H}_{i}^{n}$ and the surface polynomial $\mathcal{P}_i\left(\textbf{x}\right)$ computed from the level set field, the only unknown  $D$ in the non-linear algebraic equation \eqref{THINCf-int3} is solved by Newton-Raphson iterative method in the following form
\begin{equation}
D_{k+1}=D_{k}-\frac{f\left(D_{k}\right)}{f^{'}\left(D_{k}\right)},
\label{THINCf-NI}    
\end{equation}
with $D_{k}$ being the approximation solution at the $k$th step of iteration. 

{\color{black} This appendix is organised as follows: we first study under which condition the Newton-Raphson algorithm for an equation of type \eqref{THINCf-int3} converges, and then we show a simple modification of \eqref{THINCf-integration}, still of the form \eqref{THINCf-int3}, that satisfies these conditions. Using these conditions, we show that the Newton-Raphson algorithm converge to the unique solution, and that this solution (as well as all the terms of the sequence) stays in $[-1,1]$. These conditions are summarized in the following lemma \ref{THINC:Remi}}.
\begin{lemma}\label{THINC:Remi}
If we define $A_g=\tanh(\beta P_g+\gamma )$ and $D=\tanh(\beta\phi_i^{\Delta} -\gamma)$ such that either the condition 
\begin{subequations}
\begin{equation} \label{cond1}
\gamma< \min\limits_g (-\beta P_g)=-\max_g \big (\beta P_g)
\end{equation}
 or 
\begin{equation} \label{cond2}
 \gamma> \max\limits_g (-\beta P_g)=-\min_g \big (\beta P_g)
\end{equation}
\end{subequations}
holds true, then the Newton-Raphson method \eqref{THINCf-NI} with $f$ defined by \eqref{THINCf-int3}  converges to the unique solution of \eqref{THINCf-int2} {\color{black}with the initial solution chosen as $D_0=1$ for \eqref{cond1} and  $D_0=-1$} for \eqref{cond2}. In addition, all the terms of the sequence are in $]-1,1[$,  the sequence is monotonically decreasing (resp. increasing) in the case \eqref{cond1} (resp. \eqref{cond2}). Last the sequence converges quadratically.
\end{lemma}
\subsection{Discussion of \eqref{THINCf-NI} for \eqref{THINCf-int3}}
We define 
\begin{subequations}\label{THINC-Remi:1}
\begin{equation}
\label{THINC-Remi:1:1}
\varphi(D):=\sum_g \omega_g \dfrac{A_g+D}{1+A_gD}-C
\end{equation}
and assume that 
\begin{equation}
\label{THINC-Remi:2}A_g\in [-1,1], \qquad D\in [-1,1].
\end{equation}
\end{subequations}
Our goal is to find the sufficient conditions on $A_g$ such that the solution of $\varphi(D)=0$ is unique in $[-1,1]$ and that the Newton-Raphson algorithm converges to this unique solution. We will proceed as follows: first we will show that under the condition of $A_g\in [-1,1]$, there is a unique solution $D^\star$ in $[-1,1]$. Then we will study the Newton-Raphson algorithm, and give a condition on $A_g$ such that $D_k$ are always in $[-1,1]$ and converges quadraticaly to $D^\star$. This amounts to studying the concavity/convexity of $\varphi$ in $[-1,1]$.

\bigskip
In order to have a solution $D^\star$, a sufficient condition is that 
$$\varphi(-1) \varphi(1)\leq 0.$$
Since $C\in [-1,1]$, this condition is always met:  we have
$$\varphi(-1) =\sum_g \omega_g \dfrac{A_g-1}{1-A_g}-C=-1-C, \varphi(1) =\sum_g \omega_g \dfrac{A_g+1}{1+A_g}-C=1-C$$
and then $\varphi(-1) \varphi(1)=-(1+C)(1-C)=-(1-C^2)\leq 0$.

The second step is about the derivative of $\varphi$ with respect to $D$. We have
$$\varphi'(D)=\sum_g \omega_g \bigg ( \dfrac{1}{1+A_g D}-\dfrac{A_g(A_g+D)}{(1+A_gD)^2}\bigg )=\sum_g \omega_g \dfrac{1-A_g^2}{(1+A_gD)^2}\geq 0$$
so that $\varphi$ is an increasing function, again because $A_g\in [-1,1]$. It shows that the solution $D^\star$ of $\varphi(D^\star)=0$ is \emph{unique} in $[-1,1]$.
Note that $\varphi$ is strictly monotone if at least one $A_g$ is not equal to $\pm 1$, which is assumed next. 

Let's take the second derivative of $\varphi$:
$$\varphi''(D)=-2\sum_g \omega_g \dfrac{A_g(1-A_g^2)}{(1+A_gD)^3}.$$
Since we look for $D\in [-1,1]$  and  since $A_g\in [-1,1]$, we have $1+A_gD\geq 0$ and then:
\begin{itemize}
\item if $A_g\geq 0$  for all $g$, $\varphi''(D)\leq 0$: $\varphi$ is concave;
\item if $A_g\leq 0$  for all $g$, $\varphi''(D)\geq 0$: $\varphi$ is convex;
\item if the $A_g$ are of both side, we cannot decide like this.
\end{itemize}
{From now on, we assume that  the $A_g$ are all strictly  positive or negative.} This is our additional condition to guarantee that the Newton-Raphson algorithm converges.  Later we show how to achieve this for problem \eqref{THINCf-int3}-\eqref{THINCf-NI}.

%First, under this assumption, unless all the $A_g=0$ there exists a unique solution, $x^\star\in [-1,1]$ such that $\varphi(x^\star)=0$.

The Newton-Raphson algorithm is 
$$\text{ Given} \ D_0, \ D_{k+1}=D_k-\dfrac{\varphi(D_k)}{\varphi'(D_k)}.$$
We can assume that $D_k$ is never equal to $D^\star$, then
\begin{itemize}
\item In case of $A_g<0$ for all $g$, i.e. $\varphi$ is strictly convex:
\begin{itemize}
\item If $D_k\geq D^\star$, there exists $\xi_k\in ]D^\star, D_k[$ such that 
$$0\leq \varphi'(\xi_k)\leq \varphi'(D_k)$$ because $\varphi'$ is monotonicaly increasing. Recall that $\varphi$ is increasing (so $\varphi'\geq 0$) as well, then,
$$D_{k+1}-D^\star= (D_k-D^\star)\bigg ( 1-\dfrac{\varphi'(\xi_k)}{\varphi'(D_k)}\bigg )\geq 0,$$ 
hence $D_{k+1}\geq D^\star$.
In addition, if we assume that $D_k\leq 1$, we see that
$$D_{k+1}-1=D_k-1-\dfrac{\varphi(D_k)}{\varphi'(D_k)}\leq D_k-1\leq 0$$ because
 $\varphi$ is monotonicaly increasing, $0=\varphi(D^\star)\leq \varphi(D_k)$, and since $\varphi'(D_k)\geq 0$ the ratio
$$\dfrac{\varphi(D_k)}{\varphi'(D_k)}$$ is positive. This  also shows that the sequence $\{D_k\}$ is monotone decreasing, and we always have $D^\star\leq D_k\leq 1$ for all $k$.
\item If  $D_k\leq D^\star$, the situation is less favourable. First, using the same idea, we easily see that $D^\star\leq D_{k+1}$, but in this case, looking at the sign of $D_{k+1}-1$, we can see that we can find situations where $D_{k+1}>1$, and the whole reasoning falls apart.
\end{itemize}
The discussion above suggests that we should always consider the case of $D_k\geq D^\star$ when $A_g<0$ holds for all $g$. It can be obtained by setting the initial guess as $D_0=1$.

\item Assume that $\varphi$ is concave, we have the symmetric situations:
\begin{itemize}
\item if $D_k\leq D^\star$ then $D_{k+1}\leq D^\star$, with a monotone increasing sequence. In this case, we initialise with $D_0=-1$.
\item if $D_k\geq D^\star$ then $D_{k+1}\leq D^\star$ and as before the situation is much less favorable.
\end{itemize}
\end{itemize}

{\color{black}
Furthermore, we can prove that the convergence is quadratic as well known for the Newton-Raphson method. From Taylor expansion, we have 
\begin{equation}
D_{k+1}-D^\star=\rho \left(D_{k}-D^\star \right)^2,
\label{nt-3}    
\end{equation}
where the amplification factor reads
\begin{equation}
\rho=\frac{\varphi^{\prime \prime}\left(\xi_k\right)}{2 \varphi^{\prime}\left(D_{k}\right)}=-\frac{(A_g D_{k}+1)^{2} A_g}{(A_g \xi_k +1)^{3}}, 
\label{nt-4}    
\end{equation}
where $\xi_k$ is between $D_{k}$ and $D^\star $.  
We see that 
\begin{equation}
|D_{k+1}-D^\star|\leq |\rho| \left(D_{k}-D^\star \right)^2,
\label{nt-5}    
\end{equation}
and will show that $|\rho| \leq 1$ holds under the conditions suggested above. 
\begin{itemize}
\item In case of $A_g < 0$: as the iteration solutions are on the right side of $D^\star$, we have $D^\star\leq \xi_{k} \leq D_{k}$, and thus $A_g D_{k}\leq A_g \xi_{k} $. \eqref{nt-4}  implies  $|\rho| \leq 1$. 
\item In case of $A_g > 0$: as the iteration solutions are on the left side of $D^\star$, we have $D_{k}\leq \xi_{k} \leq D^\star$, and thus $A_g D_{k}\leq A_g \xi_{k} $. Again, \eqref{nt-4}  implies  $|\rho| \leq 1$. 
\end{itemize}

%[Xiao comment: it can be easily proved under the suggested condition.]

}
\subsection{Application to  \eqref{THINCf-NI} for \eqref{THINCf-int3}}
In \eqref{THINCf-int3}, we can also choose $A_g=\tanh(\beta P_g+\gamma )$ and  $D=\tanh(\beta\phi_i^{\Delta} -\gamma)$ for a well chosen $\gamma$ that makes all $A_g$ either positive or negative.

To make them all negative, we just need to have $\beta P_g+\gamma\leq 0$, so 
\begin{equation}
\label{THINCf-REMI-solution1}
\gamma< \min\limits_g (-\beta P_g)=-\max_g \big (\beta P_g).
\end{equation}
Then a possible $D_0$ is 
\begin{equation}
\label{THINCf-REMI-solution1:1}
D_0=1.
\end{equation}

Another solution is to choose $\beta P_g+\gamma\geq 0$, i.e.
\begin{equation}
\label{THINCf-REMI-solution3} 
\gamma> \max\limits_g (-\beta P_g)=-\min_g \big (\beta P_g).
\end{equation}
Then a possible $D_0$ is 
\begin{equation}
\label{THINCf-REMI-solution3:1}
D_0=-1.
\end{equation}

To choose between the two possible cases, it might be safe to avoid situations where $\gamma=0$.
This ends the proof and implementation of lemma \ref{THINC:Remi}.

{\color{black}
In the present work, we simply implement \eqref{THINCf-REMI-solution3} as 
\begin{equation}\label{cd2-1}
\gamma = -\min_g \big (\beta P_g)+\epsilon
\end{equation}
with $\epsilon=10^{-8}$, and set $D_0=-1$ following \eqref{THINCf-REMI-solution3:1}.

Once having $D$ solved, we finally get 
\begin{equation}
\phi_i^{\Delta}=\frac{1}{\beta}\left(\tanh^{-1}\left(D\right)+\gamma\right).
\end{equation}
}

{\color{black}

\section{Comparison between THINC-scaling method and THINC/LS method\cite{qian2018}}
\numberwithin{equation}{section}
\makeatletter 

This appendix provides a detail comparison between the present THINC-scaling method and the THINC/LS method\cite{qian2018}. They are different in both concept and solution algorithm.

THINC/LS follow the concept of the conventional CLSVOF method, i.e. VOF and LS are treated as two independent fields, which are updated separately and modified through coupling information from each other. In contrast to this, THINC-scaling method sees the VOF and LS as the two faces of the THINC reconstruction function. We can retrieve one from another by scaling or inverse-scaling with the THINC function. In this sense, the VOF and LS have no substantial difference. In practice, if an accurate and reliable THINC reconstruction can be maintained, it provides fidelity solutions for both VOF and level set fields. 
\begin{figure}[h]
\centering
\includegraphics[width=0.95\textwidth]{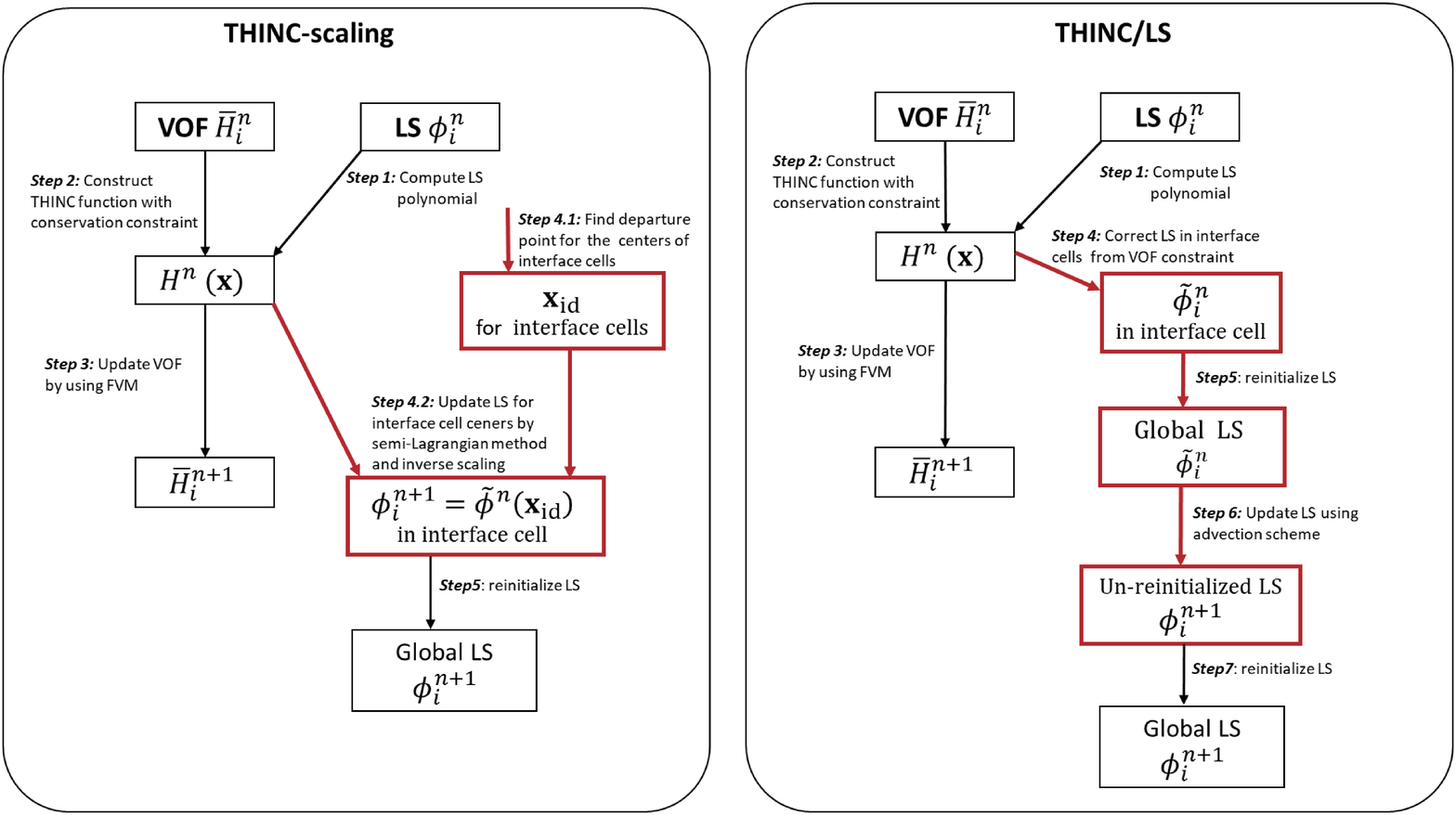}\hspace{0.3cm}
\caption{\small{Solution procedures of THINC-scaling method (left) and THINC/LS method (right). }}
\label{flowchart}
\end{figure}

Regarding solution procedure,  we  show the flowcharts of the solution procedure for the two methods in Fig.\ref{flowchart}.  THINC/LS needs to advect LS for whole domain, and the numerical scheme for advection affects the solutions. Moreover, one has to conduct a reinitialization  before the advection, and  another  reinitialization  step is also required after the advection. The THINC scaling provides a LS function for the interface cells without any interpolation, which facilitates the semi-Lagrangian step to get the LS values at the centers of the interface cells on new time level. Note that no reinitialization nor conventional Eulerian advection scheme is needed here. It's worth noting that the LS value is available everywhere in the interface cells by scaling the THINC reconstruction function to a LS field, thus no any interpolation is required to find the value at the departure point in the semi-Lagrangian computation. 

The THINC/LS scheme presented in \cite{qian2018} is limited to Cartesian grid. It needs a reliable advection scheme to transport the level set field when implemented on unstructured grids.  Whereas, the THINC-scaling scheme  does not require the advection computation for the global LS. The algorithmic simplicity of THINC scaling method eases its implementation  on unstructured grids.  

}
\clearpage{}
\bibliographystyle{unsrt}%elsarticle-num-names}
\bibliography{THINC_LS} 

\end{document}